%% file: Main_02.tex
\documentclass[twocolumn]{autart}
\input{Packages}    % all of the packages used in the document
\input{CustomCommands} 							% your macros

\begin{document}

% =========================================================
\input{Frontmatter}
\section{Introduction} \label{sec:introduction}
\input{Introduction}
	
% =========================================================
\section{Stochastic Dissipativity Theory} \label{sec:dissipativity}
\input{StochasticDissipativity}
	
% =========================================================
\section{Identifying Stochastic Dissipativity for Probabilistic Input Delays} \label{sec:general}
\input{StochasticKYP}

% =========================================================
\section{Numerical Experiments} %\label{sec:experiment}
\input{Numerical}

% =========================================================
\bibliographystyle{plain}
\bibliography{IJRNCreferences}%

% =========================================================
%\input{Backmatter}

% =========================================================
\end{document}

%% file: Packages.tex
\usepackage{multicol} 					% for tables with multiple columns
\usepackage{color} 							% for coloured text
\usepackage{pgfplots} 					% for graphics in latex
\usepackage{tikz} 							% for graphics in latex
\usepackage{amsmath,bm,amsfonts} % for math symbols and notation
\usepackage{graphicx,epstopdf}	% for figures
\usepackage{mathtools} % for equal by definition
\usepackage{pdfpages}						% include pages from another pdf. \includepdf[pages={-}]{myfile.pdf}
\usepackage{hyperref} 					% hyperlinks in text
\usepackage{subfigure}					% for subfigures

\usepackage{float}
\usepackage{array}
\usepackage{nicefrac}
\newcolumntype{L}[1]{>{\raggedright\let\newline\\\arraybackslash\hspace{0pt}}m{#1}}
\usepackage{centernot}

\usepackage{acronym}% heard about a great LaTeX package and thought I'd share. It's the acronym package. It lets you define a macro for each acronym, e.g.\acrodef{uav}[{\sc uav}]{unmanned aerial vehicle}. Then, when you want to use the acronym, you type \ac{uav} for singular and \acp{uav} for plural. The package will figure out when you first use the acronym and introduce it properly.

\usepackage{verbatim}

\usepackage{cuted}

%% file: CustomCommands.tex
\newcommand{\red}[1]{\textcolor{red}{#1}} % red writing
\newcommand{\blue}[1]{\textcolor{blue}{#1}} % red writing
\newcommand{\violet}[1]{\textcolor{violet}{#1}} % red writing

\newtheorem{define}{Definition}
\newtheorem{theorem}{Theorem}
\newtheorem{lemma}{Lemma}
\makeatletter
\newcommand\Autoref[1]{\@first@ref#1,@}
\def\@throw@dot#1.#2@{#1}% discard everything after the dot
\def\@set@refname#1{%    % set \@refname to autoefname+s using \getrefbykeydefault
	\edef\@tmp{\getrefbykeydefault{#1}{anchor}{}}%
	\xdef\@tmp{\expandafter\@throw@dot\@tmp.@}%
	\ltx@IfUndefined{\@tmp autorefnameplural}%
	{\def\@refname{\@nameuse{\@tmp autorefname}s}}%
	{\def\@refname{\@nameuse{\@tmp autorefnameplural}}}%
}
\def\@first@ref#1,#2{%
	\ifx#2@\autoref{#1}\let\@nextref\@gobble% only one ref, revert to normal \autoref
	\else%
	\@set@refname{#1}%  set \@refname to autoref name
	\@refname~\ref{#1}% add autoefname and first reference
	\let\@nextref\@next@ref% push processing to \@next@ref
	\fi%
	\@nextref#2%
}
\def\@next@ref#1,#2{%
	\ifx#2@ and~\ref{#1}\let\@nextref\@gobble% at end: print and+\ref and stop
	\else, \ref{#1}% print  ,+\ref and continue
	\fi%
	\@nextref#2%
}

\hypersetup{nolinks=true}

\newcommand{\bb}[1]{\mathbb{#1}} % double bar R
\newcommand{\bmat}[1]{\begin{bmatrix} #1 \end{bmatrix}} % bracketed matrix
\newcommand{\pmat}[1]{\begin{pmatrix} #1 \end{pmatrix}} % parenthesis matrix
\newcommand{\w}{\mathbf{w}} % bold w
\newcommand{\He}{\mbox{He}}

\newcommand{\Expect}[1]{\mathbb{E}\Big[ #1 \Big]}
\newcommand{\expect}[1]{\mathbb{E}[ #1 ]}

\newcommand{\trunc}{{[T]}}

\newcommand{\Uspace}{\mathcal{U}}
\newcommand{\Space}{\mathcal{X}}
\newcommand{\Fcal}{\mathcal{F}}
\newcommand{\bbP}{\mathbb{P}}

\acrodef{iid}[i.i.d.]{independent, identically distributed}

\newcommand{\defeq}{\vcentcolon=}

%% file: Frontmatter.tex
\begin{frontmatter}
\runtitle{Dissipativity for probabilistic input delays}  % Running title for regular 

\title{Stochastic Dissipativity for Systems with Probabilistic Input Delays} % Title, preferably not more 

\thanks[footnoteinfo]{This work is supported by NSF GRFP Grant No. 1644868, the Alfred P. Sloan Foundation, and ONR Grant No. N00014-23-1-2043. Corresponding author Ethan J. LoCicero.}

\author{Ethan J. LoCicero, }\ead{ethan.locicero@duke.edu}
\author{Amy K. Strong, }\ead{amy.k.strong@duke.edu}
\author{Leila J. Bridgeman}\ead{leila.bridgeman@duke.edu}

\address{Thomas Lord Department of Mechanical Engineering and Materials Science, Duke University, North Carolina, United States of America}

\begin{keyword}                           % Five to ten keywords,  
time-delays, stochastic control, robust control of nonlinear systems, linear matrix inequalities, passivity             % chosen from the IFAC 
\end{keyword}                             % keyword list or with the 
                                          % help of the Automatica 
                                          % keyword wizard

\begin{abstract}
This work considers stochastic operators in general inner-product spaces, and in particular, systems with stochastically time-varying input delays of a known probability distribution. Stochastic dissipativity and stability are defined from an operator-theoretic perspective, and the well-known open-loop dissipativity conditions for closed-loop/network stability are extended to the stochastic case. Criteria are derived to identify dissipative nonlinear systems with stochastic input delays, and this result is used to find delay-distribution-dependent linear matrix inequality conditions for stochastic dissipativity of a linear system with input delays of a known probability distribution. A numerical experiment demonstrates the utility of the resulting criteria for robust plant analysis and controller design, highlighting significantly reduced conservatism compared to deterministic methods.
\end{abstract}

\end{frontmatter}

%% file: Introduction.tex
Signal delays are pervasive in a wide variety of fields, including chemical systems, lasers, aircraft, the internet, robotic tele-operation, traffic control, biology, and medicine \cite{Anderson1988,Hale1993,Kolmanovskii1999,Seuret2016}. Delays are particularly notable in large-scale, distributed systems where control and measurement signals are communicated over vast distances and unreliable networks. Such delays can adversely affect performance and stability, and they greatly complicate the mathematical analysis of otherwise simple problems \cite{Fridman2014,Papachristodoulou2005,Seuret2016}. 
	
Stability analysis of delayed systems is a challenging problem. For instance, the necessary and sufficient form of a Lyapunov-Krasovskii (L-K) functional for proving stability of a linear time invariant (LTI) system with constant state delays results in infinite-dimensional linear matrix inequality (LMI) conditions \cite{Papachristodoulou2005,Seuret2016}. Moreover, no L-K functional form is yet known to be necessary and sufficient for proving the stability of an LTI system with time-varying delays. As a result, there is a vast literature on deriving variations of L-K functionals \cite{Fridman2009,Gao2007,Guo2022,Kwon2013,Papachristodoulou2005,Seuret2016}. 

An alternative approach to stability analysis has been to use input-output methods. Input-output stability notions, including passivity, gain, and conic sectors \cite{Zames1966} have been used to analyze the stability of delayed systems \cite{Cheng2012,Fridman2002,Fridman2006,Li2002,Mahmoud2009,Mahmoud2013,Zhang2022}. However, existing results for time-varying delays are often limited and overly-conservative. Moreover, these results often don't account for the inherent stochasticity in many delay processes. However, input-output methods are often thought of as special cases of QSR-dissipativity as introduced in Reference~\cite{Willems1972}, so the recent (and long foreshadowed \cite{Safonov1980}) development of stochastic dissipativity \cite{Haddad2022} provides an opportunity for less conservative analysis of delay systems. Consequently, this paper develops dissipativity theory for systems with stochastic input delays.
	
In the framework of Reference~\cite{Haddad2022}, a stochastic system need only satisfy a dissipativity property \textit{in expectation}, which is then used to show Lyapunov stability \textit{in probability}. This has the potential to significantly reduce the conservatism of stability results for stochastic systems. While the stochastic dissipativity definition and stability theorem from Reference~\cite{Haddad2022} are broad, their criteria to identify dissipativity, which comes in the form of a KYP-style lemma, is restricted to nonlinear systems with affine dependence on the random process. This restriction precludes delays from the analysis.
	
In this paper, the notion of stochastic dissipativity is revisited with two major changes compared to Reference~\cite{Haddad2022}. First, dissipativity and stability are approached from an operator-theoretic perspective. In many cases, operator-theoretic perspectives, which are centered on input-output stability, are equivalent to state-space perspectives, which are centered on Lyapunov stability. However, in certain cases, input-output stability implies Lyapunov stability \cite{Hill1980}. Second, a delay-distribution-dependent KYP-style lemma is derived to characterize stochastic, QSR-dissipative systems. Along with a new stochastic dissipativity definition, stability theorem, and general dissipativity criteria for nonlinear systems with stochastic input delays, a KYP-style lemma is provided for LTI systems with stochastic input delays of bounded probability distribution, which results in LMI conditions for stochastic dissipativity. The utility of these LMI conditions is illustrated in a numerical example. 

%==============================
\subsection{Preliminaries} \label{sec:preliminaries}
	
In this section, notation is established and some properties of inner product spaces and probability spaces are reviewed. Let element-wise inequalities be denoted $<$, $>$, $\leq$, and $\geq$, and let matrix definiteness be denoted $\succ$, $\prec$, $\succeq$, and $\preceq$. An overline, $\overline{(\cdot)}$, denotes a bound above in the sense of $a\leq \overline{a}$ for scalars and $A\preceq \overline{A}$ for matrices. The $n$-dimensional identity and $n\times m$-dimensional zero matrices are $I_n$ and $0_{n\times m}$, respectively. Subscripts are dropped when evident from context. Asterisks denote duplicate blocks in symmetric matrices, and $(\cdot)^T$ denotes the transpose of a matrix. To preserve space, $\He(M) := M+M^T$, and $\mathrm{Sq}(M,x) \defeq x^TMx$. The set of real numbers over a closed interval, $[a,b]$, is denoted $\mathbb{R}_{[a,b]}$. Square brackets are replaced with parenthesis to denote the open interval, as in $\mathbb{R}_{(a,b)}$. Integers are denoted likewise with $\mathbb{Z}$. The $n$-dimensional real vectors and $n\times m$-dimensional real matrices are denoted $\mathbb{R}^n$ and $\mathbb{R}^{n\times m}$, respectively. The set of $n\times n$ symmetric matrices is denoted $\mathbb{S}^{n \times n}$. The set of square-summable sequences is denoted $\ell_2$, and the $\ell_2$-norm is denoted $||\cdot||_{\ell_2}$.
		
A set $\Space$ is a real inner product space if it is closed under addition, $+:\Space\times \Space\rightarrow\Space$, and scalar multiplication, $\cdot:\mathbb{R}\times\Space\rightarrow\Space$, and is equipped with an inner product mapping $\langle(\cdot),(\cdot)\rangle:\Space\times\Space\rightarrow\mathbb{R}$. If the elements of $\Space$ are $n$-dimensional vector sequences, i.e. $x:\mathbb{Z}\rightarrow\mathbb{R}^n$ for all $x\in\Space$, then the space is denoted $\Space^n$. The space $\Space$ is \textit{complete} if $||x||^2_\Space:= \langle x,x\rangle <\infty$ for all $x\in\Space$. 
The extension of $\Space$, denoted $\Space_e$, satisfies $||x_\trunc||^2_\Space <\infty$ for all $T\in\mathbb{R}_{[0,\infty)}$, where $x_\trunc$ is the truncation of $x(k)$ at $k=T$, defined as $x(k) = x(k)$ for $k\leq T$ and $x(k) = 0$ for $k>T$. The truncated inner product is also sometimes denoted $\langle x,y\rangle_\trunc$. A complete probability space is defined by $(\Omega,\mathcal{F},\mathbb{P})$, where $\Omega$ is the sample space, $\mathcal{F}$ is the event space (a $\sigma$-algebra of subsets of $\Omega$), and $\mathbb{P}:\mathcal{F}\rightarrow [0,1]$ is a probability measure (a nonnegative countably additive set function satisfying $\mathbb{P}(\Omega) = 1$). A single sample from $\Omega$ at time $T$ is denoted $w_T$, while sequence of samples from time $0$ to time $T$ is denoted with brackets, as in $w_{[T]}\in\mathcal{F}$. Independent and identically distributed samples are denoted i.i.d. The notation $\mathbb{E}_{w_{[k]}\sim\mathbb{P}}(f(w_{[k]}))$ denotes the expected value of the function $f(w_{[k]})$, where $w_{[k]}\in\mathcal{F}$ is distributed according to $\mathbb{P}$.

\subsection{A Motivating Counterexample} \label{sec:counterexample}
	
It is well known that the gain of a system is invariant with respect to constant input delays. Therefore, when investigating the stability of a network using the Small Gain Theorem \cite{Desoer1975}, constant input delays can be neglected, begging the question: why is so much effort taken here to establish input-output properties? The answer is, for time-varying input delays, even establishing gain becomes complicated.

Consider, for instance, the linear memory-less system $y=Ku$, where $u\in\ell_{2}$ and $K\in\mathbb{R}$. The gain of this system is $K$. Now consider its delayed counterpart, $y_d = Ku(k-w_k)$, where $0\leq w_k\leq 2$. A valid trajectory of this system is given by the sequences $\{u(k)\} = \{2,1,0,0,...,0\}$ and $\{w_k\} = \{0,1,2,0,...,0\}$. (Note, this delay sequence represents a realistic case of communication dropout over two time steps). The input norm is $||u||^2_{\ell_2} = 5$, and consequently, $||y||^2_{\ell_2} = 5K$, as expected. However, $||y_d||^2_{\ell_2} = 12K$, which is a gain of $\frac{12}{5}K>K$. If the maximum delay increases, the gain can actually become arbitrarily large, by selecting $w_k = k$. Therefore, time-varying input delays can have a significant impact on gain and more general dissipativity properties.
	
There have been several recent investigations into gain bounds for systems with time-varying input delays \cite{Cheng2012,Zhang2022,Zheng2021}, as well as a broad literature on L-K functions for systems with time-varying delays \cite{Gao2007,Guo2022,Kwon2013}. However, these results presume no additional information about the delay beyond its maximum and minimum value. Here, stochastic information about the delay is used to recover stochastic dissipativity properties and ultimately less conservative stochastic stability guarantees.

%% file: StochasticDissipativity.tex
This section defines a stochastic version dissipativity and input-output stability on general inner-product spaces. Then, Theorem~\ref{thm:OLstable} shows when stochastic dissipativity implies stochastic input-output stability, and Theorem~\ref{thm:CLstable} builds on this to derive stability conditions for networks of dissipative operators. Both theorems are direct extensions of Vidyasagar's Dissipativity Theorem \cite{Vidyasagar1981} to the realm of stochastic dissipativity.
	
Dissipativity was introduced as an extension of Lyapunov theory to open systems (systems with inputs and outputs) \cite{Willems1972}. However, when viewed from an operator-theoretic perspective \cite{Vidyasagar1981}, it can be seen as a multi-dimensional generalization of conic sectors \cite{Zames1966}. The operator-theoretic perspective foregoes the concept of a state to show input-output stability of a more general operator on an inner product space. Though dissipativity was originally developed for continuous-time systems, it has long been extended to discrete systems \cite{Byrnes1994}. More recently, the notion of stochastic dissipativity, or dissipativity-in-expectation, has been introduced in the Lyapunov framework \cite{Haddad2022,Rajpurohit2017}. In Definition~\ref{def:QSR}, a similar definition that provides several advantages is introduced from the operator-theoretic perspective.

\begin{define} \label{def:QSR}
    Given a complete real inner product space, $\mathcal{X}$, and a complete probability space, $(\Omega,\mathcal{F},\mathbb{P})$, where any two samples, $w_i,w_j\in\Omega$, are mutually independent and identically distributed (i.i.d.), a system $G:\Space^n_e \times \mathcal{F} \rightarrow \Space^m_e$ is $(Q,S,R)$-dissipative-in-expectation (or stochastically dissipative) with $Q\in\mathbb{S}^{m}$, $S\in\mathbb{R}^{m\times n}$, and $R\in\mathbb{S}^n$ if there exists $\beta\in\mathbb{R}$ (possibly a function of initial conditions) such that for all $T\in\mathbb{Z}_{(0,\infty)}$ and $ u\in\Space^n_e$,
	\begin{comment}
		\begin{align}
			\Expect{\langle Gu_\trunc,QGu_\trunc\rangle + \langle Gu_\trunc,2Su_\trunc\rangle + \langle u_\trunc,Ru_\trunc\rangle} \geq \beta \; \forall \; T>0, \; \forall \, u\in\Space^n_e.
		\end{align}
	\end{comment}
	\begin{align}
	\underset{{w_{[T]}\sim\mathbb{P}}}{ \mathbb{E}} \Big(\Big\langle \xi,\bmat{Q & S \\ S^T & R}\xi\Big\rangle\Big) \geq \beta, \label{eqn:QSR}
	\end{align}
 where $\xi^T = \bmat{G^T(u_\trunc,w_\trunc) & \; \; u_\trunc^T}$.
\end{define}

This definition is similar to Reference \cite{Haddad2022}, Equation 30, except it is restricted to the QSR-type supply rate and relaxed to general inner product spaces. The inclusion of nonzero $\beta$ is also an important relaxation for handling initial conditions in the operator-theoretic framework. Next, a stochastic version of input-output stability is presented.

\begin{define} \label{def:IO}
    Given a complete real inner product space, $\mathcal{X}$, and a complete probability space, $(\Omega,\mathcal{F},\mathbb{P})$, where any two samples, $w_i,w_j\in\Omega$, are i.i.d., an operator $G:\Space_e^n\times \mathcal{F} \rightarrow\Space_e^m$ is input-output stable-in-expectation (or stochastically input-output stable) on $\Space$ if there exists $\beta\in\mathbb{R}$ (possibly depending on initial conditions) and $\gamma\in\mathbb{R}_{(0,\infty)}$ such that for all $u\in\Space^n$, $\mathbb{E}_{w_\trunc\sim\mathbb{P}}(||G(u_\trunc,w_\trunc)||^2_\Space) \leq \gamma ||u_\trunc||^2_\Space + \beta$ for all $T\in\mathbb{Z}_{(0,\infty)}$.
\end{define}

Naturally, one might ask whether this is a strong enough statement of stability to be practically useful. In the following lemma, Markov's inequality demonstrates that stochastic input-output stability does in fact ensure that the output gain is almost surely bounded.

\begin{lemma} \label{lem:almostSurely}
	If $G:\mathcal{X}_e^n\times \mathcal{F}\rightarrow \mathcal{X}_e^m$ is stochastically input-output stable on the complete, real inner product space, $\mathcal{X}$, and the complete probability space, $(\Omega,\mathcal{F},\mathbb{P})$, then $\|G(u_{[T]},w_{[T]})\|_\mathcal{X}^2<\infty$ almost surely.
\end{lemma}
\begin{pf}
	By Markov's Inequality \cite{Markov1884,Huber2019}, $\mathbb{P}(X{\geq} a) \leq \frac{\mathbb{E}(X)}{a}$ for any non-negative random variable ,$X$, with finite expectation, and for any $a>0$. Letting $X = \|G(u_{[T]},w_{[T]})\|_\mathcal{X}^2$ and substituting the upper bound on its expectation from Definition~\ref{def:IO}, there exists some $\gamma\in\mathbb{R}_>$ and $\beta\in\mathbb{R}$ such that $\mathbb{P}(\|G(u_{[T]},w_{[T]})\|_\mathcal{X}^2\geq a) \leq \frac{\gamma\|u_{[T]}\|_\mathcal{X}^2 + \beta}{a}$ for all $a>0$, $u\in\mathcal{X}^n$, $T\in\mathbb{Z}_>$. Taking the limit, $\lim_{a\rightarrow\infty}\mathbb{P}(\|G(u_{[T]},w_{[T]})\|_\mathcal{X}^2\geq a) = 0$. Therefore, $\|G(u_{[T]},w_{[T]})\|_\mathcal{X}^2$ is unbounded almost never.
\end{pf}

From Lemma~\ref{lem:almostSurely}, this definition of stochastic input-output stability is about as strong as deterministic input-output stability. While the stochastic version allows for some events to cause the output norm to be unbounded, these events happen with zero probability. Now, Theorem~\ref{thm:OLstable} relates stochastic dissipativity to stochastic input-output stability.

\begin{theorem} \label{thm:OLstable}
    Consider a complete real inner product space, $\mathcal{X}$, and a complete probability space, $(\Omega,\mathcal{F},\mathbb{P})$, where any two samples, $w_i,w_j\in\Omega$, are i.i.d. If $G:\Space^n_e \times \mathcal{F} \rightarrow \Space^m_e$ is stochastically dissipative with $Q\prec 0$, then $G$ is stochastically input-output stable on $\Space$.
\end{theorem}

\begin{pf}
	The proof employs the S-procedure \cite{Derinkuyu2006}, which is to say $M\preceq 0$ is implied by the existence of some $\alpha\in\mathbb{R}_{(0,\infty)}$ such that $M+\alpha F\preceq 0$ for some $F\succeq 0$. In this case, the aim is to show there exists $\beta_{io}\in\mathbb{R}$ and $\gamma\in\mathbb{R}_{(0,\infty)}$ such that $M = \mathbb{E}_{w_\trunc\sim\mathbb{P}}(||G(u_\trunc,w_\trunc)||^2_\Space) - \gamma||u_\trunc||^2_\Space - \beta_{io} \leq 0$ for all $T\in\mathbb{R}_{(0,\infty)}$. To do so, let $F$ be given by \autoref{eqn:QSR}. Rearranging and letting $\beta_{io} = -\alpha\beta$ yields the condition
	\begin{align*}
		\underset{w_\trunc\sim\mathbb{P}}{\mathbb{E}} \Big(\Big\langle\xi,\bmat{\alpha Q + I & \alpha S \\ \alpha S^T & \alpha R - \gamma I }\xi\Big\rangle\Big) \leq 0.
	\end{align*}
	This condition is implied by the matrix inside being negative definite, which by Schur complement is equivalent to $\alpha Q + I \prec 0$ and $\alpha R - \gamma I - \alpha S^T (\alpha Q + I)^{-1}\alpha S \prec 0$. There exists $\alpha$ large enough that the first inequality is satisfied if and only if $Q\prec 0$. Further, for any $S$, $R$, and compatible $\alpha$ and $Q$, there exists a finite $\gamma$ large enough to satisfy the second inequality.
\end{pf}

\begin{theorem} \label{thm:CLstable}
    Consider a complete real inner product space, $\mathcal{X}$, $N\in\mathbb{Z}_{(0,\infty)}$ complete probability spaces, $(\Omega_i,\mathcal{F}_i,\mathbb{P}_i)$, where any two samples, $w^i_j,w^i_k\in\Omega_i$, are i.i.d., and $N$ operators $G_i:\Space_e^{n_i}\times \mathcal{F}_i \rightarrow\Space_e^{m_i}$, where $G_i$ is $(Q_i,S_i,R_i)$-dissipative-on-average for $i\in1,\dots,N$. Let these operators be connected by $y^i{=}G_i(u^i,w^i)$ and $u {=} Hy {+} \eta$, where $u^T {=} [{u^i}^T\dots {u^N}^T]$, $y^T = [{y^1}^T\dots {y^N}^T]$, $\eta\in\Space_e^n$, and $H$ is a real static matrix of appropriate dimensions. The composite network $G:\Space_e^n \times \mathcal{F}_1 \times \dots \times \mathcal{F}_N \rightarrow\Space_e^m$ defined by $y = G(\eta,w^i,\dots,w^N)$ is $(Q,S,R)$-dissipative-on-average with $Q = Q_\Lambda + \He(S_\Lambda H) + H^TR_\Lambda H$, $S = S_\Lambda + H^TR_\Lambda$, and $R = R_\Lambda$, where $Q_\Lambda = \mbox{diag}(\lambda_1Q_1,\dots,\lambda_NQ_N)$, $S_\Lambda = \mbox{diag}(\lambda_1S_1,\dots,\lambda_NS_N)$, $R_\Lambda = \mbox{diag}(\lambda_1R_1,\dots,\lambda_NR_N)$, and $\lambda_1,\dots,\lambda_N\in\mathbb{R}_{(0,\infty)}$. 
\end{theorem}

\begin{pf}
	By definition, $(\lambda Q,\lambda S, \lambda R)$-dissipativity-on-average is identical for any $\lambda\in\mathbb{R}_{(0,\infty)}$. Therefore, from Definition~\ref{def:QSR}, it is known that for all $i\in1,\dots,N$, there exists $\beta_i\in\mathbb{R}$ such that for all $T>0$ and $u^i\in\Space^{n_i}_e$,
	\begin{align*}
		\underset{w_\trunc\sim\mathbb{P}}{\mathbb{E}}&\Big(\langle G_i(u^i,w^i),\lambda_iQ_iG_i(u^i,w^i)\rangle_\trunc \\ & \hspace{-1mm} + \langle G_i(u^i,w^i),2\lambda_iS_iu^i\rangle_\trunc + \langle u^i,\lambda_iR_iu^i\rangle_\trunc\Big) \geq \beta_i. 
	\end{align*}
	For more compact notation, let $G = \mbox{diag}(G_1,\dots,G_N)$, $u^T = [{u^1}^T \dots {u^N}^T]$, and $\beta^T = [\beta_1 \dots \beta_N]$, then the $N$ known inequalities can be combined as
		\begin{align*}
		\underset{w_\trunc\sim\mathbb{P}}{\mathbb{E}}&\Big(\langle G(u,w^i,\dots,w^N),Q_\Lambda G(u,w^i,\dots,w^N)\rangle_\trunc \\ &\hspace{-2mm} + \langle G(u,w^i,\dots,w^N),2S_\Lambda u\rangle_\trunc + \langle u,R_\Lambda u\rangle_\trunc\Big) \geq \beta 
	\end{align*}
for all $T>0$ and $u \in \Space^n_e$. Substituting the interconnection equations and rearranging yields
\begin{align*}
	\underset{w_\trunc\sim\mathbb{P}}{\mathbb{E}}\Big(&\langle y,\big(Q_\Lambda + \He(S_\Lambda H) + H^TR_\Lambda H\big) y \rangle_\trunc \\ &+ \langle y,2\big( S_\Lambda H + R_\Lambda\big) \eta \rangle_\trunc + \langle \eta,R_\Lambda \eta \rangle_\trunc\Big) \geq \beta
\end{align*}
for all $T{>}0$ and $\eta {\in} \Space^n_e$.
The result follows by Definition~\ref{def:QSR}.
\end{pf}

Together, Theorems~\ref{thm:OLstable} and \ref{thm:CLstable} provide open-loop dissipativity conditions for closed-loop and network stability.

%% file: StochasticKYP.tex
\subsection{General Criteria}
	In this section, stochastic dissipativity criteria are derived for nonlinear, input-affine systems with stochastic input delays. This general lemma stops several steps short of the system of equations that result from the standard KYP-style lemma in Reference~\cite{Hill1976} and the stochastic variation in Reference~\cite{Haddad2022}. Those existing results provide systems of equations to solve for a Lyapunov-like function. For delay systems, Lyapunov functions are generally too restrictive to derive delay-dependent stability criteria. For LTI systems with constant state delays, the complete L-K functional is known to be necessary and sufficient for Lyapunov stability, though its form is usually too complex to deal with directly \cite{Seuret2016}. For nonlinear systems with time-varying state delays, no such criteria has yet been discovered. Consequently, the literature on L-K functionals is vast, and the best functional form is an open question. Therefore, this lemma is presented up to the point where a particular L-K functional must be chosen. In the next section, the continuation of this lemma is demonstrated with a particular choice of L-K functional for LTI systems with bounded input delays.
	
	\begin{theorem} \label{thm:genKYP}
        Consider a complete real inner product space, $\mathcal{X}$, and a complete probability space, $(\Omega,\mathcal{F},\mathbb{P})$, where $\Omega = \mathbb{Z}_{[0,\infty)}$ and any two samples, $w_i,w_j\in\Omega$, are i.i.d. Consider the operator $G:\Space_e^m \times \mathcal{F} \rightarrow \Space_e^p$, where $y = G(u,w)$ is given by
		\begin{align}
			G:\begin{cases}
			x(k+1) = f(x(k)) + g(x(k))u(k-w_k)\\
			y(k) = h(x(k)) + j(x(k))u(k-w_k), 
			\end{cases} \label{eqn:inputdelay}
		\end{align}
		where $x\in\Space^n_e$ is the state, $f{:}\mathbb{R}^n\rightarrow\mathbb{R}^n$, $g{:}\mathbb{R}^n\rightarrow\mathbb{R}^{n\times m}$, $h{:}\mathbb{R}^n\rightarrow\mathbb{R}^p$, and $j{:}\mathbb{R}^n\rightarrow\mathbb{R}^{p\times m}$. The system $G$ is $(Q,S,R)$-dissipative-in-expectation if there exists some function $V{:} \mathbb{Z} \times \Space^n \times \Space^m \times \mathbb{Z}_{[0,\infty)} \rightarrow \mathbb{R}$, satisfying for all $k\in\mathbb{Z}$, $x(0)\in\mathbb{R}^n$, and $u_{[k]}\in\mathcal{X}^m$,
		\begin{align}
			\underset{w_{[k+1]} \sim\mathbb{P}}{\mathbb{E}}\Big(&\mathrm{Sq}\big(Q, h(x(k)) + j(x(k))u(k-w_k)\big) \nonumber \\ 
            &\hspace{-8mm} + 2\left(h(x(k)) + j(x(k))u(k-w_k)\right)^TSu(k) \nonumber \\ &\hspace{-8mm}+ \mathrm{Sq}\big(R,u(k)\big) - \overline{\Delta V(k,x_{[k]},u_{[k]},w_{[k]})}\Big) \geq 0, \label{eqn:KYPgen}
		\end{align}
		where 		
		$V(0,x(0),u_{[0]},w(0))$ depends on initial conditions, $V(k,x_{[k]},u_{[k]},w_{[k]})\geq 0$, and $\Delta V(k,x_{[k]},u_{[k]},w_{[k]}) := V(k{+}1,x_{[k+1]},u_{[k+1]},w_{[k+1]})-V(k,x_{[k]},u_{[k]},w_{[k]})$ is the  difference of $V$ at $k$.
	\end{theorem}
	\begin{pf}
		Define the expected difference of $V$ at $k$ as 
        \begin{align*}
		\underset{\mathbb{E}}{\Delta} V(k,x_{[k]},u_{[k]},w_{[k]}) \defeq \hspace{-2mm}\underset{w_{[k{+}1]}\sim\mathbb{P}}{\mathbb{E}}\big(\Delta V(k,x_{[k]},u_{[k]},w_{[k]})\big).    
	\end{align*} 
        Using the properties of $V(k,x_{[k]},u_{[k]},w_{[k]})$ and a telescoping sum, we have for all $T>0$,
		\begin{align}
			0 &= \underset{w_{[T+1]}\sim\mathbb{P}}{\mathbb{E}}\Big(V(T{+}1,x_{[T+1]},u_{[T+1]},w_{[T+1]}) \nonumber \\ & \hspace{.7cm} - V(0,x_{[0]},u_{[0]},w_{[0]})\Big) - \sum_{k=0}^{T} \underset{\mathbb{E}}{\Delta} V(k,x_{[k]},u_{[k]},w_{[k]})  \nonumber \\ &\geq \alpha- \sum_{k=0}^{T} \overline{\underset{\mathbb{E}}{\Delta} V(k,x_{[k]},u_{[k]},w_{[k]})}
			\label{eqn:S}
		\end{align}
		where $\alpha \defeq -\mathbb{E}_{w_{[0]}\sim\mathbb{P}}(V(0,x_{[0]},u_{[0]},w_{[0]}))$ is a function of initial conditions. Applying the S-procedure to \autoref{eqn:QSR} and \autoref{eqn:S} shows that
		\begin{align*}
			&\underset{w_\trunc\sim\mathbb{P}}{\mathbb{E}}\Big(\sum_{k=0}^{T} \Big(y^T(k)Qy(k) + 2y^T(k)Su(k) \\ & \hspace{1cm} + u^T(k)Ru(k)\Big)\Big) \geq \beta \\
			\impliedby& \underset{w_\trunc\sim\mathbb{P}}{\mathbb{E}}\Big(\sum_{k=0}^{T} \Big(y^T(k)Qy(k) + 2y^T(k)Su(k) \\ & \hspace{1cm} + u^T(k)Ru(k)\Big)\Big) \\ & \hspace{1cm}+ \alpha- \sum_{k=0}^{T} \overline{\underset{\mathbb{E}}{\Delta} V(k,x_{[k]},u_{[k]},w_{[k]})} \geq \beta
		\end{align*}
		Substituting \autoref{eqn:inputdelay}, selecting $\beta = \alpha$, and rearranging yields
		\begin{align*}
			\underset{w_{[T+1]}\sim\mathbb{P}}{\mathbb{E}}&\Big(\sum_{k=0}^{T}\Big( \mathrm{Sq}\big(Q,h(x(k)) + j(x(k))u(k-w_k)\big) \\ & + 2(h(x(k)) + j(x(k))u(k-w_k))^TSu(k) \\ & + \mathrm{Sq}\big(R,u(k)\big) - \overline{\Delta V(k,x_{[k]},u_{[k]},w_{[k]})}\Big)\Big) \geq 0.
		\end{align*}
		By linearity of the expectation operator, this is  implied by \autoref{eqn:KYPgen}.
	\end{pf}

	\subsection{LMI Criteria} \label{sec:particular}
	
	This section restricts the scope of discussion to LTI systems with stochastic input delays of a known, bounded distribution. First, the following lemma will be useful.

    \begin{lemma} \label{lem:RecursiveSchur}
    For any real matrices $A$, $B$, and $C$ of appropriate dimensions, where $A$ and $C$ are symmetric,
    \begin{align}
        \bmat{A & B \\ B^T & C} \succeq 0 &\iff \bmat{A & B & nB \\ B^T & C & nC \\ nB^T & nC & n^2C} \succeq 0
        \label{eqn:SchurLemma}
    \end{align}
    is true for all $n\in\mathbb{R}_{\geq 0}$.
    \end{lemma}
    \begin{pf}
        By Schur complement, the right-hand side of \autoref{eqn:SchurLemma} is equivalent to 
        \begin{align*}
        \bmat{A & B \\ B^T & C} - \bmat{nB \\ nC}(n^2C)^{-1}\bmat{nB^T & nC} \succeq 0
        \end{align*}
        and $n^2C\succeq 0$. The former expression resolves to 
        \begin{align*}
        \bmat{A - BC^{-1}B^T & 0 \\ 0 & 0} \succeq 0.
        \end{align*}
        The result follows by applying the Schur complement again to recover the left-hand side of \autoref{eqn:SchurLemma}.
    \end{pf}
 Next, \autoref{thm:KYPparticular} provides sufficient delay-distribution-dependent LMI conditions for an input-delayed LTI system to be stochastically $(Q,S,R)$-dissipative.
	\begin{theorem} \label{thm:KYPparticular}
		Consider a complete probability space $(\Omega,\mathcal{F},\mathbb{P})$ where $\Omega = \mathbb{Z}_{[w_m,w_M]}$ and any two samples, $w_i,w_j\in\Omega$, are i.i.d. Further, consider the operator $G:\mathcal{X}^m_e \times \mathcal{F} \rightarrow \mathcal{X}^p_e$, where $y=G(u,w)$ is given by
	\begin{align}
		G:\begin{cases}
		x(k+1) = Ax(k) + Bu(k-w_k) \\
		y(k) = Cx(k) + Du(k-w_k), 
		\end{cases} \label{eqn:linear}
	\end{align}
	where $A\in \bb{R}^{n\times n}$, $B\in\bb{R}^{n\times m}$, $C\in\bb{R}^{p\times n}$ and $D\in\bb{R}^{p\times m}$ are known constant matrices. The system, $G$, is QSR-dissipative-in-expectation if there exist  $M_1$, $M_2$, $M_3$, $N_1$, $N_2$, $N_3$, $W_1$, $W_2$, $W_3 \in  \bb{R}^{m \times m}$, positive semi-definite $P\in\bb{S}^{(n+m)\times(n+m)}$, $X$, $Y \in  \bb{S}^{m \times m}$, and positive definite $Z_1$, $Z_2 \in \bb{S}^{m \times m}$ such that
		\begin{align}
            \underset{w_k\sim\mathbb{P}}{\mathbb{E}}\left( \; \bmat{\widetilde{\Pi} & \widetilde{N} & \widetilde{M}& \widetilde{W} \\ \widetilde{N}^T & Z_2 & 0 & 0 \\ \widetilde{M}^T & 0 & Z_1 & 0 \\ \widetilde{W}^T& 0 & 0 & Z_1 } \; \right) \succeq 0 \label{eqn:KYPLMI}
		\end{align}
		where
		\begin{align*}
			\widetilde{\Pi} &= \bmat{\Pi_{11} & \Pi_{12} {+} \Pi_{13} & (w_M{-}w_k) \Pi_{12} \\ * & \Pi_{22} {+} \Pi_{33} {+} \He(\Pi_{23}) & (w_M{-}w_k)(\Pi_{22} {+} \Pi_{23}) \\ * & * &  (w_M{-}w_k)^2\Pi_{22}}, \\
			\Pi &= \bmat{\Pi_{11} {\in}\mathbb{S}^{n+2m} & \Pi_{12} {\in} \mathbb{R}^{(n+2m) \times m} & \Pi_{13} {\in} \mathbb{R}^{(n+2m) \times m} \\ * & \Pi_{22}{\in}\mathbb{S}^{m} & \Pi_{23} {\in}\mathbb{R}^{m\times m} \\ * & * & \Pi_{33}{\in} \mathbb{S}^m } \\ &= \Pi_0 - \Pi_1 - \Pi_2 - \Pi_3 - \Pi_4, \\
			\Pi_0 &= \bmat{C^TQC & C^TS & 0 & C^TQD & 0 \\ * & R & 0 & S^TD & 0 \\ * & * & 0_{m\times m} & 0 & 0 \\ * & * & * & D^TQD & 0 \\ * & * & * & * & 0_{m\times m}},
		\end{align*}
		\begin{align*}
			\Pi_1 &= \bmat{A^T & 0_{n\times m} \\ 0_{m\times n} & I_m \\ 0_{m \times n} & 0_{m \times m} \\ B^T & 0_{m\times m} \\ 0_{m\times n} & 0_{m\times m}}P\bmat{A^T & 0_{n\times m} \\ 0_{m\times n} & I_m \\ 0_{m \times n} & 0_{m \times m} \\ B^T & 0_{m\times m} \\ 0_{m\times n} & 0_{m\times m}}^T 
            \\ &\hspace{.5cm} - 
            \bmat{I_n & 0_{n\times m} \\ 0_{m\times n} & 0_{m\times m} \\ 0_{m\times n} & I_m \\ 0_{m\times n} & 0_{m\times m} \\ 0_{m\times n} & 0_{m\times m}}P\bmat{I_n & 0_{n\times m} \\ 0_{m\times n} & 0_{m\times m} \\ 0_{m\times n} & I_m \\ 0_{m\times n} & 0_{m\times m} \\ 0_{m\times n} & 0_{m\times m}}^T,\\	
            \Pi_2 &= \bmat{0_{(n+m)\times(n+m)} & * & * & * \\ 0 & (w_M{-}w_m{+}1)X{+}Y & * & * \\ 0 & 0 & -X & * \\ 0 & 0 & 0 & -Y},\\
			\Pi_3 &= (w_M{-}1)\bmat{0_{n\times m} \\ I_m \\ -I_m \\ 0_{2m\times m}}(Z_1 + Z_2)\bmat{0_{n\times m} \\ I_m \\ -I_m \\ 0_{2m\times m}}^T,  \\
			\Pi_4 &= \bmat{0_{(n+m)\times (n+m)} & * \\ 0 & \He(\bmat{M{+}N & W{-}M & {-}W{-}N})}, \\
			\widetilde{N} &= \sqrt{w_M-1}\bmat{0_{(n+m)\times m} \\ N_1 \\ N_2+N_3 \\ (w_M{-}w_k)N_2}, \; N = \bmat{N_1 \\ N_2 \\ N_3}
			\\
			\widetilde{M} &= \sqrt{w_k-1}\bmat{0_{(n+m)\times m} \\ M_1 \\ M_2+M_3 \\ (w_M{-}w_k)M_2}, \; M =  \bmat{M_1 \\ M_2 \\ M_3},
			\\
			\widetilde{W} &= \sqrt{w_M-w_k}\bmat{0_{(n+m)\times m} \\  W_1 \\ W_2-W_3 \\ (w_M{-}w_k)W_2}, \; W = \bmat{W_1 \\ W_2 \\ W_3}.
		\end{align*}
	\end{theorem}
	
	\begin{pf}
		The outline of the proof is as follows. First, an L-K-like function is proposed based on existing literature, and its difference operator is bounded to achieve a quadratic expression in $\theta^T(k,w_k) \defeq \bmat{x^T(k) & u^T(k) & u^T(k{-}1) & u^T(k{-}w_k) & u^T(k{-}w_M)}$, where Lyapunov-like variables $P$, $X$, $Y$, $Z_1$, and $Z_2$ relate each of these terms to one another. Then, \autoref{thm:genKYP} is applied to derive an inequality condition for stochastic $(Q,S,R)$-dissipativity using a bound on the L-K function's expected difference. The resulting condition is impractical to verify directly, so it is reformulated using Newton's series of finite differences and Lemma~\ref{lem:RecursiveSchur}, which leads to the sufficient matrix inequality condition in \autoref{eqn:KYPLMI}. Intuition for the choice of $\theta(k,w_k)$ is as follows. First, the terms $x(k)$, $u(k)$, and $u(k-w_k)$ appear explicitly in \autoref{thm:genKYP}. Further, applying the difference operator to $u(k)$ results in $u(k+1)-u(k)$, or, to obey causality, the difference operator on $u(k-1)$ results in $u(k)-u(k-1)$. Lastly using the maximum delay $u(k-w_M)$ will turn out to be useful at several points.
  
        Motivated by Reference~\cite{Gao2007}, consider the L-K-like function $V(k,x_{[k]},u_{[k]},w_{[k]}) {=} \sum_{i=1}^5 V_i(k,x_{[k]},u_{[k]},w_{[k]}) {\geq} 0$, where
        \begin{align*}
            V_1(k,x_{[k]},u_{[k]},w_{[k]}) &= \bmat{x^T(k) & u^T(k{-}1)}P\bmat{x(k) \\ u(k{-}1)},\\
			V_2(k,x_{[k]},u_{[k]},w_{[k]}) &= \hspace{-4mm} \sum_{i=k-w_k+1}^{k-1} u^T(i{-}1)Xu(i{-}1), \\
			V_3(k,x_{[k]},u_{[k]},w_{[k]}) &= \hspace{-4mm}\sum_{i=k-w_M+1}^{k-1} u^T(i{-}1)Yu(i{-}1), \\
			V_4(k,x_{[k]},u_{[k]},w_{[k]}) &= \hspace{-4mm}\sum_{j=-w_M+2}^{-w_m+1}\sum_{i=k+j}^{k-1} u^T(i{-}1)Xu(i{-}1), \\
			V_5(k,x_{[k]},u_{[k]},w_{[k]}) &= \hspace{-4mm}\sum_{i=-w_M+1}^{-1}\sum_{m=k+i}^{k-1} \hspace{-2mm}\eta^T(m)(Z_1{+}Z_2)\eta(m),
		\end{align*}
		 and $\eta(k) \defeq \Delta u(k-1) = u(k) - u(k{-}1)$, and $P$, $X$, $Y$, $Z_1$, and $Z_2$ are defined in \autoref{thm:KYPparticular}. Its expected difference is 
		\begin{align*}
            \underset{\mathbb{E}}{\Delta} V&(k,x_{[k]},u_{[k]},w_{[k]}) \\ &= \sum_{i=1}^5 \underset{\mathbb{E}}{\Delta} V_i(k,x_{[k]},u_{[k]},w_{[k]}) \\ &\defeq \sum_{i=1}^5 \underset{w_{[k+1]}\sim\mathbb{P}}{\mathbb{E}}\Big(V_i(k{+}1,x_{[k+1]},u_{[k+1]},w_{[k+1]}) \\ & \hspace{2.4cm}- V_i(k,x_{[k]},u_{[k]},w_{[k]})\Big).
            \end{align*}
            Dropping the arguments of $\Delta_{\mathbb{E}}V_i$,
            \begin{align*}
                \underset{\mathbb{E}}{\Delta}V_1 &=\underset{w_k\sim\mathbb{P}}{\mathbb{E}}\Big(\mathrm{Sq}\Big(P,\bmat{Ax(k) + Bu(k-w_k)\\ u(k)}\Big) \\ & \hspace{1.3cm}- \mathrm{Sq}\Big(P,\bmat{x(k) \\ u(k{-}1)}\Big)\Big), \\
                \underset{\mathbb{E}}{\Delta}V_2 &= \underset{w_{[k+1]}\sim\mathbb{P}}{\mathbb{E}}\Big(\sum_{i=k-w_{k+1}+2}^{k} \mathrm{Sq}(X,u(i-1)) \\ & \hspace{1.8cm} - \sum_{i=k-w_k+1}^{k-1} \mathrm{Sq}(X,u(i-1))\Big)  \\
			&\leq \underset{w_k\sim\mathbb{P}}{\mathbb{E}}\Big(\mathrm{Sq}(X,u(k-1)) - \mathrm{Sq}(X,u(k-w_k)) \\ & \hspace{1.3cm}+ \sum_{i=k-w_M+2}^{k-w_m+1} \mathrm{Sq}(X,u(i-1))\Big), \\
            \underset{\mathbb{E}}{\Delta}V_3 &= \underset{w_{[k+1]}\sim\mathbb{P}}{\mathbb{E}}\Big(\sum_{i=k-w_M+2}^{k} \mathrm{Sq}(Y,u(i-1))  \\ & \hspace{1.8cm} - \sum_{i=k-w_M+1}^{k-1} \mathrm{Sq}(Y,u(i-1)) \Big) \\
			&= \mathrm{Sq}(Y,u(k-1)) - \mathrm{Sq}(Y,u(k-w_M)),
            \end{align*}
            \begin{align*}
			\underset{\mathbb{E}}{\Delta} V_4 &=  \underset{w_{[k+1]}\sim\mathbb{P}}{\mathbb{E}}\Big(\sum_{j=-w_M+2}^{-w_m+1}\sum_{i=k+1+j}^{k} \mathrm{Sq}(X,u(i-1)) \\ & \hspace{1.8cm} -  \sum_{j=-w_M+2}^{-w_m+1}\sum_{i=k+j}^{k-1} \mathrm{Sq}(X,u(i-1))\Big)  \\
			&= (w_M {-} w_m) \mathrm{Sq}(X,u(k{-}1)) {-}\hspace{-3mm} \sum_{i=k{-}w_M{+}2}^{k{-}w_m{+}1} \hspace{-4mm}\mathrm{Sq}(X,u(i{-}1)), \\
			\underset{\mathbb{E}}{\Delta} V_5 &= \underset{w_{[k+1]}\sim\mathbb{P}}{\mathbb{E}}\Big(\sum_{i=-w_M+1}^{-1}\sum_{m=k+1+i}^{k}\mathrm{Sq}((Z_1+Z_2),\eta(m)) \\ & \hspace{1.8cm} - \hspace{-3mm} \sum_{i=-w_M+1}^{-1}\sum_{m=k+i}^{k-1}\mathrm{Sq}((Z_1+Z_2),\eta(m))\Big) \\
			&= \underset{w_{k}\sim\mathbb{P}}{\mathbb{E}}\Big((w_M{-}1)\mathrm{Sq}((Z_1 {+} Z_2),\eta(k)) \\ & \hspace{1.4cm}- \sum_{i{=}k{-}w_M{+}1}^{k{-}1}\mathrm{Sq}(Z_2,\eta(i)) \\ & \hspace{1.4cm} - \sum_{i{=}k{-}w_k{+}1}^{k{-}1}\mathrm{Sq}(Z_1,\eta(i)) \\ & \hspace{1.4cm} - \sum_{i{=}k{-}w_M{+}1}^{k{-}w_k}\mathrm{Sq}(Z_1,\eta(i))\Big).
		\end{align*}
		To see how the bound on $\Delta_\mathbb{E} V_2$ is derived, consider that $\Delta_\mathbb{E} V_2 = \mathbb{E}_{w_k\sim\mathbb{P}}\big(\mathrm{Sq}(X,u(k-1)) - \mathrm{Sq}(X,u(k-w_k)) + R\big)$, where $R$ is a remainder term that depends on $w_{k+1}$ and $w_k$. This remainder is largest when $w_{k+1} = w_M$ and $w_k = w_m$, which yields the expression above. This bound is useful because the summations in $\Delta_\mathbb{E} V_2$ and $\Delta_\mathbb{E} V_4$ cancel, which means $\sum_{i=1}^4 \Delta_\mathbb{E} V_i$ is quadratic in $\theta(k)$. To achieve the same form with $\Delta_\mathbb{E} V_5$, the summations are eliminated by adding zero, completing the square, and over-bounding, as follows. First, add zero to $\Delta_\mathbb{E} V_5$ with 
        \begin{align*}
			0 &= 2\zeta^T(k,w_k)\bmat{M{+}N, & W{-}M, & {-}W{-}N}\zeta(k,w_k)\\
			& \hspace{1cm} - 2\zeta^T(k,w_k)M\sum\limits_{i=k-w_k+1}^{k-1}\eta(i) \\ & \hspace{1cm} - 2\zeta^T(k,w_k)W\sum\limits_{i=k-w_M+1}^{k-w_k}\eta(i) \\ & \hspace{1cm} - 2\zeta^T(k,w_k)N\sum\limits_{i=k-w_M+1}^{k-1}\eta(i),
	\end{align*}
        where $\zeta^T(k,w_k) \defeq \bmat{u^T(k{-}1) & u^T(k{-}w_k) & u^T(k{-}w_M)}$, and $M$, $N$, and $W$ are matrices defined in \autoref{eqn:KYPLMI}. To see that this expression is zero, substitute $\bmat{u^T(k{-}1) & u^T(k{-}w_k) & u^T(k{-}w_M)}^T$ for $\zeta(k,w_k)$ on the first row. Then group terms containing $M$, $N$, and $W$ separately. Last, use the identity $u(k-b)-u(k-a)=\sum_{i=k-a+1}^{k-b} \eta(i)$. After adding this zero, $\Delta_\mathbb{E} V_5$ contains three summations that can be bounded above by completing the square. For example, consider the term
        \begin{align*}
            &-\sum_{i=k-w_k+1}^{k-1} \Big(\mathrm{Sq}(Z_1,\eta(i))+2\zeta^T(k,w_k)M\eta(i)\Big) \\ 
            =& -\sum_{i=k-w_k+1}^{k-1} \Big(\mathrm{Sq}(Z_1,\eta(i))+2\zeta^T(k,w_k)M\eta(i) \\ & \hspace{.5cm} + \mathrm{Sq}(Z_1^{-1},M^T\zeta(k,w_k))\Big) + \mathrm{Sq}(Z_1^{-1},M^T\zeta(k,w_k)) \\
            =& -\sum_{i=k-w_k+1}^{k-1} \Big(\mathrm{Sq}\big(Z_1,\eta(i) + Z_1^{-1}M\zeta(k,w_k)\big)\Big) \\ & \hspace{.5cm} + \mathrm{Sq}(Z_1^{-1},M^T\zeta(k,w_k)) \\
            \leq&\hspace{2mm} \mathrm{Sq}(Z_1^{-1},M^T\zeta(k,w_k)).
        \end{align*}
        Applying this reasoning to all three sums results in
	\begin{align*}
		\underset{\mathbb{E}}{\Delta} V_5 &=  \underset{w_{k}\sim\mathbb{P}}{\mathbb{E}}\Big((w_M{-}1)\mathrm{Sq}((Z_1 {+} Z_2),\eta(k)) \\ & \hspace{5mm}+ 2\mathrm{Sq}(\Phi, \zeta(k,w_k)) \\
			& \hspace{5mm} + (w_k-1) \mathrm{Sq}(Z_1^{-1},M^T \zeta(k,w_k)) \\
                &\hspace{5mm} + (w_M-w_k) \mathrm{Sq}(Z_1^{-1},W^T 
            \zeta(k,w_k)) \\
			&\hspace{5mm} + (w_M-1) \mathrm{Sq}(Z_2^{-1},N^T \zeta(k,w_k)) \\
	           &\hspace{5mm} - \sum_{i=k-w_k+1}^{k-1}\Big( \mathrm{Sq}(Z_1,      \eta(i)) + 2\zeta^T(k,w_k)M \eta(i) \\ & \hspace{2.5cm} + \mathrm{Sq}(Z_1^{-1},M^T \zeta(k,w_k))\Big) \\
			     &\hspace{5mm} - \sum_{i=k-w_M+1}^{k-w_k} \Big( \mathrm{Sq}(Z_1,\eta(i)) + 2\zeta^T(k,w_k)W\eta(i) \\ & \hspace{2.5cm} + \mathrm{Sq}(Z_1^{-1},W^T\zeta(k,w_k))\Big)\Big) \\
                  &\hspace{5mm} - \sum_{i=k-w_M+1}^{k-1} \Big( \mathrm{Sq}(Z_2,\eta(i)) + 2\zeta^T(k,w_k) N\eta(i) \\ & \hspace{2.5cm} + \mathrm{Sq}(Z_2^{-1},N^T\zeta(k,w_k))\Big) \\ 
		&\leq \underset{w_{k}\sim\mathbb{P}}{\mathbb{E}}\Big((w_M-1)\mathrm{Sq}((Z_1 + Z_2)\eta(k)) \\ & \hspace{5mm} + \mathrm{Sq}(\He(\Phi) {+} \Xi,\zeta(k,w_k))\Big),
	\end{align*}
		where $\Phi = \bmat{M+N, & W-M, & -W-N}$ and $\Xi(w_k) = (w_k-1)MZ_1^{-1}M^T + (w_M-w_k)WZ_1^{-1}W^T+ (w_M-1)NZ_2^{-1}N^T$. Equation~\ref{eqn:KYPgen} from Theorem~\ref{thm:genKYP} is now rewritten as $\mathbb{E}_{w_{k}\sim\mathbb{P}}\big(\theta^T(k,w_k)\Pi_0\theta(k,w_k)\big) - \overline{\Delta_\mathbb{E} V} \geq 0$, where $\overline{\Delta_\mathbb{E} V} = \Delta_\mathbb{E} V_1 + \overline{\Delta_\mathbb{E} V_2} + \Delta_\mathbb{E} V_3 + \Delta_\mathbb{E} V_4+ \overline{\Delta_\mathbb{E} V_5}$, $\Delta_\mathbb{E} V_1 = \theta^T(k,w_k)\Pi_1\theta(k,w_k)$, $\overline{\Delta_\mathbb{E} V_2} + \Delta_\mathbb{E} V_3+ \Delta_\mathbb{E} V_4 = \theta^T(k,w_k)\Pi_2\theta(k,w_k)$, $\overline{\Delta_\mathbb{E} V_5} =\mathrm{Sq}(\Pi_3 + \Pi_4 + \Psi(w_k),$ $\theta(k,w_k))$, and $\Psi(w_k) = \mbox{diag}(0_{(n+m)\times(n+m)},\Xi(w_k))$. Letting $\Pi = \Pi_0- \Pi_1 - \Pi_2 - \Pi_3 - \Pi_4$ and $\Upsilon(w_k) = \Pi - \Psi(w_k)$, Equation~\ref{eqn:KYPgen} reduces to, for all $k\in\mathbb{Z}$, $x(0)\in\mathbb{R}^n$, and $u_{[k]}\in\mathcal{X}^m$,
		\begin{align}
			\mathbb{E}_{w_{k}\sim\mathbb{P}}\big(\theta^T(k,w_k)\Upsilon(w_k)\theta(k,w_k)\big) &\geq 0.\label{eqn:KYPLMIproof}
		\end{align}

Equation~\ref{eqn:KYPLMIproof} is a sufficient condition for $(Q,S,R)$-dissipativity, but it is not easily verifiable. To convert this to an LMI, the vectors $\theta$ must be ``pulled out'' of the expectation; however, $\theta$ is a function of the random variable $w_k$. To separate the vectors from the random variable, consider Newton's Series of finite differences \cite[pg.357]{Jordan1939}\cite{Konig2021}, which gives the expression
\begin{align*}
    u(k_1) = \sum_{i=0}^\infty \begin{pmatrix} k_1{-}k_2 \\ i \end{pmatrix} \Delta^i [u](k_2),
\end{align*}
where $\begin{pmatrix} \alpha \\ \beta \end{pmatrix} = \frac{\alpha!}{\beta!(\alpha-\beta)!}$ is the binomial, and $\Delta^i [u](k_2) = \sum_{l=0}^i (-1)^{i-l}\begin{pmatrix} i \\ l \end{pmatrix}u(k_2+l)$ is the $i^{th}$ order finite difference of $u$ at $k_2$. Applying this with $k_1 = k-w_k$, and $k_2 = k-w_M$ yields
\begin{align*}
    u(k{-}w_k) &= \sum_{i=0}^\infty \begin{pmatrix} w_M{-}w_k \\ i \end{pmatrix} \Delta^i[u](k-w_M) 
    \\ &= u(k{-}w_M) + \hspace{-4mm} \sum_{i=1}^{w_M-w_m} \hspace{-2mm}\begin{pmatrix} w_M{-}w_k \\ i \end{pmatrix} \Delta^i[u](k{-}w_M).
\end{align*}
The series is finite because for all $i> w_M-w_m$, the binomial term equals zero for all $w_k$. This may also be viewed as the $(w_M{-}w_m)^{th}$-order Lagrange interpolation polynomial, which exactly reproduces $u(n)$ at the integer values $n=k{-}w_M,k{-}w_M{+}1,\dots,k{-}w_m$. To simplify notation, denote $L_i(u,k) \defeq \Delta^i[u](k-w_M)$. Then 
\begin{align*}
    \theta(k,w_k) &= \bmat{x(k) \\ \hspace{.2cm}u(k) \\ u(k-1) \\  u(k{-}w_M) {+} \hspace{-3mm} \sum\limits_{i=1}^{w_M-w_m} \begin{pmatrix} w_M{-}w_k \\ i \end{pmatrix} L_i(u,k) \\ u(k-w_M)}.
\end{align*}
Now the random variable, $w_k$, appears as a coefficient in $\theta$, so $w_k$ can be grouped instead with the matrix, $\Upsilon$, by defining the new vector 
\begin{align*} 
z(k) = \bmat{x(k) \\ u(k) \\ u(k-1) \\ u(k-w_M) \\ L_1(u,k) \\ L_2(u,k) \\ \vdots \\ L_{w_M-w_m}(u,k)}
\end{align*}
and rewriting \autoref{eqn:KYPLMIproof} as $ z^T(k)\mathbb{E}_{w_k\sim \mathbb{P}}\big(\widetilde{\Upsilon}(w_k)\big)z(k) {\geq} 0$ for all $k\in\mathbb{Z}$, $x(0)\in\mathbb{R}^n$, and $u_{[k]}\in\mathcal{X}^m$. This is now implied by $\mathbb{E}_{w_k\sim \mathbb{P}}\big(\widetilde{\Upsilon}(w_k)\big) \succeq 0$
 where
 \begin{align*}
  \widetilde{\Upsilon}(w_k) &= \bmat{\widetilde{\Upsilon}_{1,1} & \widetilde{\Upsilon}_{1,2} & \widetilde{\Upsilon}_{1,3} & \widetilde{\Upsilon}_{1,4} & \dots & \widetilde{\Upsilon}_{1,\tilde{n}{+}2} \\ * & \widetilde{\Upsilon}_{2,2} & \widetilde{\Upsilon}_{2,3} & \widetilde{\Upsilon}_{2,4} & \dots & \widetilde{\Upsilon}_{2,\tilde{n}{+}2} \\ * & * & \widetilde{\Upsilon}_{3,3} & \widetilde{\Upsilon}_{3,4} & \dots & \widetilde{\Upsilon}_{3,\tilde{n}{+}2} \\ * & * & * & \widetilde{\Upsilon}_{4,4} & \dots & \widetilde{\Upsilon}_{4,\tilde{n}{+}2} \\ \vdots & \vdots & \vdots & \vdots & \ddots & \vdots \\ * & * & * & * & \dots & \widetilde{\Upsilon}_{\tilde{n}{+}2,\tilde{n}{+}2}  },
  \\
  \tilde{n} &= w_M-w_m,
  \\
  \widetilde{\Upsilon}_{1,1} & = \Upsilon_{11},
  \\ 
  \widetilde{\Upsilon}_{1,2} & = \Upsilon_{12}+\Upsilon_{13},
  \\
  \widetilde{\Upsilon}_{1,3} & = \pmat{w_M{-}w_k \\ 1}\Upsilon_{12},  
  \\
  \widetilde{\Upsilon}_{1,4} & = \pmat{w_M{-}w_k \\ 2}\Upsilon_{12},
  \\
  \widetilde{\Upsilon}_{1,\tilde{n}+2} & = \pmat{w_M{-}w_k \\ \tilde{n}}\Upsilon_{12},
  \\
  \widetilde{\Upsilon}_{2,2} & = \Upsilon_{22}+ \Upsilon_{33} + \He(\Upsilon_{23}), 
  \\
  \widetilde{\Upsilon}_{2,3} & = \pmat{w_M{-}w_k \\ 1}(\Upsilon_{22} + \Upsilon_{23}), 
  \\
  \widetilde{\Upsilon}_{2,4} & = \pmat{w_M{-}w_k \\ 2}(\Upsilon_{22} + \Upsilon_{23}), \\
    \widetilde{\Upsilon}_{2,\tilde{n}+2} & = \pmat{w_M{-}w_k \\ \tilde{n}}(\Upsilon_{22} + \Upsilon_{23}), \\
      \widetilde{\Upsilon}_{3,3} & = \pmat{w_M{-}w_k \\ 1}^2\Upsilon_{22}, \\
  \widetilde{\Upsilon}_{3,4} & = \pmat{w_M{-}w_k \\ 1}\pmat{w_M{-}w_k \\ 2}\Upsilon_{22}, \\
  \widetilde{\Upsilon}_{3,\tilde{n}+2} & = \pmat{w_M{-}w_k \\ 1}\pmat{w_M{-}w_k \\ \tilde{n}}\Upsilon_{22},
  \\
  \widetilde{\Upsilon}_{4,4} & = \pmat{w_M{-}w_k \\ 2}^2\Upsilon_{22},
  \\
  \widetilde{\Upsilon}_{4,\tilde{n}+2} & = \pmat{w_M{-}w_k \\ 2}\pmat{w_M{-}w_k \\ \tilde{n}}\Upsilon_{22},
  \\
  \widetilde{\Upsilon}_{\tilde{n}{+}2,\tilde{n}{+}2} & = \pmat{w_M{-}w_k \\ \tilde{n}}^2\Upsilon_{22}, \mbox{ and}
\end{align*}
\begin{align*}
  \Upsilon(w_k) &= \bmat{\Upsilon_{11} {\in}\mathbb{S}^{n+2m} & * & * \\ \Upsilon_{12}^T {\in} \mathbb{R}^{m \times (n+2m)} & \Upsilon_{22}{\in}\mathbb{S}^{m} & * \\ \Upsilon_{13}^T {\in} \mathbb{R}^{m \times (n+2m)} & \Upsilon_{23}^T {\in}\mathbb{R}^{m\times m} & \Upsilon_{33}{\in} \mathbb{S}^m }.
        \end{align*}
    Applying Lemma~\ref{lem:RecursiveSchur} recursively to $\widetilde{\Upsilon}(w_k)$ reveals that $\mathbb{E}_{w_k\sim \mathbb{P}}\big(\widetilde{\Upsilon}(w_k)\big) \succeq 0$ is equivalent to $\mathbb{E}_{w_k\sim \mathbb{P}}\big(\widetilde{\Upsilon}'(w_k) \big) \succeq 0$, where
    \begin{align*}
        \widetilde{\Upsilon}&'(w_k) = \\ & \bmat{\Upsilon_{11} & \Upsilon_{12} + \Upsilon_{13} & (w_M{-}w_k)\Upsilon_{12} \\ * & \Upsilon_{22} + \Upsilon_{33} + \He(\Upsilon_{23}) & (w_M{-}w_k)(\Upsilon_{22} + \Upsilon_{23})  \\ * & * &  (w_M{-}w_k)^2\Upsilon_{22} }.
    \end{align*}
    Therefore, $\mathbb{E}_{w_{k}\sim\mathbb{P}}\big(\widetilde{\Upsilon}'(w_k)\big) \succeq 0$ implies \autoref{eqn:KYPLMIproof} which implies \autoref{eqn:KYPgen}. By \autoref{thm:genKYP}, this implies that $G:\Uspace^m_e \times \mathcal{F}\rightarrow\Uspace^p_e$ is $(Q,S,R)$-dissipative-in-expectation. Furthermore, $\widetilde{\Upsilon}' = \widetilde{\Pi}(w_k) - \widetilde{\Psi}(w_k)$, where $\widetilde{\Psi}(w_k) = \widetilde{M}Z_1^{-1}\widetilde{M}^T + \widetilde{W}Z_1^{-1}\widetilde{W}^T + \widetilde{N}Z_2^{-1}\widetilde{N}^T$, and $\widetilde{\Pi}$, $\widetilde{M}$, $\widetilde{W}$, and $\widetilde{N}$ are given in \autoref{eqn:KYPLMI}. Applying the Schur complement results in the equivalent LMI given by \autoref{eqn:KYPLMI}.
	\end{pf}
	
	\subsection{Comparison with Deterministic Results} \label{sec:compare}
	
	Theorem~\ref{thm:KYPparticular} is the first delay-distribution-dependent LMI condition for stochastic dissipativity of discrete LTI systems with time-varying input delay. To the authors' knowledge, there is also no comparable deterministic result. The closest result may be Reference~\cite{ElAiss2017}, which uses the average of the upper and lower state delay bounds to construct a less conservative L-K function. Because we consider input delays here, even the work of Reference~\cite{ElAiss2017} cannot be compared directly with the current result. To provide a basis for comparison, Theorem~\ref{thm:KYPparticular} can be used to derive delay-bound-dependent conditions for deterministic dissipativity, which is more conservative but relies on less information about the delay. Here, this condition is derived for comparison. The proof is identical to that of Theorem~\ref{thm:KYPparticular} through completing the square on $\Delta_\mathbb{E} V_5$, then $\Delta_\mathbb{E} V_5$ is bounded with $\Xi \preceq \overline{\Xi} = (w_M-1)NZ_2^{-1}N^T + (w_M-1)MZ_1^{-1}M^T + (w_M-w_m)WZ_1^{-1}W^T$, as in Reference~\cite{Gao2007}, so $\Psi \preceq \overline{\Psi} = \mbox{diag}(0_{n+m},\overline{\Xi})$. The term $\overline{\Psi}$ is no longer a function of $w_k$, so $\Pi - \overline{\Psi} \succeq 0$ provides a sufficient condition for deterministic dissipativity. This is linearized using Schur complement as
	\begin{align}
		\bmat{\Pi & * & * & * \\ \bmat{0 & \sqrt{w_M-1}N^T} & Z_2 & * & * \\ \bmat{0 & \sqrt{w_M-1}M^T} & 0 & Z_1 & * \\ \bmat{0 & \sqrt{w_M-w_m}W^T} & 0 & 0 & Z_1} \succeq 0, \label{eqn:deterministicKYP}
	\end{align}
which replaces \autoref{eqn:KYPLMI} in Theorem~\ref{thm:KYPparticular}.

%% file: Numerical.tex
In this section, Theorem~\ref{thm:KYPparticular} is used to characterize the dissipativity of a simple numerical system. The results are compared to the deterministic LMI from Equation~\ref{eqn:deterministicKYP} in order to illustrate some key advantages of the stochastic method. The stochastic dissipativity characterization is also used to design a static output feedback (SOF) controller with the maximum allowable gain that preserves stochastic input-output stability on $\Uspace$ according to Theorems~\ref{thm:OLstable} and \ref{thm:CLstable}.
	
	This section focuses the analysis on single-input-single-output (SISO) LTI systems, because the representation of dissipativity is easily visualized on a Nyquist plot, which allows a clear illustration of the advantages and disadvantages of Theorem~\ref{thm:KYPparticular}. However, this theorem can also be used for non-square multiple-input-multiple output (MIMO) systems, and Theorems~\ref{thm:OLstable}, \ref{thm:CLstable}, and \ref{thm:genKYP} hold nonlinear systems as well.
	
	For SISO systems, QSR-dissipativity is equivalent to Zames' conic sectors \cite{Zames1966}. For SISO LTI systems with constant time delay, conic sectors can be visualized as the smallest circle that circumscribes the Nyquist plot. This circle can be characterized by its center, $c$, and radius, $r$, with the equivalent dissipativity characterization $Q = -1$, $S = c$, $R = r^2-c^2$. Alternatively, it can be parameterized by the circle's minimum and maximum values on the real axis, $a$ and $b$. The equivalent dissipativity characterization for this is $Q = -1$, $S = \frac{a+b}{2}$, and $R = -ab$. Gain and passivity are special cases of conic sectors, where the former is a circle centered at the origin ($Q=-1$, $S = 0$, $R = r^2$), and the latter is the degenerate circle with minimum value at 0 and maximum value approaching infinity ($Q=0$, $S = \frac{1}{2}$, $R = 0$). When identifying conic sectors, there are two well known methods for achieving a ``tight" characterization. One method \cite{Joshi2002} is to minimize the conic radius by minimizing $R+S^TS$. An alternative \cite{Bridgeman2014} is to set $b\rightarrow \infty$, which yields the dissipative representation $Q=0$, $S = \frac{1}{2}$, $R = -a$,  and maximize $a$; then fix $a$ and minimize $b$. If a finite $b$ is possible, this dissipative representation returns to $Q = -1$, $S = \frac{a+b}{2}$, $R = -ab$. (In some cases, numerical problems ensue in maximizing $a$, which can be remedied by setting $b$ to a large constant). 

\begin{figure} 
\centering
\input{Figures/PlotStochasticNyquist.tex}
	\caption{Nyquist plot of Equation~\ref{eqn:example} with various constant time delays.} \label{fig:Nyquist}
\end{figure}
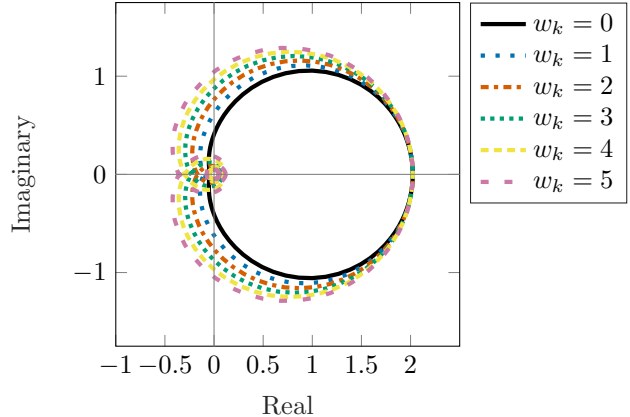

	Consider the SISO LTI system 
\begin{align}
	x(k{+}1) &= \bmat{0.6024 & -0.0038 \\ 0.00381 & 0.9451}x(k) {+} \bmat{0.1647 \\ 0.0960}u(k{-}w_k) \nonumber \\ y(k) &=  \bmat{0 & 1}x(k), \label{eqn:example}
\end{align}
which is modified from the example in Reference~\cite{Lima2021} to be Lyapunov stable. The Nyquist plot of Equation~\ref{eqn:example} with no delay ($w_k=0 \; \forall \, k$) has a maximum gain of 2 and small passivity violation. As the constant delay is increased, the gain remains constant and the passivity violation grows, as demonstrated in Figure~\ref{fig:Nyquist} for constant time delays between 0 and 5.

To demonstrate Theorem~\ref{thm:KYPparticular}, the stochastic dissipativity of Equation~\ref{eqn:example} is analyzed given $w_m = 1$, $w_M = 5$, and various probability density functions for $w_k$.  These cases are enumerated in Table~\ref{tbl:cases} and illustrated in Figure~\ref{fig:cones}. The ``tightest'' conic characterization is identified using three methods: minimizing the gain, minimizing the conic radius, and maximizing $a$ then minimizing $b$. The stochastic dissipativity is compared alongside the deterministic dissipativity for a time-varying delay with the same bounds and unknown probability distribution, using Equation~\ref{eqn:deterministicKYP}.

Several interesting trends are evident. First, all of the stochastic dissipative characterizations are tighter than the deterministic case, even when the delay distribution is uniform (Distribution 3) or concentrated around the largest delay (Distribution 1). Second, the loosest stochastic characterization is provided by the uniform delay distribution, and the characterizations get tighter as the delay distribution becomes more concentrated. This is expected since the more certain information about the delays should correspond to tighter system characterizations. Further, when the delays are concentrated around the minimum value, the tightest lower conic bound is achieved, while the gain characterization is the same for delays concentrated around the minimum and maximum, which makes sense in the context of Figure~\ref{fig:Nyquist}. 

\begin{figure} 
\centering
	\begin{tabular}{|c|c||c|c|c|}
		\hline
		$w_k\sim\dots$ & LMI & Gain & $\min r$ & $\max a$ \\
		\hline
		Unknown & Eqn.~\ref{eqn:deterministicKYP} & 4.5 & $c =0$, $r = 4.5$ & $a = -4.1$ \\
		$\mathbb{P}_1$ & Eqn.~\ref{eqn:KYPLMI} & 2.1 & $c = 0$, $r = 2.1$ & $a = -2.0$ \\
		$\mathbb{P}_2$ & Eqn.~\ref{eqn:KYPLMI} & 2.8 & $c = 0.1$, $r = 2.8$ & $a = -1.9$ \\
		$\mathbb{P}_3$ & Eqn.~\ref{eqn:KYPLMI} & 3.9 & $c=0$, $r = 3.9$ & $a = -2.8$ \\
		$\mathbb{P}_4$ & Eqn.~\ref{eqn:KYPLMI} & 3.6 & $c = 0.9$, $r = 3.4$ & $a = -0.9$\\
		$\mathbb{P}_5$ & Eqn.~\ref{eqn:KYPLMI} & 2.1 & $c = 0.9$, $r = 1.5$ & $a = -0.2$\\
		\hline
	\end{tabular}
 \begin{align*}
     \mathbb{P}_1(w_k=[1, \, 2, \, 3, \, 4, \, 5]) &= [0.01,\, 0.01, \, 0.01, \, 0.01, \, 0.96] \\
     \mathbb{P}_2(w_k=[1, \, 2, \, 3, \, 4, \, 5]) &= [0.05,\, 0.05, \, 0.05, \, 0.1, \, 0.75] \\
     \mathbb{P}_3(w_k=[1, \, 2, \, 3, \, 4, \, 5]) &= [0.2, \, 0.2, \, 0.2, \, 0.2, \, 0.2] \\ 
     \mathbb{P}_4(w_k=[1, \, 2, \, 3, \, 4, \, 5]) &= [0.75,\, 0.1, \, 0.05, \, 0.05, \, 0.05]\\
     \mathbb{P}_5(w_k=[1, \, 2, \, 3, \, 4, \, 5]) &= [0.96,\, 0.01,\, 0.01,\, 0.01,\, 0.01]
 \end{align*}
 \vspace{-1cm}
\caption{Conic sectors identified for Equation~\ref{eqn:example} under different delay probability distributions, $\mathbb{P}_i$. The ``Unknown'' row assumes a bounded but otherwise unknown delay distribution, and the identified conic sectors are deterministic. The remaining 5 cases are stochastic conic sectors identified for various delay distributions. The notation $\mathbb{P}(w_k {=} [a,b,\dots]) = [x,y,\dots]$ denotes that $\mathbb{P}(w_k{=}a) = x$, $\mathbb{P}(w_k{=}b)=y$, etc. The second column shows the LMI used to identify the conic bounds, and the last three columns show the conic descriptions for three different objectives.} \label{tbl:cases}
\end{figure}

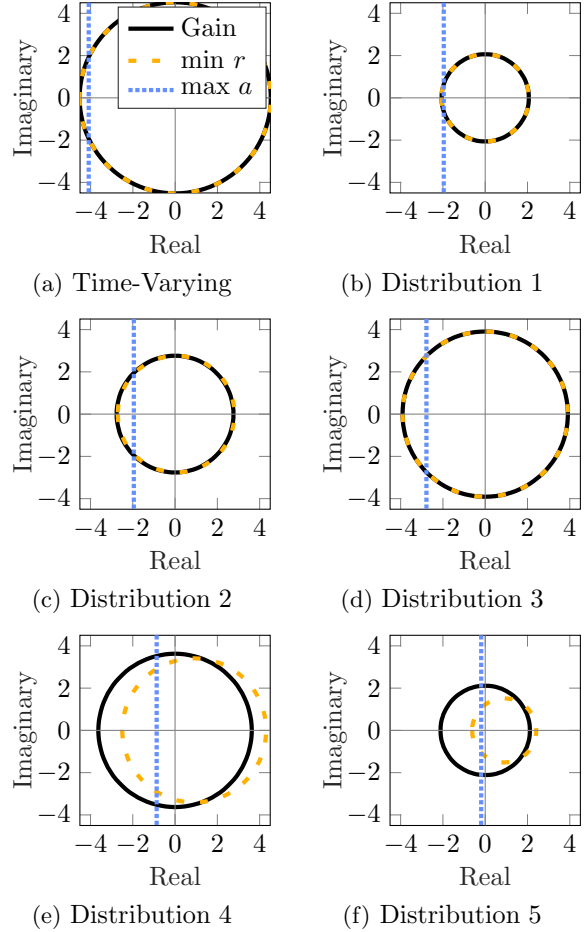
\begin{figure}%
\centering
\subfigure[\normalsize{Time-Varying}]{\input{Figures/PlotTimeVaryingCone.tex}}
\subfigure[\normalsize{Distribution 1}]{\input{Figures/PlotDistribution1Cones.tex}}
\subfigure[\normalsize{Distribution 2}]{\input{Figures/PlotDistribution2Cones.tex}}
\subfigure[\normalsize{Distribution 3}]{\input{Figures/PlotDistribution3Cones.tex}}
\subfigure[\normalsize{Distribution 4}]{\input{Figures/PlotDistribution4Cones.tex}}
\subfigure[\normalsize{Distribution 5}]{\input{Figures/PlotDistribution5Cones.tex}}
\caption{Nyquist plots of (stochastic) conic sectors of Equation~\ref{eqn:example} from Table~\ref{tbl:cases}.}
\label{fig:cones}
\end{figure}

To demonstrate the utility of stochastic dissipativity for controller design, a SOF controller $y_c = Ky$, is designed in negative feedback ($u=-y_c$) with Equation~\ref{eqn:example}, where the plant input delays are distributed according to Distribution 4 in Table~\ref{tbl:cases}. First, the stochastic dissipative characterization from Theorem~\ref{thm:KYPparticular} is employed with the maximum lower conic bound ($a = -.9$, $b=2e4$). According to Theorems~\ref{thm:OLstable} and \ref{thm:CLstable}, the closed loop is stochastically input-output stable on $\mathcal{X}$ if the controller is dissipative-in-expectation with conic bounds $a = 0$, $b = 1.1$. Since the controller has no delay, stochastic dissipativity is equivalent to deterministic dissipativity, so $K=1.1$ is the  SOF controller with the largest permitted gain. For comparison, two other controllers are designed. The first is designed using the deterministic KYP Lemma from Equation~\ref{eqn:deterministicKYP}, which yields a dissipative characterization ($a=-4.1$, $b=2e5$). According to the deterministic Dissipativity Theorem \cite{Vidyasagar1981}, the controller is required to be in the cone $a = 0$, $b = .22$, which is achieved by $K = .22$. Last, the results are compared to a dissipativity-based controller that neglects delays entirely. The undelayed plant satisfies the conic bounds $a = -.055$, $b\rightarrow\infty$, which requires the controller to be in $a = 0$, $b = 18$. Therefore, the second comparison controller is SOF with $K=18$. 

To simulate the SOF controllers, a three-step square wave disturbance is injected to the plant and controller at the first step. The open-loop plant is stable, but converges slowly. The controller designed for deterministic time-varying delays is overly-conservative, and does not provide a significant impact, while the dissipativity-based controller that neglects delays is overly-aggressive and causes instability. In contrast, the stochastic dissipativity-based design provides the desired fast but stable response.

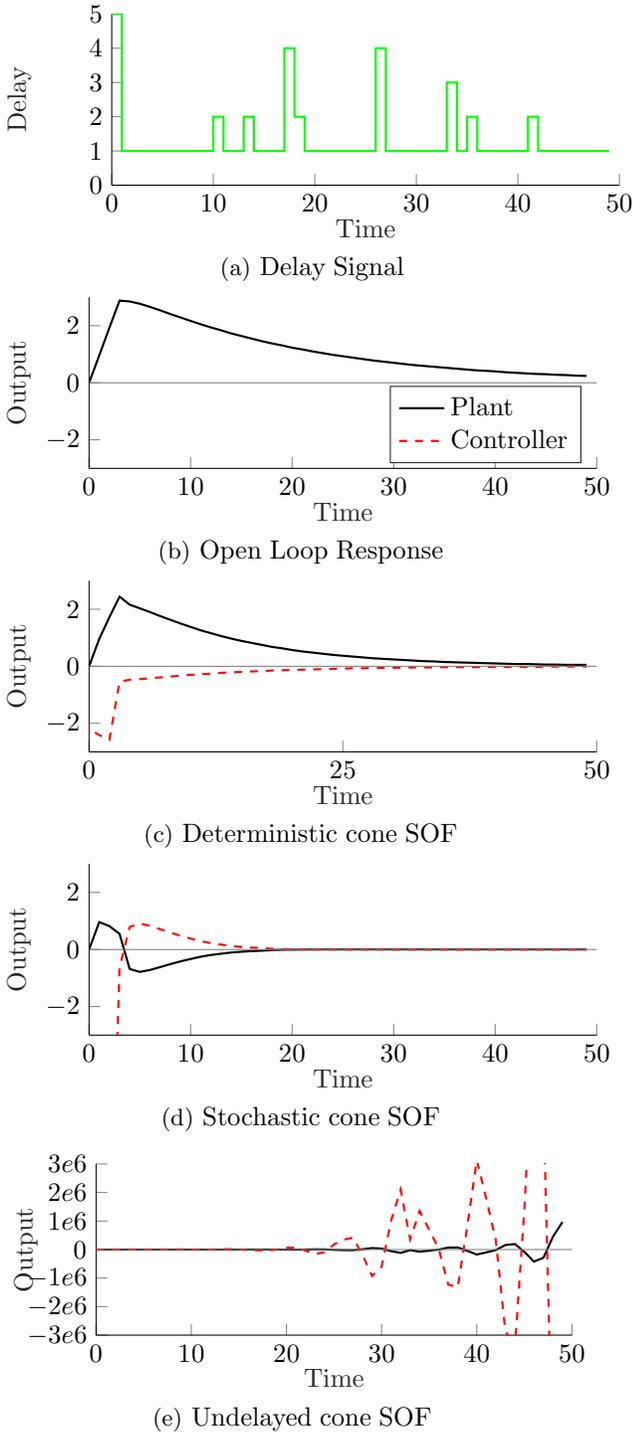
\begin{figure} \label{fig:sim}
\centering
\subfigure[\normalsize{Delay Signal}]{\input{Figures/PlotSimulationDelay.tex}}
\subfigure[\normalsize{Open Loop Response}]{\input{Figures/PlotSimulationOpenLoop}}
\subfigure[\normalsize{Deterministic cone SOF}]{\input{Figures/PlotSimulationDeterministic}}
\subfigure[\normalsize{Stochastic cone SOF}]{\input{Figures/PlotSimulationStochastic}}
\subfigure[\normalsize{Undelayed cone SOF}]{\input{Figures/PlotSimulationUnstable}}
\caption{Simulation of Equation~\ref{eqn:example} with delay distribution 4 (Table~\ref{tbl:cases}) and various SOF controllers in response to three step square wave disturbance of amplitude 10 applied to both plant and controller. Plot (a) gives the particular delay time-series, plot (b) is the system response with no controller, and the controller gain for plot (c) is $K=0.22$, for (d) is $K=1.1$, and (e) is $K=18$.}
\end{figure}

\section{Conclusion}\label{sec:conclusion}

In this paper, stochastic dissipativity and input-output stability was defined for operators on complete inner product spaces. Theorem~\ref{thm:OLstable} related stochastic dissipativity to stochastic input-output stability, and Theorem~\ref{thm:CLstable} derived open-loop stochastic dissipativity conditions for the closed-loop stochastic input-output stability of a network. Subsequently, a KYP Lemma was developed in two steps to verify the stochastic dissipativity of a system with time-varying input delays of a known probability distribution. First, Theorem~\ref{thm:genKYP} demonstrated how to use an L-K functional to verify the stochastic dissipativity of a nonlinear, input-affine system. Then, Theorem \ref{thm:KYPparticular} demonstrated the use of a particular L-K functional to derive delay-distribution-dependent LMI conditions for stochastic dissipativity of a linear system with bounded stochastic input delays. A numerical example shows that the resulting LMI can be used for robust plant analysis and controller design, and exploiting probabilistic delay information significantly reduces conservatism over deterministic methods.

%% file: Figures/PlotStochasticNyquist.tex
\definecolor{mycolorgray}{HTML}{B5B5B5}%
\definecolor{mycolor2}{HTML}{D55E00}%
\definecolor{mycolor1}{HTML}{0072B2}%
\definecolor{mycolor3}{HTML}{009E73}%
\definecolor{mycolor4}{HTML}{F0E442}%
\definecolor{mycolor5}{HTML}{CC79A7}%
\def \w {4}
\begin{tikzpicture}

\begin{axis}[%
width=1.8in,
height=1.8in,
at={(0in,0in)},
scale only axis,
bar shift auto,
xmin=-1,
xmax=2.5,
xlabel style={font=\color{white!15!black}},
xlabel={Real},
xtick={-1, -0.5, 0, 0.5, 1, 1.5, 2},
ymin=-1.75,
ymax=1.75,
ylabel style={font=\color{white!15!black}, align=center},
ylabel={Imaginary},
axis background/.style={fill=white},
%axis x line*=bottom,
%axis y line*=left,
legend style={legend cell align=left, align=left, draw=white!15!black},
legend pos = outer north east
]
\addplot[draw= black,ultra thick] table[row sep=crcr] {
	 -0.0473    0.0000		\\
	-0.0473    0.0019		\\
	-0.0474    0.0038		\\
	-0.0474    0.0058		\\
	-0.0474    0.0077		\\
	-0.0474    0.0097		\\
	-0.0475    0.0117		\\
	-0.0475    0.0138		\\
	-0.0476    0.0159		\\
	-0.0476    0.0181		\\
	-0.0477    0.0203		\\
	-0.0478    0.0227		\\
	-0.0479    0.0251		\\
	-0.0480    0.0277		\\
	-0.0481    0.0304		\\
	-0.0482    0.0333		\\
	-0.0484    0.0363		\\
	-0.0486    0.0396		\\
	-0.0488    0.0431		\\
	-0.0490    0.0469		\\
	-0.0496    0.0572		\\
	-0.0503    0.0687		\\
	-0.0511    0.0816		\\
	-0.0519    0.0963		\\
	-0.0525    0.1131		\\
	-0.0528    0.1324		\\
	-0.0525    0.1546		\\
	-0.0514    0.1801		\\
	-0.0490    0.2093		\\
	-0.0447    0.2426		\\
	-0.0379    0.2802		\\
	-0.0280    0.3225		\\
	-0.0140    0.3697		\\
	0.0050    0.4218		\\
	0.0300    0.4789		\\
	0.0625    0.5405		\\
	0.1037    0.6062		\\
	0.1551    0.6750		\\
	0.2182    0.7454		\\
	0.2942    0.8152		\\
	0.3838    0.8819		\\
	0.4870    0.9423		\\
	0.6029    0.9929		\\
	0.7293    1.0304		\\
	0.8631    1.0522		\\
	1.0000    1.0567		\\
	1.1356    1.0434		\\
	1.2656    1.0137		\\
	1.3863    0.9698		\\
	1.4951    0.9148		\\
	1.5907    0.8523		\\
	1.6728    0.7854		\\
	1.7418    0.7172		\\
	1.7989    0.6499		\\
	1.8455    0.5852		\\
	1.8831    0.5243		\\
	1.9131    0.4678		\\
	1.9370    0.4160		\\
	1.9558    0.3690		\\
	1.9705    0.3266		\\
	1.9821    0.2887		\\
	1.9911    0.2548		\\
	1.9981    0.2247		\\
	2.0035    0.1980		\\
	2.0078    0.1744		\\
	2.0110    0.1535		\\
	2.0136    0.1351		\\
	2.0155    0.1188		\\
	2.0171    0.1045		\\
	2.0182    0.0919		\\
	2.0191    0.0808		\\
	2.0198    0.0710		\\
	2.0204    0.0624		\\
	2.0208    0.0549		\\
	2.0211    0.0483		\\
	2.0214    0.0424		\\
	2.0222    0.0042		\\
	2.0222    0.0000		\\
	2.0222    0.0000		\\
	2.0222    0.0000		\\
	2.0222    0.0000		\\
	2.0222    0.0000		\\
	2.0222         0		\\
	NaN       NaN		\\
	2.0222         0		\\
	2.0222   -0.0000		\\
	2.0222   -0.0000		\\
	2.0222   -0.0000		\\
	2.0222   -0.0000		\\
	2.0222   -0.0000		\\
	2.0222   -0.0042		\\
	2.0214   -0.0424		\\
	2.0211   -0.0483		\\
	2.0208   -0.0549		\\
	2.0204   -0.0624		\\
	2.0198   -0.0710		\\
	2.0191   -0.0808		\\
	2.0182   -0.0919		\\
	2.0171   -0.1045		\\
	2.0155   -0.1188		\\
	2.0136   -0.1351		\\
	2.0110   -0.1535		\\
	2.0078   -0.1744		\\
	2.0035   -0.1980		\\
	1.9981   -0.2247		\\
	1.9911   -0.2548		\\
	1.9821   -0.2887		\\
	1.9705   -0.3266		\\
	1.9558   -0.3690		\\
	1.9370   -0.4160		\\
	1.9131   -0.4678		\\
	1.8831   -0.5243		\\
	1.8455   -0.5852		\\
	1.7989   -0.6499		\\
	1.7418   -0.7172		\\
	1.6728   -0.7854		\\
	1.5907   -0.8523		\\
	1.4951   -0.9148		\\
	1.3863   -0.9698		\\
	1.2656   -1.0137		\\
	1.1356   -1.0434		\\
	1.0000   -1.0567		\\
	0.8631   -1.0522		\\
	0.7293   -1.0304		\\
	0.6029   -0.9929		\\
	0.4870   -0.9423		\\
	0.3838   -0.8819		\\
	0.2942   -0.8152		\\
	0.2182   -0.7454		\\
	0.1551   -0.6750		\\
	0.1037   -0.6062		\\
	0.0625   -0.5405		\\
	0.0300   -0.4789		\\
	0.0050   -0.4218		\\
	-0.0140   -0.3697		\\
	-0.0280   -0.3225		\\
	-0.0379   -0.2802		\\
	-0.0447   -0.2426		\\
	-0.0490   -0.2093		\\
	-0.0514   -0.1801		\\
	-0.0525   -0.1546		\\
	-0.0528   -0.1324		\\
	-0.0525   -0.1131		\\
	-0.0519   -0.0963		\\
	-0.0511   -0.0816		\\
	-0.0503   -0.0687		\\
	-0.0496   -0.0572		\\
	-0.0490   -0.0469		\\
	-0.0488   -0.0431		\\
	-0.0486   -0.0396		\\
	-0.0484   -0.0363		\\
	-0.0482   -0.0333		\\
	-0.0481   -0.0304		\\
	-0.0480   -0.0277		\\
	-0.0479   -0.0251		\\
	-0.0478   -0.0227		\\
	-0.0477   -0.0203		\\
	-0.0476   -0.0181		\\
	-0.0476   -0.0159		\\
	-0.0475   -0.0138		\\
	-0.0475   -0.0117		\\
	-0.0474   -0.0097		\\
	-0.0474   -0.0077		\\
	-0.0474   -0.0058		\\
	-0.0474   -0.0038		\\
	-0.0473   -0.0019		\\
	-0.0473   -0.0000		\\
};
\addlegendentry{$w_k = 0$}

\addplot[draw=mycolor1,ultra thick,loosely dotted] table[row sep=crcr] {
 0.0473   -0.0000		\\
0.0470   -0.0058		\\
0.0461   -0.0116		\\
0.0445   -0.0172		\\
0.0423   -0.0227		\\
0.0395   -0.0279		\\
0.0362   -0.0329		\\
0.0322   -0.0375		\\
0.0287   -0.0409		\\
0.0278   -0.0417		\\
0.0228   -0.0455		\\
0.0207   -0.0469		\\
0.0173   -0.0488		\\
0.0116   -0.0516		\\
0.0115   -0.0516		\\
0.0052   -0.0538		\\
0.0014   -0.0548		\\
-0.0015   -0.0554		\\
-0.0085   -0.0562		\\
-0.0097   -0.0563		\\
-0.0158   -0.0564		\\
-0.0214   -0.0560		\\
-0.0233   -0.0558		\\
-0.0310   -0.0544		\\
-0.0337   -0.0537		\\
-0.0389   -0.0521		\\
-0.0464   -0.0492		\\
-0.0469   -0.0490		\\
-0.0593   -0.0423		\\
-0.0656   -0.0379		\\
-0.0722   -0.0326		\\
-0.0820   -0.0231		\\
-0.0850   -0.0198		\\
-0.0962   -0.0053		\\
-0.0976   -0.0031		\\
-0.1083    0.0152		\\
-0.1099    0.0185		\\
-0.1186    0.0383		\\
-0.1218    0.0469		\\
-0.1272    0.0642		\\
-0.1326    0.0857		\\
-0.1341    0.0931		\\
-0.1391    0.1254		\\
-0.1407    0.1423		\\
-0.1419    0.1614		\\
-0.1421    0.2016		\\
-0.1402    0.2328		\\
-0.1391    0.2462		\\
-0.1321    0.2956		\\
-0.1203    0.3499		\\
-0.1056    0.3999		\\
-0.1027    0.4092		\\
-0.0783    0.4734		\\
-0.0456    0.5422		\\
-0.0034    0.6150		\\
0.0499    0.6908		\\
0.1158    0.7680		\\
0.1327    0.7844		\\
0.1955    0.8443		\\
0.2898    0.9171		\\
0.3986    0.9829		\\
0.5210    1.0382		\\
0.6545    1.0795		\\
0.7959    1.1039		\\
0.9407    1.1098		\\
1.0841    1.0968		\\
1.2216    1.0663		\\
1.3493    1.0206		\\
1.4645    0.9632		\\
1.5656    0.8976		\\
1.6524    0.8274		\\
1.7255    0.7557		\\
1.7859    0.6848		\\
1.8352    0.6167		\\
1.8750    0.5526		\\
1.9068    0.4930		\\
1.9320    0.4385		\\
1.9519    0.3889		\\
1.9675    0.3443		\\
1.9797    0.3043		\\
1.9893    0.2686		\\
1.9967    0.2369		\\
2.0024    0.2087		\\
2.0069    0.1838		\\
2.0104    0.1618		\\
2.0131    0.1424		\\
2.0151    0.1253		\\
2.0167    0.1102		\\
2.0180    0.0969		\\
2.0190    0.0852		\\
2.0197    0.0749		\\
2.0203    0.0658		\\
2.0207    0.0579		\\
2.0211    0.0509		\\
2.0213    0.0447		\\
2.0222    0.0045		\\
2.0222    0.0000		\\
2.0222    0.0000		\\
2.0222    0.0000		\\
2.0222    0.0000		\\
2.0222    0.0000		\\
2.0222         0		\\
2.0222         0		\\
2.0222   -0.0000		\\
2.0222   -0.0000		\\
2.0222   -0.0000		\\
2.0222   -0.0000		\\
2.0222   -0.0000		\\
2.0222   -0.0045		\\
2.0213   -0.0447		\\
2.0211   -0.0509		\\
2.0207   -0.0579		\\
2.0203   -0.0658		\\
2.0197   -0.0749		\\
2.0190   -0.0852		\\
2.0180   -0.0969		\\
2.0167   -0.1102		\\
2.0151   -0.1253		\\
2.0131   -0.1424		\\
2.0104   -0.1618		\\
2.0069   -0.1838		\\
2.0024   -0.2087		\\
1.9967   -0.2369		\\
1.9893   -0.2686		\\
1.9797   -0.3043		\\
1.9675   -0.3443		\\
1.9519   -0.3889		\\
1.9320   -0.4385		\\
1.9068   -0.4930		\\
1.8750   -0.5526		\\
1.8352   -0.6167		\\
1.7859   -0.6848		\\
1.7255   -0.7557		\\
1.6524   -0.8274		\\
1.5656   -0.8976		\\
1.4645   -0.9632		\\
1.3493   -1.0206		\\
1.2216   -1.0663		\\
1.0841   -1.0968		\\
0.9407   -1.1098		\\
0.7959   -1.1039		\\
0.6545   -1.0795		\\
0.5210   -1.0382		\\
0.3986   -0.9829		\\
0.2898   -0.9171		\\
0.1955   -0.8443		\\
0.1327   -0.7844		\\
0.1158   -0.7680		\\
0.0499   -0.6908		\\
-0.0034   -0.6150		\\
-0.0456   -0.5422		\\
-0.0783   -0.4734		\\
-0.1027   -0.4092		\\
-0.1056   -0.3999		\\
-0.1203   -0.3499		\\
-0.1321   -0.2956		\\
-0.1391   -0.2462		\\
-0.1402   -0.2328		\\
-0.1421   -0.2016		\\
-0.1419   -0.1614		\\
-0.1407   -0.1423		\\
-0.1391   -0.1254		\\
-0.1341   -0.0931		\\
-0.1326   -0.0857		\\
-0.1272   -0.0642		\\
-0.1218   -0.0469		\\
-0.1186   -0.0383		\\
-0.1099   -0.0185		\\
-0.1083   -0.0152		\\
-0.0976    0.0031		\\
-0.0962    0.0053		\\
-0.0850    0.0198		\\
-0.0820    0.0231		\\
-0.0722    0.0326		\\
-0.0656    0.0379		\\
-0.0593    0.0423		\\
-0.0469    0.0490		\\
-0.0464    0.0492		\\
-0.0389    0.0521		\\
-0.0337    0.0537		\\
-0.0310    0.0544		\\
-0.0233    0.0558		\\
-0.0214    0.0560		\\
-0.0158    0.0564		\\
-0.0097    0.0563		\\
-0.0085    0.0562		\\
-0.0015    0.0554		\\
0.0014    0.0548		\\
0.0052    0.0538		\\
0.0115    0.0516		\\
0.0116    0.0516		\\
0.0173    0.0488		\\
0.0207    0.0469		\\
0.0228    0.0455		\\
0.0278    0.0417		\\
0.0287    0.0409		\\
0.0322    0.0375		\\
0.0362    0.0329		\\
0.0395    0.0279		\\
0.0423    0.0227		\\
0.0445    0.0172		\\
0.0461    0.0116		\\
0.0470    0.0058		\\
0.0473    0.0000		\\
};
\addlegendentry{$w_k = 1$}

\addplot[draw=mycolor2,ultra thick,dashdotted] table[row sep=crcr] {
	 -0.0473    0.0000		\\
	-0.0464    0.0097		\\
	-0.0435    0.0190		\\
	-0.0389    0.0276		\\
	-0.0327    0.0352		\\
	-0.0250    0.0414		\\
	-0.0162    0.0461		\\
	-0.0065    0.0490		\\
	0.0017    0.0500		\\
	0.0037    0.0500		\\
	0.0097    0.0496		\\
	0.0141    0.0489		\\
	0.0177    0.0481		\\
	0.0242    0.0458		\\
	0.0254    0.0453		\\
	0.0328    0.0413		\\
	0.0337    0.0407		\\
	0.0396    0.0362		\\
	0.0422    0.0337		\\
	0.0456    0.0300		\\
	0.0494    0.0250		\\
	0.0508    0.0229		\\
	0.0549    0.0148		\\
	0.0549    0.0148		\\
	0.0579    0.0061		\\
	0.0585    0.0034		\\
	0.0595   -0.0032		\\
	0.0598   -0.0089		\\
	0.0598   -0.0130		\\
	0.0588   -0.0217		\\
	0.0585   -0.0229		\\
	0.0558   -0.0328		\\
	0.0552   -0.0345		\\
	0.0515   -0.0426		\\
	0.0490   -0.0469		\\
	0.0456   -0.0519		\\
	0.0382   -0.0607		\\
	0.0292   -0.0685		\\
	0.0248   -0.0716		\\
	0.0187   -0.0753		\\
	0.0068   -0.0808		\\
	-0.0064   -0.0847		\\
	-0.0070   -0.0849		\\
	-0.0209   -0.0870		\\
	-0.0364   -0.0872		\\
	-0.0418   -0.0867		\\
	-0.0528   -0.0852		\\
	-0.0701   -0.0806		\\
	-0.0764   -0.0783		\\
	-0.0880   -0.0732		\\
	-0.1063   -0.0626		\\
	-0.1087   -0.0610		\\
	-0.1250   -0.0484		\\
	-0.1379   -0.0362		\\
	-0.1437   -0.0299		\\
	-0.1621   -0.0064		\\
	-0.1632   -0.0048		\\
	-0.1799    0.0235		\\
	-0.1844    0.0325		\\
	-0.1964    0.0612		\\
	-0.2013    0.0753		\\
	-0.2104    0.1095		\\
	-0.2134    0.1236		\\
	-0.2199    0.1722		\\
	-0.2204    0.1772		\\
	-0.2216    0.2360		\\
	-0.2204    0.2563		\\
	-0.2163    0.3002		\\
	-0.2036    0.3695		\\
	-0.2027    0.3732		\\
	-0.1826    0.4437		\\
	-0.1520    0.5225		\\
	-0.1410    0.5451		\\
	-0.1104    0.6050		\\
	-0.0565    0.6903		\\
	0.0113    0.7766		\\
	0.0452    0.8120		\\
	0.0941    0.8616		\\
	0.1926    0.9423		\\
	0.3069    1.0153		\\
	0.4356    1.0768		\\
	0.5765    1.1231		\\
	0.6941    1.1457		\\
	0.7257    1.1513		\\
	0.8786    1.1596		\\
	1.0301    1.1477		\\
	1.1754    1.1170		\\
	1.3104    1.0700		\\
	1.4322    1.0105		\\
	1.5392    0.9421		\\
	1.6310    0.8688		\\
	1.7083    0.7937		\\
	1.7722    0.7195		\\
	1.8244    0.6481		\\
	1.8665    0.5807		\\
	1.9001    0.5182		\\
	1.9268    0.4609		\\
	1.9478    0.4089		\\
	1.9643    0.3620		\\
	1.9773    0.3199		\\
	1.9873    0.2824		\\
	1.9952    0.2490		\\
	2.0013    0.2195		\\
	2.0060    0.1933		\\
	2.0097    0.1701		\\
	2.0125    0.1497		\\
	2.0147    0.1317		\\
	2.0164    0.1158		\\
	2.0177    0.1019		\\
	2.0188    0.0896		\\
	2.0196    0.0787		\\
	2.0202    0.0692		\\
	2.0206    0.0608		\\
	2.0210    0.0535		\\
	2.0213    0.0470		\\
	2.0222    0.0047		\\
	2.0222    0.0000		\\
	2.0222    0.0000		\\
	2.0222    0.0000		\\
	2.0222    0.0000		\\
	2.0222    0.0000		\\
	2.0222         0		\\
	2.0222         0		\\
	2.0222   -0.0000		\\
	2.0222   -0.0000		\\
	2.0222   -0.0000		\\
	2.0222   -0.0000		\\
	2.0222   -0.0000		\\
	2.0222   -0.0047		\\
	2.0213   -0.0470		\\
	2.0210   -0.0535		\\
	2.0206   -0.0608		\\
	2.0202   -0.0692		\\
	2.0196   -0.0787		\\
	2.0188   -0.0896		\\
	2.0177   -0.1019		\\
	2.0164   -0.1158		\\
	2.0147   -0.1317		\\
	2.0125   -0.1497		\\
	2.0097   -0.1701		\\
	2.0060   -0.1933		\\
	2.0013   -0.2195		\\
	1.9952   -0.2490		\\
	1.9873   -0.2824		\\
	1.9773   -0.3199		\\
	1.9643   -0.3620		\\
	1.9478   -0.4089		\\
	1.9268   -0.4609		\\
	1.9001   -0.5182		\\
	1.8665   -0.5807		\\
	1.8244   -0.6481		\\
	1.7722   -0.7195		\\
	1.7083   -0.7937		\\
	1.6310   -0.8688		\\
	1.5392   -0.9421		\\
	1.4322   -1.0105		\\
	1.3104   -1.0700		\\
	1.1754   -1.1170		\\
	1.0301   -1.1477		\\
	0.8786   -1.1596		\\
	0.7257   -1.1513		\\
	0.6941   -1.1457		\\
	0.5765   -1.1231		\\
	0.4356   -1.0768		\\
	0.3069   -1.0153		\\
	0.1926   -0.9423		\\
	0.0941   -0.8616		\\
	0.0452   -0.8120		\\
	0.0113   -0.7766		\\
	-0.0565   -0.6903		\\
	-0.1104   -0.6050		\\
	-0.1410   -0.5451		\\
	-0.1520   -0.5225		\\
	-0.1826   -0.4437		\\
	-0.2027   -0.3732		\\
	-0.2036   -0.3695		\\
	-0.2163   -0.3002		\\
	-0.2204   -0.2563		\\
	-0.2216   -0.2360		\\
	-0.2204   -0.1772		\\
	-0.2199   -0.1722		\\
	-0.2134   -0.1236		\\
	-0.2104   -0.1095		\\
	-0.2013   -0.0753		\\
	-0.1964   -0.0612		\\
	-0.1844   -0.0325		\\
	-0.1799   -0.0235		\\
	-0.1632    0.0048		\\
	-0.1621    0.0064		\\
	-0.1437    0.0299		\\
	-0.1379    0.0362		\\
	-0.1250    0.0484		\\
	-0.1087    0.0610		\\
	-0.1063    0.0626		\\
	-0.0880    0.0732		\\
	-0.0764    0.0783		\\
	-0.0701    0.0806		\\
	-0.0528    0.0852		\\
	-0.0418    0.0867		\\
	-0.0364    0.0872		\\
	-0.0209    0.0870		\\
	-0.0070    0.0849		\\
	-0.0064    0.0847		\\
	0.0068    0.0808		\\
	0.0187    0.0753		\\
	0.0248    0.0716		\\
	0.0292    0.0685		\\
	0.0382    0.0607		\\
	0.0456    0.0519		\\
	0.0490    0.0469		\\
	0.0515    0.0426		\\
	0.0552    0.0345		\\
	0.0558    0.0328		\\
	0.0585    0.0229		\\
	0.0588    0.0217		\\
	0.0598    0.0130		\\
	0.0598    0.0089		\\
	0.0595    0.0032		\\
	0.0585   -0.0034		\\
	0.0579   -0.0061		\\
	0.0549   -0.0148		\\
	0.0549   -0.0148		\\
	0.0508   -0.0229		\\
	0.0494   -0.0250		\\
	0.0456   -0.0300		\\
	0.0422   -0.0337		\\
	0.0396   -0.0362		\\
	0.0337   -0.0407		\\
	0.0328   -0.0413		\\
	0.0254   -0.0453		\\
	0.0242   -0.0458		\\
	0.0177   -0.0481		\\
	0.0141   -0.0489		\\
	0.0097   -0.0496		\\
	0.0037   -0.0500		\\
	0.0017   -0.0500		\\
	-0.0065   -0.0490		\\
	-0.0162   -0.0461		\\
	-0.0250   -0.0414		\\
	-0.0327   -0.0352		\\
	-0.0389   -0.0276		\\
	-0.0435   -0.0190		\\
	-0.0464   -0.0097		\\
	-0.0473   -0.0000		\\
};
\addlegendentry{$w_k = 2$}

\addplot[draw=mycolor3,ultra thick,dotted] table[row sep=crcr] {
0.0473   -0.0000		\\
0.0454   -0.0135		\\
0.0398   -0.0259		\\
0.0310   -0.0363		\\
0.0195   -0.0439		\\
0.0062   -0.0480		\\
-0.0077   -0.0483		\\
-0.0214   -0.0446		\\
-0.0314   -0.0389		\\
-0.0336   -0.0372		\\
-0.0371   -0.0341		\\
-0.0420   -0.0285		\\
-0.0435   -0.0265		\\
-0.0461   -0.0222		\\
-0.0494   -0.0154		\\
-0.0501   -0.0133		\\
-0.0516   -0.0080		\\
-0.0527   -0.0003		\\
-0.0528    0.0016		\\
-0.0528    0.0075		\\
-0.0517    0.0153		\\
-0.0513    0.0169		\\
-0.0495    0.0231		\\
-0.0461    0.0305		\\
-0.0455    0.0315		\\
-0.0417    0.0374		\\
-0.0363    0.0437		\\
-0.0356    0.0444		\\
-0.0299    0.0493		\\
-0.0226    0.0540		\\
-0.0222    0.0542		\\
-0.0147    0.0576		\\
-0.0062    0.0602		\\
-0.0061    0.0602		\\
0.0028    0.0615		\\
0.0117    0.0615		\\
0.0120    0.0615		\\
0.0213    0.0602		\\
0.0298    0.0578		\\
0.0305    0.0576		\\
0.0395    0.0536		\\
0.0469    0.0490		\\
0.0480    0.0482		\\
0.0558    0.0416		\\
0.0628    0.0338		\\
0.0688    0.0248		\\
0.0736    0.0148		\\
0.0750    0.0108		\\
0.0770    0.0039		\\
0.0790   -0.0077		\\
0.0793   -0.0198		\\
0.0779   -0.0322		\\
0.0771   -0.0363		\\
0.0747   -0.0448		\\
0.0697   -0.0573		\\
0.0627   -0.0694		\\
0.0557   -0.0785		\\
0.0538   -0.0808		\\
0.0429   -0.0915		\\
0.0301   -0.1010		\\
0.0179   -0.1079		\\
0.0155   -0.1090		\\
-0.0009   -0.1154		\\
-0.0189   -0.1198		\\
-0.0292   -0.1212		\\
-0.0385   -0.1219		\\
-0.0596   -0.1214		\\
-0.0793   -0.1184		\\
-0.0818   -0.1179		\\
-0.1052   -0.1110		\\
-0.1284   -0.1008		\\
-0.1293   -0.1004		\\
-0.1541   -0.0853		\\
-0.1736   -0.0703		\\
-0.1790   -0.0655		\\
-0.2037   -0.0400		\\
-0.2130   -0.0286		\\
-0.2276   -0.0078		\\
-0.2456    0.0229		\\
-0.2501    0.0320		\\
-0.2699    0.0814		\\
-0.2704    0.0829		\\
-0.2852    0.1426		\\
-0.2867    0.1505		\\
-0.2935    0.2189		\\
-0.2938    0.2249		\\
-0.2911    0.3054		\\
-0.2899    0.3152		\\
-0.2776    0.3914		\\
-0.2652    0.4389		\\
-0.2523    0.4821		\\
-0.2141    0.5766		\\
-0.2013    0.6013		\\
-0.1616    0.6735		\\
-0.0935    0.7710		\\
-0.0533    0.8174		\\
-0.0087    0.8666		\\
0.0934    0.9573		\\
0.2126    1.0392		\\
0.2965    1.0830		\\
0.3475    1.1084		\\
0.4955    1.1611		\\
0.6526    1.1943		\\
0.8138    1.2059		\\
0.9737    1.1959		\\
1.1265    1.1658		\\
1.1272    1.1657		\\
1.2698    1.1180		\\
1.3985    1.0567		\\
1.5115    0.9859		\\
1.6086    0.9097		\\
1.6903    0.8314		\\
1.7579    0.7539		\\
1.8130    0.6792		\\
1.8575    0.6088		\\
1.8931    0.5433		\\
1.9213    0.4833		\\
1.9435    0.4288		\\
1.9610    0.3796		\\
1.9747    0.3355		\\
1.9853    0.2962		\\
1.9936    0.2612		\\
2.0001    0.2302		\\
2.0051    0.2027		\\
2.0090    0.1784		\\
2.0120    0.1570		\\
2.0143    0.1381		\\
2.0161    0.1215		\\
2.0175    0.1068		\\
2.0186    0.0939		\\
2.0194    0.0826		\\
2.0200    0.0726		\\
2.0205    0.0638		\\
2.0209    0.0561		\\
2.0212    0.0493		\\
2.0222    0.0049		\\
2.0222    0.0000		\\
2.0222    0.0000		\\
2.0222    0.0000		\\
2.0222    0.0000		\\
2.0222    0.0000		\\
2.0222         0		\\
2.0222         0		\\
2.0222   -0.0000		\\
2.0222   -0.0000		\\
2.0222   -0.0000		\\
2.0222   -0.0000		\\
2.0222   -0.0000		\\
2.0222   -0.0049		\\
2.0212   -0.0493		\\
2.0209   -0.0561		\\
2.0205   -0.0638		\\
2.0200   -0.0726		\\
2.0194   -0.0826		\\
2.0186   -0.0939		\\
2.0175   -0.1068		\\
2.0161   -0.1215		\\
2.0143   -0.1381		\\
2.0120   -0.1570		\\
2.0090   -0.1784		\\
2.0051   -0.2027		\\
2.0001   -0.2302		\\
1.9936   -0.2612		\\
1.9853   -0.2962		\\
1.9747   -0.3355		\\
1.9610   -0.3796		\\
1.9435   -0.4288		\\
1.9213   -0.4833		\\
1.8931   -0.5433		\\
1.8575   -0.6088		\\
1.8130   -0.6792		\\
1.7579   -0.7539		\\
1.6903   -0.8314		\\
1.6086   -0.9097		\\
1.5115   -0.9859		\\
1.3985   -1.0567		\\
1.2698   -1.1180		\\
1.1272   -1.1657		\\
1.1265   -1.1658		\\
0.9737   -1.1959		\\
0.8138   -1.2059		\\
0.6526   -1.1943		\\
0.4955   -1.1611		\\
0.3475   -1.1084		\\
0.2965   -1.0830		\\
0.2126   -1.0392		\\
0.0934   -0.9573		\\
-0.0087   -0.8666		\\
-0.0533   -0.8174		\\
-0.0935   -0.7710		\\
-0.1616   -0.6735		\\
-0.2013   -0.6013		\\
-0.2141   -0.5766		\\
-0.2523   -0.4821		\\
-0.2652   -0.4389		\\
-0.2776   -0.3914		\\
-0.2899   -0.3152		\\
-0.2911   -0.3054		\\
-0.2938   -0.2249		\\
-0.2935   -0.2189		\\
-0.2867   -0.1505		\\
-0.2852   -0.1426		\\
-0.2704   -0.0829		\\
-0.2699   -0.0814		\\
-0.2501   -0.0320		\\
-0.2456   -0.0229		\\
-0.2276    0.0078		\\
-0.2130    0.0286		\\
-0.2037    0.0400		\\
-0.1790    0.0655		\\
-0.1736    0.0703		\\
-0.1541    0.0853		\\
-0.1293    0.1004		\\
-0.1284    0.1008		\\
-0.1052    0.1110		\\
-0.0818    0.1179		\\
-0.0793    0.1184		\\
-0.0596    0.1214		\\
-0.0385    0.1219		\\
-0.0292    0.1212		\\
-0.0189    0.1198		\\
-0.0009    0.1154		\\
0.0155    0.1090		\\
0.0179    0.1079		\\
0.0301    0.1010		\\
0.0429    0.0915		\\
0.0538    0.0808		\\
0.0557    0.0785		\\
0.0627    0.0694		\\
0.0697    0.0573		\\
0.0747    0.0448		\\
0.0771    0.0363		\\
0.0779    0.0322		\\
0.0793    0.0198		\\
0.0790    0.0077		\\
0.0770   -0.0039		\\
0.0750   -0.0108		\\
0.0736   -0.0148		\\
0.0688   -0.0248		\\
0.0628   -0.0338		\\
0.0558   -0.0416		\\
0.0480   -0.0482		\\
0.0469   -0.0490		\\
0.0395   -0.0536		\\
0.0305   -0.0576		\\
0.0298   -0.0578		\\
0.0213   -0.0602		\\
0.0120   -0.0615		\\
0.0117   -0.0615		\\
0.0028   -0.0615		\\
-0.0061   -0.0602		\\
-0.0062   -0.0602		\\
-0.0147   -0.0576		\\
-0.0222   -0.0542		\\
-0.0226   -0.0540		\\
-0.0299   -0.0493		\\
-0.0356   -0.0444		\\
-0.0363   -0.0437		\\
-0.0417   -0.0374		\\
-0.0455   -0.0315		\\
-0.0461   -0.0305		\\
-0.0495   -0.0231		\\
-0.0513   -0.0169		\\
-0.0517   -0.0153		\\
-0.0528   -0.0075		\\
-0.0528   -0.0016		\\
-0.0527    0.0003		\\
-0.0516    0.0080		\\
-0.0501    0.0133		\\
-0.0494    0.0154		\\
-0.0461    0.0222		\\
-0.0435    0.0265		\\
-0.0420    0.0285		\\
-0.0371    0.0341		\\
-0.0336    0.0372		\\
-0.0314    0.0389		\\
-0.0214    0.0446		\\
-0.0077    0.0483		\\
0.0062    0.0480		\\
0.0195    0.0439		\\
0.0310    0.0363		\\
0.0398    0.0259		\\
0.0454    0.0135		\\
0.0473    0.0000		\\
};
\addlegendentry{$w_k = 3$}

\addplot[draw=mycolor4,ultra thick, densely dashed] table[row sep=crcr] {
	-0.0473    0.0000	\\
	-0.0442    0.0172	\\
	-0.0350    0.0321	\\
	-0.0211    0.0428	\\
	-0.0042    0.0478	\\
	0.0136    0.0465	\\
	0.0298    0.0388	\\
	0.0423    0.0256	\\
	0.0485    0.0122	\\
	0.0494    0.0087	\\
	0.0500    0.0052	\\
	0.0505   -0.0021	\\
	0.0500   -0.0093	\\
	0.0499   -0.0100	\\
	0.0485   -0.0165	\\
	0.0460   -0.0234	\\
	0.0437   -0.0279	\\
	0.0424   -0.0300	\\
	0.0380   -0.0360	\\
	0.0327   -0.0414	\\
	0.0312   -0.0427	\\
	0.0267   -0.0460	\\
	0.0200   -0.0498	\\
	0.0139   -0.0522	\\
	0.0128   -0.0526	\\
	0.0051   -0.0544	\\
	-0.0027   -0.0551	\\
	-0.0061   -0.0550	\\
	-0.0107   -0.0547	\\
	-0.0186   -0.0531	\\
	-0.0263   -0.0505	\\
	-0.0263   -0.0504	\\
	-0.0336   -0.0467	\\
	-0.0404   -0.0419	\\
	-0.0441   -0.0386	\\
	-0.0465   -0.0361	\\
	-0.0518   -0.0294	\\
	-0.0562   -0.0220	\\
	-0.0569   -0.0206	\\
	-0.0596   -0.0139	\\
	-0.0618   -0.0053	\\
	-0.0626    0.0014	\\
	-0.0628    0.0036	\\
	-0.0625    0.0128	\\
	-0.0609    0.0219	\\
	-0.0601    0.0249	\\
	-0.0581    0.0309	\\
	-0.0539    0.0395	\\
	-0.0490    0.0469	\\
	-0.0485    0.0476	\\
	-0.0419    0.0550	\\
	-0.0342    0.0615	\\
	-0.0256    0.0670	\\
	-0.0160    0.0713	\\
	-0.0058    0.0743	\\
	0.0036    0.0757	\\
	0.0050    0.0758	\\
	0.0162    0.0759	\\
	0.0275    0.0744	\\
	0.0387    0.0713	\\
	0.0496    0.0665	\\
	0.0601    0.0601	\\
	0.0610    0.0595	\\
	0.0698    0.0522	\\
	0.0786    0.0428	\\
	0.0862    0.0319	\\
	0.0924    0.0197	\\
	0.0957    0.0108	\\
	0.0970    0.0065	\\
	0.0999   -0.0078	\\
	0.1009   -0.0228	\\
	0.0997   -0.0382	\\
	0.0975   -0.0495	\\
	0.0964   -0.0540	\\
	0.0907   -0.0698	\\
	0.0827   -0.0852	\\
	0.0723   -0.1000	\\
	0.0691   -0.1038	\\
	0.0594   -0.1139	\\
	0.0440   -0.1266	\\
	0.0263   -0.1376	\\
	0.0190   -0.1413	\\
	0.0063   -0.1466	\\
	-0.0160   -0.1533	\\
	-0.0403   -0.1572	\\
	-0.0434   -0.1574	\\
	-0.0667   -0.1578	\\
	-0.0946   -0.1549	\\
	-0.1098   -0.1517	\\
	-0.1239   -0.1478	\\
	-0.1544   -0.1360	\\
	-0.1743   -0.1257	\\
	-0.1855   -0.1190	\\
	-0.2169   -0.0960	\\
	-0.2326   -0.0820	\\
	-0.2479   -0.0663	\\
	-0.2778   -0.0290	\\
	-0.2818   -0.0232	\\
	-0.3055    0.0174	\\
	-0.3201    0.0483	\\
	-0.3298    0.0739	\\
	-0.3463    0.1303	\\
	-0.3488    0.1426	\\
	-0.3593    0.2211	\\
	-0.3595    0.2260	\\
	-0.3584    0.3190	\\
	-0.3575    0.3277	\\
	-0.3427    0.4227	\\
	-0.3350    0.4526	\\
	-0.3112    0.5305	\\
	-0.2789    0.6066	\\
	-0.2628    0.6408	\\
	-0.1965    0.7514	\\
	-0.1619    0.7970	\\
	-0.1113    0.8595	\\
	-0.0069    0.9618	\\
	0.0722    1.0222	\\
	0.1165    1.0543	\\
	0.2571    1.1328	\\
	0.4120    1.1933	\\
	0.5421    1.2249	\\
	0.5770    1.2326	\\
	0.7466    1.2487	\\
	0.9151    1.2414	\\
	1.0769    1.2123	\\
	1.2273    1.1644	\\
	1.3632    1.1018	\\
	1.3828    1.0902	\\
	1.4826    1.0289	\\
	1.5851    0.9500	\\
	1.6714    0.8687	\\
	1.7429    0.7880	\\
	1.8011    0.7102	\\
	1.8481    0.6367	\\
	1.8857    0.5683	\\
	1.9155    0.5056	\\
	1.9391    0.4486	\\
	1.9575    0.3972	\\
	1.9720    0.3511	\\
	1.9832    0.3100	\\
	1.9920    0.2734	\\
	1.9988    0.2409	\\
	2.0041    0.2122	\\
	2.0082    0.1868	\\
	2.0114    0.1643	\\
	2.0138    0.1446	\\
	2.0157    0.1272	\\
	2.0172    0.1118	\\
	2.0184    0.0983	\\
	2.0192    0.0864	\\
	2.0199    0.0760	\\
	2.0204    0.0668	\\
	2.0208    0.0587	\\
	2.0212    0.0516	\\
	2.0222    0.0052	\\
	2.0222    0.0001	\\
	2.0222    0.0000	\\
	2.0222    0.0000	\\
	2.0222    0.0000	\\
	2.0222    0.0000	\\
	2.0222         0	\\
	NaN       NaN	\\
	2.0222         0	\\
	2.0222   -0.0000	\\
	2.0222   -0.0000	\\
	2.0222   -0.0000	\\
	2.0222   -0.0000	\\
	2.0222   -0.0001	\\
	2.0222   -0.0052	\\
	2.0212   -0.0516	\\
	2.0208   -0.0587	\\
	2.0204   -0.0668	\\
	2.0199   -0.0760	\\
	2.0192   -0.0864	\\
	2.0184   -0.0983	\\
	2.0172   -0.1118	\\
	2.0157   -0.1272	\\
	2.0138   -0.1446	\\
	2.0114   -0.1643	\\
	2.0082   -0.1868	\\
	2.0041   -0.2122	\\
	1.9988   -0.2409	\\
	1.9920   -0.2734	\\
	1.9832   -0.3100	\\
	1.9720   -0.3511	\\
	1.9575   -0.3972	\\
	1.9391   -0.4486	\\
	1.9155   -0.5056	\\
	1.8857   -0.5683	\\
	1.8481   -0.6367	\\
	1.8011   -0.7102	\\
	1.7429   -0.7880	\\
	1.6714   -0.8687	\\
	1.5851   -0.9500	\\
	1.4826   -1.0289	\\
	1.3828   -1.0902	\\
	1.3632   -1.1018	\\
	1.2273   -1.1644	\\
	1.0769   -1.2123	\\
	0.9151   -1.2414	\\
	0.7466   -1.2487	\\
	0.5770   -1.2326	\\
	0.5421   -1.2249	\\
	0.4120   -1.1933	\\
	0.2571   -1.1328	\\
	0.1165   -1.0543	\\
	0.0722   -1.0222	\\
	-0.0069   -0.9618	\\
	-0.1113   -0.8595	\\
	-0.1619   -0.7970	\\
	-0.1965   -0.7514	\\
	-0.2628   -0.6408	\\
	-0.2789   -0.6066	\\
	-0.3112   -0.5305	\\
	-0.3350   -0.4526	\\
	-0.3427   -0.4227	\\
	-0.3575   -0.3277	\\
	-0.3584   -0.3190	\\
	-0.3595   -0.2260	\\
	-0.3593   -0.2211	\\
	-0.3488   -0.1426	\\
	-0.3463   -0.1303	\\
	-0.3298   -0.0739	\\
	-0.3201   -0.0483	\\
	-0.3055   -0.0174	\\
	-0.2818    0.0232	\\
	-0.2778    0.0290	\\
	-0.2479    0.0663	\\
	-0.2326    0.0820	\\
	-0.2169    0.0960	\\
	-0.1855    0.1190	\\
	-0.1743    0.1257	\\
	-0.1544    0.1360	\\
	-0.1239    0.1478	\\
	-0.1098    0.1517	\\
	-0.0946    0.1549	\\
	-0.0667    0.1578	\\
	-0.0434    0.1574	\\
	-0.0403    0.1572	\\
	-0.0160    0.1533	\\
	0.0063    0.1466	\\
	0.0190    0.1413	\\
	0.0263    0.1376	\\
	0.0440    0.1266	\\
	0.0594    0.1139	\\
	0.0691    0.1038	\\
	0.0723    0.1000	\\
	0.0827    0.0852	\\
	0.0907    0.0698	\\
	0.0964    0.0540	\\
	0.0975    0.0495	\\
	0.0997    0.0382	\\
	0.1009    0.0228	\\
	0.0999    0.0078	\\
	0.0970   -0.0065	\\
	0.0957   -0.0108	\\
	0.0924   -0.0197	\\
	0.0862   -0.0319	\\
	0.0786   -0.0428	\\
	0.0698   -0.0522	\\
	0.0610   -0.0595	\\
	0.0601   -0.0601	\\
	0.0496   -0.0665	\\
	0.0387   -0.0713	\\
	0.0275   -0.0744	\\
	0.0162   -0.0759	\\
	0.0050   -0.0758	\\
	0.0036   -0.0757	\\
	-0.0058   -0.0743	\\
	-0.0160   -0.0713	\\
	-0.0256   -0.0670	\\
	-0.0342   -0.0615	\\
	-0.0419   -0.0550	\\
	-0.0485   -0.0476	\\
	-0.0490   -0.0469	\\
	-0.0539   -0.0395	\\
	-0.0581   -0.0309	\\
	-0.0601   -0.0249	\\
	-0.0609   -0.0219	\\
	-0.0625   -0.0128	\\
	-0.0628   -0.0036	\\
	-0.0626   -0.0014	\\
	-0.0618    0.0053	\\
	-0.0596    0.0139	\\
	-0.0569    0.0206	\\
	-0.0562    0.0220	\\
	-0.0518    0.0294	\\
	-0.0465    0.0361	\\
	-0.0441    0.0386	\\
	-0.0404    0.0419	\\
	-0.0336    0.0467	\\
	-0.0263    0.0504	\\
	-0.0263    0.0505	\\
	-0.0186    0.0531	\\
	-0.0107    0.0547	\\
	-0.0061    0.0550	\\
	-0.0027    0.0551	\\
	0.0051    0.0544	\\
	0.0128    0.0526	\\
	0.0139    0.0522	\\
	0.0200    0.0498	\\
	0.0267    0.0460	\\
	0.0312    0.0427	\\
	0.0327    0.0414	\\
	0.0380    0.0360	\\
	0.0424    0.0300	\\
	0.0437    0.0279	\\
	0.0460    0.0234	\\
	0.0485    0.0165	\\
	0.0499    0.0100	\\
	0.0500    0.0093	\\
	0.0505    0.0021	\\
	0.0500   -0.0052	\\
	0.0494   -0.0087	\\
	0.0485   -0.0122	\\
	0.0423   -0.0256	\\
	0.0298   -0.0388	\\
	0.0136   -0.0465	\\
	-0.0042   -0.0478	\\
	-0.0211   -0.0428	\\
	-0.0350   -0.0321	\\
	-0.0442   -0.0172	\\
	-0.0473   -0.0000	\\
};
\addlegendentry{$w_k = 4$}

\addplot[draw=mycolor5,ultra thick, loosely dashed] table[row sep=crcr] {
0.0473	0	\\
0.0426	-0.0208	\\
0.0293	-0.0374	\\
0.0099	-0.0467	\\
-0.0116	-0.0466	\\
-0.0311	-0.0371	\\
-0.0446	-0.0199	\\
-0.0494	0.0017	\\
-0.0461	0.0194	\\
-0.0443	0.0235	\\
-0.0431	0.0258	\\
-0.0393	0.0317	\\
-0.0347	0.037	\\
-0.03	0.0412	\\
-0.0293	0.0417	\\
-0.0234	0.0456	\\
-0.0169	0.0487	\\
-0.01	0.0508	\\
-0.0091	0.051	\\
-0.0028	0.052	\\
0.0045	0.0522	\\
0.0118	0.0514	\\
0.0145	0.0508	\\
0.019	0.0496	\\
0.0259	0.0468	\\
0.0323	0.043	\\
0.0361	0.0402	\\
0.0383	0.0384	\\
0.0436	0.0329	\\
0.0481	0.0267	\\
0.0513	0.0209	\\
0.0517	0.0199	\\
0.0545	0.0126	\\
0.0561	0.0049	\\
0.0568	-0.003	\\
0.0568	-0.0038	\\
0.0563	-0.0109	\\
0.0547	-0.0188	\\
0.052	-0.0265	\\
0.0508	-0.0292	\\
0.0483	-0.0338	\\
0.0435	-0.0406	\\
0.0378	-0.0467	\\
0.034	-0.05	\\
0.0313	-0.0521	\\
0.024	-0.0565	\\
0.016	-0.06	\\
0.0089	-0.062	\\
0.0076	-0.0623	\\
-0.0012	-0.0635	\\
-0.0101	-0.0634	\\
-0.0191	-0.0621	\\
-0.0199	-0.062	\\
-0.0279	-0.0596	\\
-0.0364	-0.0559	\\
-0.0444	-0.0509	\\
-0.0469	-0.049	\\
-0.0518	-0.0448	\\
-0.0583	-0.0377	\\
-0.0639	-0.0296	\\
-0.0684	-0.0208	\\
-0.0717	-0.0112	\\
-0.0737	-0.0011	\\
-0.0742	0.0093	\\
-0.0736	0.0178	\\
-0.0734	0.0199	\\
-0.0711	0.0304	\\
-0.0673	0.0407	\\
-0.062	0.0506	\\
-0.0554	0.0598	\\
-0.0474	0.0682	\\
-0.0381	0.0755	\\
-0.0344	0.0779	\\
-0.0277	0.0817	\\
-0.0164	0.0865	\\
-0.0043	0.0898	\\
0.0085	0.0915	\\
0.0216	0.0913	\\
0.0349	0.0894	\\
0.0368	0.089	\\
0.0482	0.0857	\\
0.0612	0.08	\\
0.0736	0.0724	\\
0.0852	0.063	\\
0.0956	0.0518	\\
0.0976	0.0493	\\
0.1048	0.039	\\
0.1124	0.0246	\\
0.1182	0.0089	\\
0.1219	-0.008	\\
0.1231	-0.0198	\\
0.1233	-0.0259	\\
0.1224	-0.0445	\\
0.1189	-0.0634	\\
0.1127	-0.0824	\\
0.1078	-0.0933	\\
0.1036	-0.1011	\\
0.0917	-0.1193	\\
0.0768	-0.1364	\\
0.0591	-0.152	\\
0.0587	-0.1524	\\
0.0383	-0.166	\\
0.0148	-0.1776	\\
-0.0113	-0.1865	\\
-0.0126	-0.1869	\\
-0.0402	-0.1922	\\
-0.0712	-0.1942	\\
-0.0943	-0.1931	\\
-0.1043	-0.192	\\
-0.1392	-0.1849	\\
-0.1753	-0.1725	\\
-0.1769	-0.1719	\\
-0.2124	-0.1539	\\
-0.2496	-0.1288	\\
-0.2531	-0.126	\\
-0.2863	-0.0958	\\
-0.3183	-0.0592	\\
-0.3217	-0.0547	\\
-0.3544	-0.0037	\\
-0.3692	0.0247	\\
-0.3831	0.0579	\\
-0.4038	0.1222	\\
-0.406	0.1317	\\
-0.4199	0.22	\\
-0.4209	0.2304	\\
-0.4208	0.3252	\\
-0.4195	0.3466	\\
-0.4027	0.4502	\\
-0.3988	0.4682	\\
-0.3579	0.593	\\
-0.3552	0.5989	\\
-0.2959	0.7181	\\
-0.2591	0.7746	\\
-0.2124	0.8403	\\
-0.107	0.9558	\\
-0.084	0.976	\\
0.0194	1.0605	\\
0.165	1.1498	\\
0.2318	1.18	\\
0.3265	1.2195	\\
0.4991	1.2661	\\
0.6771	1.2877	\\
0.7767	1.2866	\\
0.8543	1.284	\\
1.0247	1.2567	\\
1.1832	1.2092	\\
1.3265	1.1457	\\
1.4525	1.0711	\\
1.5531	0.9956	\\
1.5607	0.9897	\\
1.6518	0.9055	\\
1.7272	0.8219	\\
1.7887	0.7409	\\
1.8384	0.6644	\\
1.878	0.5932	\\
1.9095	0.5278	\\
1.9344	0.4684	\\
1.9539	0.4148	\\
1.9691	0.3667	\\
1.981	0.3237	\\
1.9903	0.2855	\\
1.9975	0.2516	\\
2.0031	0.2216	\\
2.0074	0.1951	\\
2.0108	0.1717	\\
2.0134	0.151	\\
2.0154	0.1328	\\
2.0169	0.1168	\\
2.0181	0.1027	\\
2.0191	0.0903	\\
2.0198	0.0794	\\
2.0203	0.0698	\\
2.0208	0.0613	\\
2.0211	0.0539	\\
2.0222	0.0054	\\
2.0222	0.0001	\\
2.0222	0	\\
2.0222	0	\\
2.0222	0	\\
2.0222	0	\\
2.0222	0	\\
2.0222	0	\\
2.0222	0	\\
2.0222	0	\\
2.0222	0	\\
2.0222	0	\\
2.0222	-0.0001	\\
2.0222	-0.0054	\\
2.0211	-0.0539	\\
2.0208	-0.0613	\\
2.0203	-0.0698	\\
2.0198	-0.0794	\\
2.0191	-0.0903	\\
2.0181	-0.1027	\\
2.0169	-0.1168	\\
2.0154	-0.1328	\\
2.0134	-0.151	\\
2.0108	-0.1717	\\
2.0074	-0.1951	\\
2.0031	-0.2216	\\
1.9975	-0.2516	\\
1.9903	-0.2855	\\
1.981	-0.3237	\\
1.9691	-0.3667	\\
1.9539	-0.4148	\\
1.9344	-0.4684	\\
1.9095	-0.5278	\\
1.878	-0.5932	\\
1.8384	-0.6644	\\
1.7887	-0.7409	\\
1.7272	-0.8219	\\
1.6518	-0.9055	\\
1.5607	-0.9897	\\
1.5531	-0.9956	\\
1.4525	-1.0711	\\
1.3265	-1.1457	\\
1.1832	-1.2092	\\
1.0247	-1.2567	\\
0.8543	-1.284	\\
0.7767	-1.2866	\\
0.6771	-1.2877	\\
0.4991	-1.2661	\\
0.3265	-1.2195	\\
0.2318	-1.18	\\
0.165	-1.1498	\\
0.0194	-1.0605	\\
-0.084	-0.976	\\
-0.107	-0.9558	\\
-0.2124	-0.8403	\\
-0.2591	-0.7746	\\
-0.2959	-0.7181	\\
-0.3552	-0.5989	\\
-0.3579	-0.593	\\
-0.3988	-0.4682	\\
-0.4027	-0.4502	\\
-0.4195	-0.3466	\\
-0.4208	-0.3252	\\
-0.4209	-0.2304	\\
-0.4199	-0.22	\\
-0.406	-0.1317	\\
-0.4038	-0.1222	\\
-0.3831	-0.0579	\\
-0.3692	-0.0247	\\
-0.3544	0.0037	\\
-0.3217	0.0547	\\
-0.3183	0.0592	\\
-0.2863	0.0958	\\
-0.2531	0.126	\\
-0.2496	0.1288	\\
-0.2124	0.1539	\\
-0.1769	0.1719	\\
-0.1753	0.1725	\\
-0.1392	0.1849	\\
-0.1043	0.192	\\
-0.0943	0.1931	\\
-0.0712	0.1942	\\
-0.0402	0.1922	\\
-0.0126	0.1869	\\
-0.0113	0.1865	\\
0.0148	0.1776	\\
0.0383	0.166	\\
0.0587	0.1524	\\
0.0591	0.152	\\
0.0768	0.1364	\\
0.0917	0.1193	\\
0.1036	0.1011	\\
0.1078	0.0933	\\
0.1127	0.0824	\\
0.1189	0.0634	\\
0.1224	0.0445	\\
0.1233	0.0259	\\
0.1231	0.0198	\\
0.1219	0.008	\\
0.1182	-0.0089	\\
0.1124	-0.0246	\\
0.1048	-0.039	\\
0.0976	-0.0493	\\
0.0956	-0.0518	\\
0.0852	-0.063	\\
0.0736	-0.0724	\\
0.0612	-0.08	\\
0.0482	-0.0857	\\
0.0368	-0.089	\\
0.0349	-0.0894	\\
0.0216	-0.0913	\\
0.0085	-0.0915	\\
-0.0043	-0.0898	\\
-0.0164	-0.0865	\\
-0.0277	-0.0817	\\
-0.0344	-0.0779	\\
-0.0381	-0.0755	\\
-0.0474	-0.0682	\\
-0.0554	-0.0598	\\
-0.062	-0.0506	\\
-0.0673	-0.0407	\\
-0.0711	-0.0304	\\
-0.0734	-0.0199	\\
-0.0736	-0.0178	\\
-0.0742	-0.0093	\\
-0.0737	0.0011	\\
-0.0717	0.0112	\\
-0.0684	0.0208	\\
-0.0639	0.0296	\\
-0.0583	0.0377	\\
-0.0518	0.0448	\\
-0.0469	0.049	\\
-0.0444	0.0509	\\
-0.0364	0.0559	\\
-0.0279	0.0596	\\
-0.0199	0.062	\\
-0.0191	0.0621	\\
-0.0101	0.0634	\\
-0.0012	0.0635	\\
0.0076	0.0623	\\
0.0089	0.062	\\
0.016	0.06	\\
0.024	0.0565	\\
0.0313	0.0521	\\
0.034	0.05	\\
0.0378	0.0467	\\
0.0435	0.0406	\\
0.0483	0.0338	\\
0.0508	0.0292	\\
0.052	0.0265	\\
0.0547	0.0188	\\
0.0563	0.0109	\\
0.0568	0.0038	\\
0.0568	0.003	\\
0.0561	-0.0049	\\
0.0545	-0.0126	\\
0.0517	-0.0199	\\
0.0513	-0.0209	\\
0.0481	-0.0267	\\
0.0436	-0.0329	\\
0.0383	-0.0384	\\
0.0361	-0.0402	\\
0.0323	-0.043	\\
0.0259	-0.0468	\\
0.019	-0.0496	\\
0.0145	-0.0508	\\
0.0118	-0.0514	\\
0.0045	-0.0522	\\
-0.0028	-0.052	\\
-0.0091	-0.051	\\
-0.01	-0.0508	\\
-0.0169	-0.0487	\\
-0.0234	-0.0456	\\
-0.0293	-0.0417	\\
-0.03	-0.0412	\\
-0.0347	-0.037	\\
-0.0393	-0.0317	\\
-0.0431	-0.0258	\\
-0.0443	-0.0235	\\
-0.0461	-0.0194	\\
-0.0494	-0.0017	\\
-0.0446	0.0199	\\
-0.0311	0.0371	\\
-0.0116	0.0466	\\
0.0099	0.0467	\\
0.0293	0.0374	\\
0.0426	0.0208	\\
0.0473	0	\\
};
\addlegendentry{$w_k = 5$}

\addplot[draw=gray,thin] table[row sep=crcr] {
	-1.5 0 \\
	2.5 0 \\
};
\addplot[draw=gray,thin] table[row sep=crcr] {
	0 -2 \\
	0 2 \\
};

\end{axis}
\end{tikzpicture}

%% file: Figures/PlotTimeVaryingCone.tex
\definecolor{mycolorgray}{HTML}{B5B5B5}%
\definecolor{mycolor1}{HTML}{FFB000}%
\definecolor{mycolor2}{HTML}{648FFF}%
\definecolor{mycolor3}{HTML}{009E73}%
\definecolor{mycolor4}{HTML}{F0E442}%
\definecolor{mycolor5}{HTML}{CC79A7}%
\def \w {4}
\def \fsize {\footnotesize}
\def \tsize {\footnotesize}
\begin{tikzpicture}

\begin{axis}[%
width=.3\columnwidth,
height=.3\columnwidth,
at={(0in,0in)},
scale only axis,
bar shift auto,
xmin=-4.5,
xmax=4.5,
xlabel style={font=\color{white!15!black},yshift=1mm},%,font=\fsize},
xlabel={Real},
xtick={-4,  -2,  0,  2,  4},
ymin=-4.5,
ymax=4.5,
ylabel style={font=\color{white!15!black},align=center,yshift=-5mm},%,font=\fsize},
ylabel={Imaginary},
axis background/.style={fill=white},
%axis x line*=bottom,
%axis y line*=left,
legend style={legend cell align=left, align=left, draw=white!15!black}%, font = \fsiz
]
\addplot[draw=black,ultra thick] table[row sep=crcr] {
  -4.5219         0	\\
-4.5129   -0.2839	\\
-4.4862   -0.5667	\\
-4.4418   -0.8473	\\
-4.3798   -1.1245	\\
-4.3006   -1.3973	\\
-4.2043   -1.6646	\\
-4.0915   -1.9253	\\
-3.9625   -2.1784	\\
-3.8179   -2.4229	\\
-3.6583   -2.6579	\\
-3.4842   -2.8823	\\
-3.2963   -3.0954	\\
-3.0954   -3.2963	\\
-2.8823   -3.4842	\\
-2.6579   -3.6583	\\
-2.4229   -3.8179	\\
-2.1784   -3.9625	\\
-1.9253   -4.0915	\\
-1.6646   -4.2043	\\
-1.3973   -4.3006	\\
-1.1245   -4.3798	\\
-0.8473   -4.4418	\\
-0.5667   -4.4862	\\
-0.2839   -4.5129	\\
0.0000   -4.5219	\\
0.2839   -4.5129	\\
0.5667   -4.4862	\\
0.8473   -4.4418	\\
1.1245   -4.3798	\\
1.3973   -4.3006	\\
1.6646   -4.2043	\\
1.9253   -4.0915	\\
2.1784   -3.9625	\\
2.4229   -3.8179	\\
2.6579   -3.6583	\\
2.8823   -3.4842	\\
3.0954   -3.2963	\\
3.2963   -3.0954	\\
3.4842   -2.8823	\\
3.6583   -2.6579	\\
3.8179   -2.4229	\\
3.9625   -2.1784	\\
4.0915   -1.9253	\\
4.2043   -1.6646	\\
4.3006   -1.3973	\\
4.3798   -1.1245	\\
4.4418   -0.8473	\\
4.4862   -0.5667	\\
4.5129   -0.2839	\\
4.5219   -0.0000	\\
4.5129    0.2839	\\
4.4862    0.5667	\\
4.4418    0.8473	\\
4.3798    1.1245	\\
4.3006    1.3973	\\
4.2043    1.6646	\\
4.0915    1.9253	\\
3.9625    2.1784	\\
3.8179    2.4229	\\
3.6583    2.6579	\\
3.4842    2.8823	\\
3.2963    3.0954	\\
3.0954    3.2963	\\
2.8823    3.4842	\\
2.6579    3.6583	\\
2.4229    3.8179	\\
2.1784    3.9625	\\
1.9253    4.0915	\\
1.6646    4.2043	\\
1.3973    4.3006	\\
1.1245    4.3798	\\
0.8473    4.4418	\\
0.5667    4.4862	\\
0.2839    4.5129	\\
0.0000    4.5219	\\
-0.2839    4.5129	\\
-0.5667    4.4862	\\
-0.8473    4.4418	\\
-1.1245    4.3798	\\
-1.3973    4.3006	\\
-1.6646    4.2043	\\
-1.9253    4.0915	\\
-2.1784    3.9625	\\
-2.4229    3.8179	\\
-2.6579    3.6583	\\
-2.8823    3.4842	\\
-3.0954    3.2963	\\
-3.2963    3.0954	\\
-3.4842    2.8823	\\
-3.6583    2.6579	\\
-3.8179    2.4229	\\
-3.9625    2.1784	\\
-4.0915    1.9253	\\
-4.2043    1.6646	\\
-4.3006    1.3973	\\
-4.3798    1.1245	\\
-4.4418    0.8473	\\
-4.4862    0.5667	\\
-4.5129    0.2839	\\
-4.5219    0.0000	\\
};
\addlegendentry{Gain}

\addplot[draw=mycolor1,ultra thick, loosely dashed] table[row sep=crcr] {
  -4.5219         0	\\
-4.5129   -0.2839	\\
-4.4862   -0.5667	\\
-4.4418   -0.8473	\\
-4.3798   -1.1245	\\
-4.3006   -1.3973	\\
-4.2043   -1.6646	\\
-4.0915   -1.9253	\\
-3.9625   -2.1784	\\
-3.8179   -2.4229	\\
-3.6583   -2.6579	\\
-3.4842   -2.8823	\\
-3.2963   -3.0954	\\
-3.0954   -3.2963	\\
-2.8823   -3.4842	\\
-2.6579   -3.6583	\\
-2.4229   -3.8179	\\
-2.1784   -3.9625	\\
-1.9253   -4.0915	\\
-1.6646   -4.2043	\\
-1.3973   -4.3006	\\
-1.1245   -4.3798	\\
-0.8473   -4.4418	\\
-0.5667   -4.4862	\\
-0.2839   -4.5129	\\
0.0000   -4.5219	\\
0.2839   -4.5129	\\
0.5667   -4.4862	\\
0.8473   -4.4418	\\
1.1245   -4.3798	\\
1.3973   -4.3006	\\
1.6646   -4.2043	\\
1.9253   -4.0915	\\
2.1784   -3.9625	\\
2.4229   -3.8179	\\
2.6579   -3.6583	\\
2.8823   -3.4842	\\
3.0954   -3.2963	\\
3.2963   -3.0954	\\
3.4842   -2.8823	\\
3.6583   -2.6579	\\
3.8179   -2.4229	\\
3.9625   -2.1784	\\
4.0915   -1.9253	\\
4.2043   -1.6646	\\
4.3006   -1.3973	\\
4.3798   -1.1245	\\
4.4418   -0.8473	\\
4.4862   -0.5667	\\
4.5129   -0.2839	\\
4.5219   -0.0000	\\
4.5129    0.2839	\\
4.4862    0.5667	\\
4.4418    0.8473	\\
4.3798    1.1245	\\
4.3006    1.3973	\\
4.2043    1.6646	\\
4.0915    1.9253	\\
3.9625    2.1784	\\
3.8179    2.4229	\\
3.6583    2.6579	\\
3.4842    2.8823	\\
3.2963    3.0954	\\
3.0954    3.2963	\\
2.8823    3.4842	\\
2.6579    3.6583	\\
2.4229    3.8179	\\
2.1784    3.9625	\\
1.9253    4.0915	\\
1.6646    4.2043	\\
1.3973    4.3006	\\
1.1245    4.3798	\\
0.8473    4.4418	\\
0.5667    4.4862	\\
0.2839    4.5129	\\
0.0000    4.5219	\\
-0.2839    4.5129	\\
-0.5667    4.4862	\\
-0.8473    4.4418	\\
-1.1245    4.3798	\\
-1.3973    4.3006	\\
-1.6646    4.2043	\\
-1.9253    4.0915	\\
-2.1784    3.9625	\\
-2.4229    3.8179	\\
-2.6579    3.6583	\\
-2.8823    3.4842	\\
-3.0954    3.2963	\\
-3.2963    3.0954	\\
-3.4842    2.8823	\\
-3.6583    2.6579	\\
-3.8179    2.4229	\\
-3.9625    2.1784	\\
-4.0915    1.9253	\\
-4.2043    1.6646	\\
-4.3006    1.3973	\\
-4.3798    1.1245	\\
-4.4418    0.8473	\\
-4.4862    0.5667	\\
-4.5129    0.2839	\\
-4.5219    0.0000	\\
};
\addlegendentry{min $r$}

\addplot[draw=mycolor2,ultra thick, densely dotted] table[row sep=crcr] {
      -4.09     -10		\\
-4.09          10		\\
};
\addlegendentry{max $a$}

\addplot[draw=gray,thin] table[row sep=crcr] {
	-4.5 0 \\
	4.5 0 \\
};
\addplot[draw=gray,thin] table[row sep=crcr] {
	0 -4.5 \\
	0 4.5 \\
};

\end{axis}
\end{tikzpicture}

%% file: Figures/PlotDistribution1Cones.tex
\definecolor{mycolorgray}{HTML}{B5B5B5}%
\definecolor{mycolor1}{HTML}{FFB000}%
\definecolor{mycolor2}{HTML}{648FFF}%
\definecolor{mycolor3}{HTML}{009E73}%
\definecolor{mycolor4}{HTML}{F0E442}%
\definecolor{mycolor5}{HTML}{CC79A7}%
\def \w {4}
\def \fsize {\footnotesize}
\def \tsize {\footnotesize}
\begin{tikzpicture}

\begin{axis}[%
width=.3\columnwidth,
height=.3\columnwidth,
at={(0in,0in)},
scale only axis,
bar shift auto,
xmin=-4.5,
xmax=4.5,
xlabel style={font=\color{white!15!black},yshift=1mm},%,font=\fsize},
xlabel={Real},
xtick={-4,  -2,  0,  2,  4},
ymin=-4.5,
ymax=4.5,
ylabel style={font=\color{white!15!black},align=center,yshift=-5mm},%,font=\fsize},
ylabel={Imaginary},
axis background/.style={fill=white},
%axis x line*=bottom,
%axis y line*=left,
legend style={legend cell align=left, align=left, draw=white!15!black}%, font = \fsiz
]
\addplot[draw=black,ultra thick] table[row sep=crcr] {
         -2.06             0		\\
-2.06         -0.13		\\
-2.05         -0.26		\\
-2.03         -0.39		\\
-2.00         -0.51		\\
-1.96         -0.64		\\
-1.92         -0.76		\\
-1.87         -0.88		\\
-1.81         -0.99		\\
-1.74         -1.10		\\
-1.67         -1.21		\\
-1.59         -1.31		\\
-1.50         -1.41		\\
-1.41         -1.50		\\
-1.31         -1.59		\\
-1.21         -1.67		\\
-1.10         -1.74		\\
-0.99         -1.81		\\
-0.88         -1.87		\\
-0.76         -1.92		\\
-0.64         -1.96		\\
-0.51         -2.00		\\
-0.39         -2.03		\\
-0.26         -2.05		\\
-0.13         -2.06		\\
0.00         -2.06		\\
0.13         -2.06		\\
0.26         -2.05		\\
0.39         -2.03		\\
0.51         -2.00		\\
0.64         -1.96		\\
0.76         -1.92		\\
0.88         -1.87		\\
0.99         -1.81		\\
1.10         -1.74		\\
1.21         -1.67		\\
1.31         -1.59		\\
1.41         -1.50		\\
1.50         -1.41		\\
1.59         -1.31		\\
1.67         -1.21		\\
1.74         -1.10		\\
1.81         -0.99		\\
1.87         -0.88		\\
1.92         -0.76		\\
1.96         -0.64		\\
2.00         -0.51		\\
2.03         -0.39		\\
2.05         -0.26		\\
2.06         -0.13		\\
2.06         -0.00		\\
2.06          0.13		\\
2.05          0.26		\\
2.03          0.39		\\
2.00          0.51		\\
1.96          0.64		\\
1.92          0.76		\\
1.87          0.88		\\
1.81          0.99		\\
1.74          1.10		\\
1.67          1.21		\\
1.59          1.31		\\
1.50          1.41		\\
1.41          1.50		\\
1.31          1.59		\\
1.21          1.67		\\
1.10          1.74		\\
0.99          1.81		\\
0.88          1.87		\\
0.76          1.92		\\
0.64          1.96		\\
0.51          2.00		\\
0.39          2.03		\\
0.26          2.05		\\
0.13          2.06		\\
0.00          2.06		\\
-0.13          2.06		\\
-0.26          2.05		\\
-0.39          2.03		\\
-0.51          2.00		\\
-0.64          1.96		\\
-0.76          1.92		\\
-0.88          1.87		\\
-0.99          1.81		\\
-1.10          1.74		\\
-1.21          1.67		\\
-1.31          1.59		\\
-1.41          1.50		\\
-1.50          1.41		\\
-1.59          1.31		\\
-1.67          1.21		\\
-1.74          1.10		\\
-1.81          0.99		\\
-1.87          0.88		\\
-1.92          0.76		\\
-1.96          0.64		\\
-2.00          0.51		\\
-2.03          0.39		\\
-2.05          0.26		\\
-2.06          0.13		\\
-2.06          0.00		\\
};

\addplot[draw=mycolor1,ultra thick, loosely dashed] table[row sep=crcr] {
    -2.06             0		\\
-2.06         -0.13		\\
-2.05         -0.26		\\
-2.03         -0.39		\\
-2.00         -0.51		\\
-1.96         -0.64		\\
-1.92         -0.76		\\
-1.87         -0.88		\\
-1.81         -0.99		\\
-1.74         -1.10		\\
-1.67         -1.21		\\
-1.59         -1.31		\\
-1.51         -1.41		\\
-1.41         -1.50		\\
-1.32         -1.59		\\
-1.21         -1.67		\\
-1.11         -1.74		\\
-1.00         -1.81		\\
-0.88         -1.87		\\
-0.76         -1.92		\\
-0.64         -1.96		\\
-0.51         -2.00		\\
-0.39         -2.03		\\
-0.26         -2.05		\\
-0.13         -2.06		\\
-0.00         -2.06		\\
0.13         -2.06		\\
0.26         -2.05		\\
0.38         -2.03		\\
0.51         -2.00		\\
0.64         -1.96		\\
0.76         -1.92		\\
0.88         -1.87		\\
0.99         -1.81		\\
1.10         -1.74		\\
1.21         -1.67		\\
1.31         -1.59		\\
1.41         -1.50		\\
1.50         -1.41		\\
1.59         -1.31		\\
1.67         -1.21		\\
1.74         -1.10		\\
1.80         -0.99		\\
1.86         -0.88		\\
1.92         -0.76		\\
1.96         -0.64		\\
2.00         -0.51		\\
2.02         -0.39		\\
2.04         -0.26		\\
2.06         -0.13		\\
2.06         -0.00		\\
2.06          0.13		\\
2.04          0.26		\\
2.02          0.39		\\
2.00          0.51		\\
1.96          0.64		\\
1.92          0.76		\\
1.86          0.88		\\
1.80          0.99		\\
1.74          1.10		\\
1.67          1.21		\\
1.59          1.31		\\
1.50          1.41		\\
1.41          1.50		\\
1.31          1.59		\\
1.21          1.67		\\
1.10          1.74		\\
0.99          1.81		\\
0.88          1.87		\\
0.76          1.92		\\
0.64          1.96		\\
0.51          2.00		\\
0.38          2.03		\\
0.26          2.05		\\
0.13          2.06		\\
-0.00          2.06		\\
-0.13          2.06		\\
-0.26          2.05		\\
-0.39          2.03		\\
-0.51          2.00		\\
-0.64          1.96		\\
-0.76          1.92		\\
-0.88          1.87		\\
-1.00          1.81		\\
-1.11          1.74		\\
-1.21          1.67		\\
-1.32          1.59		\\
-1.41          1.50		\\
-1.51          1.41		\\
-1.59          1.31		\\
-1.67          1.21		\\
-1.74          1.10		\\
-1.81          0.99		\\
-1.87          0.88		\\
-1.92          0.76		\\
-1.96          0.64		\\
-2.00          0.51		\\
-2.03          0.39		\\
-2.05          0.26		\\
-2.06          0.13		\\
-2.06          0.00		\\
};

\addplot[draw=mycolor2,ultra thick, densely dotted] table[row sep=crcr] {
      -1.96    -10		\\
-1.96          10		\\
};

\addplot[draw=gray,thin] table[row sep=crcr] {
	-4.5 0 \\
	4.5 0 \\
};
\addplot[draw=gray,thin] table[row sep=crcr] {
	0 -4.5 \\
	0 4.5 \\
};

\end{axis}
\end{tikzpicture}

%% file: Figures/PlotDistribution2Cones.tex
\definecolor{mycolorgray}{HTML}{B5B5B5}%
\definecolor{mycolor1}{HTML}{FFB000}%
\definecolor{mycolor2}{HTML}{648FFF}%
\definecolor{mycolor3}{HTML}{009E73}%
\definecolor{mycolor4}{HTML}{F0E442}%
\definecolor{mycolor5}{HTML}{CC79A7}%
\def \w {4}
\def \fsize {\footnotesize}
\def \tsize {\footnotesize}
\begin{tikzpicture}

\begin{axis}[%
width=.3\columnwidth,
height=.3\columnwidth,
at={(0in,0in)},
scale only axis,
bar shift auto,
xmin=-4.5,
xmax=4.5,
xlabel style={font=\color{white!15!black},yshift=1mm},%,font=\fsize},
xlabel={Real},
xtick={-4,  -2,  0,  2,  4},
ymin=-4.5,
ymax=4.5,
ylabel style={font=\color{white!15!black},align=center,yshift=-5mm},%,font=\fsize},
ylabel={Imaginary},
axis background/.style={fill=white},
%axis x line*=bottom,
%axis y line*=left,
legend style={legend cell align=left, align=left, draw=white!15!black}%, font = \fsiz
]
\addplot[draw=black,ultra thick] table[row sep=crcr] {
        -2.76             0		\\
    -2.76         -0.17		\\
    -2.74         -0.35		\\
    -2.71         -0.52		\\
    -2.67         -0.69		\\
    -2.63         -0.85		\\
    -2.57         -1.02		\\
    -2.50         -1.18		\\
    -2.42         -1.33		\\
    -2.33         -1.48		\\
    -2.23         -1.62		\\
    -2.13         -1.76		\\
    -2.01         -1.89		\\
    -1.89         -2.01		\\
    -1.76         -2.13		\\
    -1.62         -2.23		\\
    -1.48         -2.33		\\
    -1.33         -2.42		\\
    -1.18         -2.50		\\
    -1.02         -2.57		\\
    -0.85         -2.63		\\
    -0.69         -2.67		\\
    -0.52         -2.71		\\
    -0.35         -2.74		\\
    -0.17         -2.76		\\
    0.00         -2.76		\\
    0.17         -2.76		\\
    0.35         -2.74		\\
    0.52         -2.71		\\
    0.69         -2.67		\\
    0.85         -2.63		\\
    1.02         -2.57		\\
    1.18         -2.50		\\
    1.33         -2.42		\\
    1.48         -2.33		\\
    1.62         -2.23		\\
    1.76         -2.13		\\
    1.89         -2.01		\\
    2.01         -1.89		\\
    2.13         -1.76		\\
    2.23         -1.62		\\
    2.33         -1.48		\\
    2.42         -1.33		\\
    2.50         -1.18		\\
    2.57         -1.02		\\
    2.63         -0.85		\\
    2.67         -0.69		\\
    2.71         -0.52		\\
    2.74         -0.35		\\
    2.76         -0.17		\\
    2.76         -0.00		\\
    2.76          0.17		\\
    2.74          0.35		\\
    2.71          0.52		\\
    2.67          0.69		\\
    2.63          0.85		\\
    2.57          1.02		\\
    2.50          1.18		\\
    2.42          1.33		\\
    2.33          1.48		\\
    2.23          1.62		\\
    2.13          1.76		\\
    2.01          1.89		\\
    1.89          2.01		\\
    1.76          2.13		\\
    1.62          2.23		\\
    1.48          2.33		\\
    1.33          2.42		\\
    1.18          2.50		\\
    1.02          2.57		\\
    0.85          2.63		\\
    0.69          2.67		\\
    0.52          2.71		\\
    0.35          2.74		\\
    0.17          2.76		\\
    0.00          2.76		\\
    -0.17          2.76		\\
    -0.35          2.74		\\
    -0.52          2.71		\\
    -0.69          2.67		\\
    -0.85          2.63		\\
    -1.02          2.57		\\
    -1.18          2.50		\\
    -1.33          2.42		\\
    -1.48          2.33		\\
    -1.62          2.23		\\
    -1.76          2.13		\\
    -1.89          2.01		\\
    -2.01          1.89		\\
    -2.13          1.76		\\
    -2.23          1.62		\\
    -2.33          1.48		\\
    -2.42          1.33		\\
    -2.50          1.18		\\
    -2.57          1.02		\\
    -2.63          0.85		\\
    -2.67          0.69		\\
    -2.71          0.52		\\
    -2.74          0.35		\\
    -2.76          0.17		\\
    -2.76          0.00		\\
};

\addplot[draw=mycolor1,ultra thick, loosely dashed] table[row sep=crcr] {
  -2.71             0		\\
-2.70         -0.17		\\
-2.69         -0.35		\\
-2.66         -0.52		\\
-2.62         -0.69		\\
-2.57         -0.85		\\
-2.51         -1.02		\\
-2.44         -1.17		\\
-2.37         -1.33		\\
-2.28         -1.48		\\
-2.18         -1.62		\\
-2.07         -1.76		\\
-1.96         -1.89		\\
-1.84         -2.01		\\
-1.71         -2.13		\\
-1.57         -2.23		\\
-1.43         -2.33		\\
-1.28         -2.42		\\
-1.12         -2.50		\\
-0.96         -2.57		\\
-0.80         -2.62		\\
-0.63         -2.67		\\
-0.46         -2.71		\\
-0.29         -2.74		\\
-0.12         -2.75		\\
0.05         -2.76		\\
0.23         -2.75		\\
0.40         -2.74		\\
0.57         -2.71		\\
0.74         -2.67		\\
0.91         -2.62		\\
1.07         -2.57		\\
1.23         -2.50		\\
1.38         -2.42		\\
1.53         -2.33		\\
1.67         -2.23		\\
1.81         -2.13		\\
1.94         -2.01		\\
2.06         -1.89		\\
2.18         -1.76		\\
2.28         -1.62		\\
2.38         -1.48		\\
2.47         -1.33		\\
2.55         -1.17		\\
2.62         -1.02		\\
2.68         -0.85		\\
2.72         -0.69		\\
2.76         -0.52		\\
2.79         -0.35		\\
2.81         -0.17		\\
2.81         -0.00		\\
2.81          0.17		\\
2.79          0.35		\\
2.76          0.52		\\
2.72          0.69		\\
2.68          0.85		\\
2.62          1.02		\\
2.55          1.17		\\
2.47          1.33		\\
2.38          1.48		\\
2.28          1.62		\\
2.18          1.76		\\
2.06          1.89		\\
1.94          2.01		\\
1.81          2.13		\\
1.67          2.23		\\
1.53          2.33		\\
1.38          2.42		\\
1.23          2.50		\\
1.07          2.57		\\
0.91          2.62		\\
0.74          2.67		\\
0.57          2.71		\\
0.40          2.74		\\
0.23          2.75		\\
0.05          2.76		\\
-0.12          2.75		\\
-0.29          2.74		\\
-0.46          2.71		\\
-0.63          2.67		\\
-0.80          2.62		\\
-0.96          2.57		\\
-1.12          2.50		\\
-1.28          2.42		\\
-1.43          2.33		\\
-1.57          2.23		\\
-1.71          2.13		\\
-1.84          2.01		\\
-1.96          1.89		\\
-2.07          1.76		\\
-2.18          1.62		\\
-2.28          1.48		\\
-2.37          1.33		\\
-2.44          1.17		\\
-2.51          1.02		\\
-2.57          0.85		\\
-2.62          0.69		\\
-2.66          0.52		\\
-2.69          0.35		\\
-2.70          0.17		\\
-2.71          0.00		\\
};

\addplot[draw=mycolor2,ultra thick, densely dotted] table[row sep=crcr] {
-1.95             -10 \\
-1.95          10.00 \\
};

\addplot[draw=gray,thin] table[row sep=crcr] {
	-4.5 0 \\
	4.5 0 \\
};
\addplot[draw=gray,thin] table[row sep=crcr] {
	0 -4.5 \\
	0 4.5 \\
};

\end{axis}
\end{tikzpicture}

%% file: Figures/PlotDistribution3Cones.tex
\definecolor{mycolorgray}{HTML}{B5B5B5}%
\definecolor{mycolor1}{HTML}{FFB000}%
\definecolor{mycolor2}{HTML}{648FFF}%
\definecolor{mycolor3}{HTML}{009E73}%
\definecolor{mycolor4}{HTML}{F0E442}%
\definecolor{mycolor5}{HTML}{CC79A7}%
\def \w {4}
\def \fsize {\footnotesize}
\def \tsize {\footnotesize}
\begin{tikzpicture}

\begin{axis}[%
width=.3\columnwidth,
height=.3\columnwidth,
at={(0in,0in)},
scale only axis,
bar shift auto,
xmin=-4.5,
xmax=4.5,
xlabel style={font=\color{white!15!black},yshift=1mm},%,font=\fsize},
xlabel={Real},
xtick={-4,  -2,  0,  2,  4},
ymin=-4.5,
ymax=4.5,
ylabel style={font=\color{white!15!black},align=center,yshift=-5mm},%,font=\fsize},
ylabel={Imaginary},
axis background/.style={fill=white},
%axis x line*=bottom,
%axis y line*=left,
legend style={legend cell align=left, align=left, draw=white!15!black}%, font = \fsiz
]
\addplot[draw=black,ultra thick] table[row sep=crcr] {
        -3.91             0		\\
   -3.90         -0.25		\\
   -3.88         -0.49		\\
   -3.84         -0.73		\\
   -3.79         -0.97		\\
   -3.72         -1.21		\\
   -3.64         -1.44		\\
   -3.54         -1.66		\\
   -3.43         -1.88		\\
   -3.30         -2.10		\\
   -3.16         -2.30		\\
   -3.01         -2.49		\\
   -2.85         -2.68		\\
   -2.68         -2.85		\\
   -2.49         -3.01		\\
   -2.30         -3.16		\\
   -2.10         -3.30		\\
   -1.88         -3.43		\\
   -1.66         -3.54		\\
   -1.44         -3.64		\\
   -1.21         -3.72		\\
   -0.97         -3.79		\\
   -0.73         -3.84		\\
   -0.49         -3.88		\\
   -0.25         -3.90		\\
   0.00         -3.91		\\
   0.25         -3.90		\\
   0.49         -3.88		\\
   0.73         -3.84		\\
   0.97         -3.79		\\
   1.21         -3.72		\\
   1.44         -3.64		\\
   1.66         -3.54		\\
   1.88         -3.43		\\
   2.10         -3.30		\\
   2.30         -3.16		\\
   2.49         -3.01		\\
   2.68         -2.85		\\
   2.85         -2.68		\\
   3.01         -2.49		\\
   3.16         -2.30		\\
   3.30         -2.10		\\
   3.43         -1.88		\\
   3.54         -1.66		\\
   3.64         -1.44		\\
   3.72         -1.21		\\
   3.79         -0.97		\\
   3.84         -0.73		\\
   3.88         -0.49		\\
   3.90         -0.25		\\
   3.91         -0.00		\\
   3.90          0.25		\\
   3.88          0.49		\\
   3.84          0.73		\\
   3.79          0.97		\\
   3.72          1.21		\\
   3.64          1.44		\\
   3.54          1.66		\\
   3.43          1.88		\\
   3.30          2.10		\\
   3.16          2.30		\\
   3.01          2.49		\\
   2.85          2.68		\\
   2.68          2.85		\\
   2.49          3.01		\\
   2.30          3.16		\\
   2.10          3.30		\\
   1.88          3.43		\\
   1.66          3.54		\\
   1.44          3.64		\\
   1.21          3.72		\\
   0.97          3.79		\\
   0.73          3.84		\\
   0.49          3.88		\\
   0.25          3.90		\\
   0.00          3.91		\\
   -0.25          3.90		\\
   -0.49          3.88		\\
   -0.73          3.84		\\
   -0.97          3.79		\\
   -1.21          3.72		\\
   -1.44          3.64		\\
   -1.66          3.54		\\
   -1.88          3.43		\\
   -2.10          3.30		\\
   -2.30          3.16		\\
   -2.49          3.01		\\
   -2.68          2.85		\\
   -2.85          2.68		\\
   -3.01          2.49		\\
   -3.16          2.30		\\
   -3.30          2.10		\\
   -3.43          1.88		\\
   -3.54          1.66		\\
   -3.64          1.44		\\
   -3.72          1.21		\\
   -3.79          0.97		\\
   -3.84          0.73		\\
   -3.88          0.49		\\
   -3.90          0.25		\\
   -3.91          0.00		\\
};

\addplot[draw=mycolor1,ultra thick, loosely dashed] table[row sep=crcr] {
   -3.91             0		\\
-3.90         -0.25		\\
-3.88         -0.49		\\
-3.84         -0.73		\\
-3.79         -0.97		\\
-3.72         -1.21		\\
-3.64         -1.44		\\
-3.54         -1.66		\\
-3.43         -1.88		\\
-3.30         -2.10		\\
-3.16         -2.30		\\
-3.01         -2.49		\\
-2.85         -2.68		\\
-2.68         -2.85		\\
-2.49         -3.01		\\
-2.30         -3.16		\\
-2.10         -3.30		\\
-1.88         -3.43		\\
-1.66         -3.54		\\
-1.44         -3.64		\\
-1.21         -3.72		\\
-0.97         -3.79		\\
-0.73         -3.84		\\
-0.49         -3.88		\\
-0.25         -3.90		\\
0.00         -3.91		\\
0.25         -3.90		\\
0.49         -3.88		\\
0.73         -3.84		\\
0.97         -3.79		\\
1.21         -3.72		\\
1.44         -3.64		\\
1.66         -3.54		\\
1.88         -3.43		\\
2.10         -3.30		\\
2.30         -3.16		\\
2.49         -3.01		\\
2.68         -2.85		\\
2.85         -2.68		\\
3.01         -2.49		\\
3.16         -2.30		\\
3.30         -2.10		\\
3.43         -1.88		\\
3.54         -1.66		\\
3.64         -1.44		\\
3.72         -1.21		\\
3.79         -0.97		\\
3.84         -0.73		\\
3.88         -0.49		\\
3.90         -0.25		\\
3.91         -0.00		\\
3.90          0.25		\\
3.88          0.49		\\
3.84          0.73		\\
3.79          0.97		\\
3.72          1.21		\\
3.64          1.44		\\
3.54          1.66		\\
3.43          1.88		\\
3.30          2.10		\\
3.16          2.30		\\
3.01          2.49		\\
2.85          2.68		\\
2.68          2.85		\\
2.49          3.01		\\
2.30          3.16		\\
2.10          3.30		\\
1.88          3.43		\\
1.66          3.54		\\
1.44          3.64		\\
1.21          3.72		\\
0.97          3.79		\\
0.73          3.84		\\
0.49          3.88		\\
0.25          3.90		\\
0.00          3.91		\\
-0.25          3.90		\\
-0.49          3.88		\\
-0.73          3.84		\\
-0.97          3.79		\\
-1.21          3.72		\\
-1.44          3.64		\\
-1.66          3.54		\\
-1.88          3.43		\\
-2.10          3.30		\\
-2.30          3.16		\\
-2.49          3.01		\\
-2.68          2.85		\\
-2.85          2.68		\\
-3.01          2.49		\\
-3.16          2.30		\\
-3.30          2.10		\\
-3.43          1.88		\\
-3.54          1.66		\\
-3.64          1.44		\\
-3.72          1.21		\\
-3.79          0.97		\\
-3.84          0.73		\\
-3.88          0.49		\\
-3.90          0.25		\\
-3.91          0.00		\\
};

\addplot[draw=mycolor2,ultra thick, densely dotted] table[row sep=crcr] {
-2.78             -10 \\
-2.78          10.00 \\
};

\addplot[draw=gray,thin] table[row sep=crcr] {
	-4.5 0 \\
	4.5 0 \\
};
\addplot[draw=gray,thin] table[row sep=crcr] {
	0 -4.5 \\
	0 4.5 \\
};

\end{axis}
\end{tikzpicture}

%% file: Figures/PlotDistribution4Cones.tex
%le was created by matlab2tikz.
\definecolor{mycolorgray}{HTML}{B5B5B5}%
\definecolor{mycolor1}{HTML}{FFB000}%
\definecolor{mycolor2}{HTML}{648FFF}%
\definecolor{mycolor3}{HTML}{009E73}%
\definecolor{mycolor4}{HTML}{F0E442}%
\definecolor{mycolor5}{HTML}{CC79A7}%
\def \w {4}
\def \fsize {\footnotesize}
\def \tsize {\footnotesize}
\begin{tikzpicture}

\begin{axis}[%
width=.3\columnwidth,
height=.3\columnwidth,
at={(0in,0in)},
scale only axis,
bar shift auto,
xmin=-4.5,
xmax=4.5,
xlabel style={font=\color{white!15!black},yshift=1mm},%,font=\fsize},
xlabel={Real},
xtick={-4,  -2,  0,  2,  4},
ymin=-4.5,
ymax=4.5,
ylabel style={font=\color{white!15!black},align=center,yshift=-5mm},%,font=\fsize},
ylabel={Imaginary},
axis background/.style={fill=white},
%axis x line*=bottom,
%axis y line*=left,
legend style={legend cell align=left, align=left, draw=white!15!black}%, font = \fsiz
]
\addplot[draw=black,ultra thick] table[row sep=crcr] {
         -3.63             0		\\
    -3.62         -0.23		\\
    -3.60         -0.45		\\
    -3.56         -0.68		\\
    -3.51         -0.90		\\
    -3.45         -1.12		\\
    -3.37         -1.34		\\
    -3.28         -1.54		\\
    -3.18         -1.75		\\
    -3.06         -1.94		\\
    -2.94         -2.13		\\
    -2.80         -2.31		\\
    -2.64         -2.48		\\
    -2.48         -2.64		\\
    -2.31         -2.80		\\
    -2.13         -2.94		\\
    -1.94         -3.06		\\
    -1.75         -3.18		\\
    -1.54         -3.28		\\
    -1.34         -3.37		\\
    -1.12         -3.45		\\
    -0.90         -3.51		\\
    -0.68         -3.56		\\
    -0.45         -3.60		\\
    -0.23         -3.62		\\
    0.00         -3.63		\\
    0.23         -3.62		\\
    0.45         -3.60		\\
    0.68         -3.56		\\
    0.90         -3.51		\\
    1.12         -3.45		\\
    1.34         -3.37		\\
    1.54         -3.28		\\
    1.75         -3.18		\\
    1.94         -3.06		\\
    2.13         -2.94		\\
    2.31         -2.80		\\
    2.48         -2.64		\\
    2.64         -2.48		\\
    2.80         -2.31		\\
    2.94         -2.13		\\
    3.06         -1.94		\\
    3.18         -1.75		\\
    3.28         -1.54		\\
    3.37         -1.34		\\
    3.45         -1.12		\\
    3.51         -0.90		\\
    3.56         -0.68		\\
    3.60         -0.45		\\
    3.62         -0.23		\\
    3.63         -0.00		\\
    3.62          0.23		\\
    3.60          0.45		\\
    3.56          0.68		\\
    3.51          0.90		\\
    3.45          1.12		\\
    3.37          1.34		\\
    3.28          1.54		\\
    3.18          1.75		\\
    3.06          1.94		\\
    2.94          2.13		\\
    2.80          2.31		\\
    2.64          2.48		\\
    2.48          2.64		\\
    2.31          2.80		\\
    2.13          2.94		\\
    1.94          3.06		\\
    1.75          3.18		\\
    1.54          3.28		\\
    1.34          3.37		\\
    1.12          3.45		\\
    0.90          3.51		\\
    0.68          3.56		\\
    0.45          3.60		\\
    0.23          3.62		\\
    0.00          3.63		\\
    -0.23          3.62		\\
    -0.45          3.60		\\
    -0.68          3.56		\\
    -0.90          3.51		\\
    -1.12          3.45		\\
    -1.34          3.37		\\
    -1.54          3.28		\\
    -1.75          3.18		\\
    -1.94          3.06		\\
    -2.13          2.94		\\
    -2.31          2.80		\\
    -2.48          2.64		\\
    -2.64          2.48		\\
    -2.80          2.31		\\
    -2.94          2.13		\\
    -3.06          1.94		\\
    -3.18          1.75		\\
    -3.28          1.54		\\
    -3.37          1.34		\\
    -3.45          1.12		\\
    -3.51          0.90		\\
    -3.56          0.68		\\
    -3.60          0.45		\\
    -3.62          0.23		\\
    -3.63          0.00		\\
};

\addplot[draw=mycolor1,ultra thick, loosely dashed] table[row sep=crcr] {
       -2.50             0		\\
 -2.49         -0.21		\\
 -2.47         -0.43		\\
 -2.44         -0.64		\\
 -2.39         -0.85		\\
 -2.33         -1.05		\\
 -2.26         -1.25		\\
 -2.18         -1.45		\\
 -2.08         -1.64		\\
 -1.97         -1.83		\\
 -1.85         -2.00		\\
 -1.72         -2.17		\\
 -1.58         -2.33		\\
 -1.43         -2.48		\\
 -1.27         -2.62		\\
 -1.10         -2.76		\\
 -0.92         -2.88		\\
 -0.74         -2.98		\\
 -0.54         -3.08		\\
 -0.35         -3.17		\\
 -0.15         -3.24		\\
 0.06         -3.30		\\
 0.27         -3.35		\\
 0.48         -3.38		\\
 0.69         -3.40		\\
 0.91         -3.41		\\
 1.12         -3.40		\\
 1.33         -3.38		\\
 1.54         -3.35		\\
 1.75         -3.30		\\
 1.96         -3.24		\\
 2.16         -3.17		\\
 2.36         -3.08		\\
 2.55         -2.98		\\
 2.73         -2.88		\\
 2.91         -2.76		\\
 3.08         -2.62		\\
 3.24         -2.48		\\
 3.39         -2.33		\\
 3.53         -2.17		\\
 3.66         -2.00		\\
 3.78         -1.83		\\
 3.89         -1.64		\\
 3.99         -1.45		\\
 4.07         -1.25		\\
 4.15         -1.05		\\
 4.20         -0.85		\\
 4.25         -0.64		\\
 4.28         -0.43		\\
 4.30         -0.21		\\
 4.31         -0.00		\\
 4.30          0.21		\\
 4.28          0.43		\\
 4.25          0.64		\\
 4.20          0.85		\\
 4.15          1.05		\\
 4.07          1.25		\\
 3.99          1.45		\\
 3.89          1.64		\\
 3.78          1.83		\\
 3.66          2.00		\\
 3.53          2.17		\\
 3.39          2.33		\\
 3.24          2.48		\\
 3.08          2.62		\\
 2.91          2.76		\\
 2.73          2.88		\\
 2.55          2.98		\\
 2.36          3.08		\\
 2.16          3.17		\\
 1.96          3.24		\\
 1.75          3.30		\\
 1.54          3.35		\\
 1.33          3.38		\\
 1.12          3.40		\\
 0.91          3.41		\\
 0.69          3.40		\\
 0.48          3.38		\\
 0.27          3.35		\\
 0.06          3.30		\\
 -0.15          3.24		\\
 -0.35          3.17		\\
 -0.54          3.08		\\
 -0.74          2.98		\\
 -0.92          2.88		\\
 -1.10          2.76		\\
 -1.27          2.62		\\
 -1.43          2.48		\\
 -1.58          2.33		\\
 -1.72          2.17		\\
 -1.85          2.00		\\
 -1.97          1.83		\\
 -2.08          1.64		\\
 -2.18          1.45		\\
 -2.26          1.25		\\
 -2.33          1.05		\\
 -2.39          0.85		\\
 -2.44          0.64		\\
 -2.47          0.43		\\
 -2.49          0.21		\\
 -2.50          0.00		\\
};

\addplot[draw=mycolor2,ultra thick, densely dotted] table[row sep=crcr] {
         -0.87            -10 \\
-0.87          10.00 \\
};

\addplot[draw=gray,thin] table[row sep=crcr] {
	-4.5 0 \\
	4.5 0 \\
};
\addplot[draw=gray,thin] table[row sep=crcr] {
	0 -4.5 \\
	0 4.5 \\
};

\end{axis}
\end{tikzpicture}

%% file: Figures/PlotDistribution5Cones.tex
\definecolor{mycolorgray}{HTML}{B5B5B5}%
\definecolor{mycolor1}{HTML}{FFB000}%
\definecolor{mycolor2}{HTML}{648FFF}%
\definecolor{mycolor3}{HTML}{009E73}%
\definecolor{mycolor4}{HTML}{F0E442}%
\definecolor{mycolor5}{HTML}{CC79A7}%
\def \w {4}
\def \fsize {\footnotesize}
\def \tsize {\footnotesize}
\begin{tikzpicture}

\begin{axis}[%
width=.3\columnwidth,
height=.3\columnwidth,
at={(0in,0in)},
scale only axis,
bar shift auto,
xmin=-4.5,
xmax=4.5,
xlabel style={font=\color{white!15!black},yshift=1mm},%,font=\fsize},
xlabel={Real},
xtick={-4,  -2,  0,  2,  4},
ymin=-4.5,
ymax=4.5,
ylabel style={font=\color{white!15!black},align=center,yshift=-5mm},%,font=\fsize},
ylabel={Imaginary},
axis background/.style={fill=white},
%axis x line*=bottom,
%axis y line*=left,
legend style={legend cell align=left, align=left, draw=white!15!black}
]
\addplot[draw=black,ultra thick] table[row sep=crcr] {
         -2.11             0		\\
    -2.11         -0.13		\\
    -2.10         -0.26		\\
    -2.08         -0.40		\\
    -2.05         -0.53		\\
    -2.01         -0.65		\\
    -1.97         -0.78		\\
    -1.91         -0.90		\\
    -1.85         -1.02		\\
    -1.78         -1.13		\\
    -1.71         -1.24		\\
    -1.63         -1.35		\\
    -1.54         -1.45		\\
    -1.45         -1.54		\\
    -1.35         -1.63		\\
    -1.24         -1.71		\\
    -1.13         -1.78		\\
    -1.02         -1.85		\\
    -0.90         -1.91		\\
    -0.78         -1.97		\\
    -0.65         -2.01		\\
    -0.53         -2.05		\\
    -0.40         -2.08		\\
    -0.26         -2.10		\\
    -0.13         -2.11		\\
    0.00         -2.11		\\
    0.13         -2.11		\\
    0.26         -2.10		\\
    0.40         -2.08		\\
    0.53         -2.05		\\
    0.65         -2.01		\\
    0.78         -1.97		\\
    0.90         -1.91		\\
    1.02         -1.85		\\
    1.13         -1.78		\\
    1.24         -1.71		\\
    1.35         -1.63		\\
    1.45         -1.54		\\
    1.54         -1.45		\\
    1.63         -1.35		\\
    1.71         -1.24		\\
    1.78         -1.13		\\
    1.85         -1.02		\\
    1.91         -0.90		\\
    1.97         -0.78		\\
    2.01         -0.65		\\
    2.05         -0.53		\\
    2.08         -0.40		\\
    2.10         -0.26		\\
    2.11         -0.13		\\
    2.11         -0.00		\\
    2.11          0.13		\\
    2.10          0.26		\\
    2.08          0.40		\\
    2.05          0.53		\\
    2.01          0.65		\\
    1.97          0.78		\\
    1.91          0.90		\\
    1.85          1.02		\\
    1.78          1.13		\\
    1.71          1.24		\\
    1.63          1.35		\\
    1.54          1.45		\\
    1.45          1.54		\\
    1.35          1.63		\\
    1.24          1.71		\\
    1.13          1.78		\\
    1.02          1.85		\\
    0.90          1.91		\\
    0.78          1.97		\\
    0.65          2.01		\\
    0.53          2.05		\\
    0.40          2.08		\\
    0.26          2.10		\\
    0.13          2.11		\\
    0.00          2.11		\\
    -0.13          2.11		\\
    -0.26          2.10		\\
    -0.40          2.08		\\
    -0.53          2.05		\\
    -0.65          2.01		\\
    -0.78          1.97		\\
    -0.90          1.91		\\
    -1.02          1.85		\\
    -1.13          1.78		\\
    -1.24          1.71		\\
    -1.35          1.63		\\
    -1.45          1.54		\\
    -1.54          1.45		\\
    -1.63          1.35		\\
    -1.71          1.24		\\
    -1.78          1.13		\\
    -1.85          1.02		\\
    -1.91          0.90		\\
    -1.97          0.78		\\
    -2.01          0.65		\\
    -2.05          0.53		\\
    -2.08          0.40		\\
    -2.10          0.26		\\
    -2.11          0.13		\\
    -2.11          0.00		\\
};

\addplot[draw=mycolor1,ultra thick, loosely dashed] table[row sep=crcr] {
      -0.62             0		\\
-0.61         -0.10		\\
-0.61         -0.19		\\
-0.59         -0.28		\\
-0.57         -0.38		\\
-0.54         -0.47		\\
-0.51         -0.56		\\
-0.47         -0.65		\\
-0.43         -0.73		\\
-0.38         -0.81		\\
-0.33         -0.89		\\
-0.27         -0.97		\\
-0.21         -1.04		\\
-0.14         -1.11		\\
-0.07         -1.17		\\
0.01         -1.23		\\
0.09         -1.28		\\
0.17         -1.33		\\
0.25         -1.37		\\
0.34         -1.41		\\
0.43         -1.44		\\
0.52         -1.47		\\
0.61         -1.49		\\
0.71         -1.50		\\
0.80         -1.51		\\
0.90         -1.52		\\
0.99         -1.51		\\
1.09         -1.50		\\
1.18         -1.49		\\
1.28         -1.47		\\
1.37         -1.44		\\
1.46         -1.41		\\
1.54         -1.37		\\
1.63         -1.33		\\
1.71         -1.28		\\
1.79         -1.23		\\
1.86         -1.17		\\
1.94         -1.11		\\
2.00         -1.04		\\
2.07         -0.97		\\
2.12         -0.89		\\
2.18         -0.81		\\
2.23         -0.73		\\
2.27         -0.65		\\
2.31         -0.56		\\
2.34         -0.47		\\
2.37         -0.38		\\
2.39         -0.28		\\
2.40         -0.19		\\
2.41         -0.10		\\
2.41         -0.00		\\
2.41          0.10		\\
2.40          0.19		\\
2.39          0.28		\\
2.37          0.38		\\
2.34          0.47		\\
2.31          0.56		\\
2.27          0.65		\\
2.23          0.73		\\
2.18          0.81		\\
2.12          0.89		\\
2.07          0.97		\\
2.00          1.04		\\
1.94          1.11		\\
1.86          1.17		\\
1.79          1.23		\\
1.71          1.28		\\
1.63          1.33		\\
1.54          1.37		\\
1.46          1.41		\\
1.37          1.44		\\
1.28          1.47		\\
1.18          1.49		\\
1.09          1.50		\\
0.99          1.51		\\
0.90          1.52		\\
0.80          1.51		\\
0.71          1.50		\\
0.61          1.49		\\
0.52          1.47		\\
0.43          1.44		\\
0.34          1.41		\\
0.25          1.37		\\
0.17          1.33		\\
0.09          1.28		\\
0.01          1.23		\\
-0.07          1.17		\\
-0.14          1.11		\\
-0.21          1.04		\\
-0.27          0.97		\\
-0.33          0.89		\\
-0.38          0.81		\\
-0.43          0.73		\\
-0.47          0.65		\\
-0.51          0.56		\\
-0.54          0.47		\\
-0.57          0.38		\\
-0.59          0.28		\\
-0.61          0.19		\\
-0.61          0.10		\\
-0.62          0.00		\\
};

\addplot[draw=mycolor2,ultra thick, densely dotted] table[row sep=crcr] {
         -0.19            -10 \\
-0.19          10.00 \\
};

\addplot[draw=gray,thin] table[row sep=crcr] {
	-4.5 0 \\
	4.5 0 \\
};
\addplot[draw=gray,thin] table[row sep=crcr] {
	0 -4.5 \\
	0 4.5 \\
};

\end{axis}
\end{tikzpicture}

%% file: Figures/PlotSimulationDelay.tex
\definecolor{mycolorgray}{HTML}{B5B5B5}%
\definecolor{mycolor2}{HTML}{D55E00}%
\definecolor{mycolor1}{HTML}{0072B2}%
\definecolor{mycolor3}{HTML}{009E73}%
\definecolor{mycolor4}{HTML}{F0E442}%
\definecolor{mycolor5}{HTML}{CC79A7}%
\def \w {4}
\def \fsize {\footnotesize}
\def \tsize {\footnotesize}
\begin{tikzpicture}

\begin{axis}[%
width=.8\columnwidth,
height=.27\columnwidth,
at={(0in,0in)},
scale only axis,
bar shift auto,
xmin=0,
xmax=50,
xlabel style={font=\color{white!15!black},yshift=2mm},%,font=\fsize},
xlabel={Time},
ymin=0,
ymax=5,
ylabel style={font=\color{white!15!black},align=center},%,font=\fsize},
ylabel={Delay},
ytick = {0,1,2,3,4,5},
axis background/.style={fill=white},
axis x line*=bottom,
axis y line*=left,
legend style={legend cell align=left, align=left, draw=white!15!black}%, font = \fsiz
]
\addplot[draw= green,thick] table[row sep=crcr] {
0	5	\\
1  5   \\
1	1	\\
2	1	\\
3	1	\\
4	1	\\
5	1	\\
6	1	\\
7	1	\\
8	1	\\
9	1	\\
10 1 \\
10	2	\\
11 2 \\
11	1	\\
12	1	\\
13 1 \\
13	2	\\
14 2 \\
14	1	\\
15	1	\\
16	1	\\
17 1 \\
17	4	\\
18 4 \\
18	2	\\
19 2 \\
19	1	\\
20	1	\\
21	1	\\
22	1	\\
23	1	\\
24	1	\\
25	1	\\
26 1 \\
26	4	\\
27 4 \\
27	1	\\
28	1	\\
29	1	\\
30	1	\\
31	1	\\
32	1	\\
33 1 \\
33	3	\\
34 3 \\
34	1	\\
35 1 \\
35	2	\\
36 2 \\
36	1	\\
37	1	\\
38	1	\\
39	1	\\
40	1	\\
41 1 \\
41	2	\\
42 2 \\
42	1	\\
43	1	\\
44	1	\\
45	1	\\
46	1	\\
47	1	\\
48	1	\\
49	1	\\
};

\end{axis}
\end{tikzpicture}

%% file: Figures/PlotSimulationOpenLoop.tex
\definecolor{mycolorgray}{HTML}{B5B5B5}%
\definecolor{mycolor2}{HTML}{D55E00}%
\definecolor{mycolor1}{HTML}{0072B2}%
\definecolor{mycolor3}{HTML}{009E73}%
\definecolor{mycolor4}{HTML}{F0E442}%
\definecolor{mycolor5}{HTML}{CC79A7}%
\def \w {4}
\def \fsize {\footnotesize}
\def \tsize {\footnotesize}
\begin{tikzpicture}

\begin{axis}[%
width=.8\columnwidth,
height=.27\columnwidth,
at={(0in,0in)},
scale only axis,
bar shift auto,
xmin=0,
xmax=50,
xlabel style={font=\color{white!15!black},yshift=2mm},%,font=\fsize},
xlabel={Time},
ymin=-3,
ymax=3,
ylabel style={font=\color{white!15!black},align=center,yshift=-3mm},%,font=\fsize},
ylabel={Output},
axis background/.style={fill=white},
axis x line*=bottom,
axis y line*=left,
legend style={legend cell align=left, align=left, draw=white!15!black},%, font = \fsiz
legend pos = south east
]
\addplot[draw= black,thick] table[row sep=crcr] {
-5	0	\\
-4	0	\\
-3	0	\\
-2	0	\\
-1	0	\\
0	0	\\
1	0.96	\\
2	1.93	\\
3	2.88	\\
4	2.85	\\
5	2.77	\\
6	2.66	\\
7	2.54	\\
8	2.41	\\
9	2.29	\\
10	2.17	\\
11	2.05	\\
12	1.94	\\
13	1.84	\\
14	1.73	\\
15	1.64	\\
16	1.55	\\
17	1.46	\\
18	1.38	\\
19	1.31	\\
20	1.23	\\
21	1.17	\\
22	1.1	\\
23	1.04	\\
24	0.98	\\
25	0.93	\\
26	0.88	\\
27	0.83	\\
28	0.78	\\
29	0.74	\\
30	0.7	\\
31	0.66	\\
32	0.62	\\
33	0.59	\\
34	0.56	\\
35	0.53	\\
36	0.5	\\
37	0.47	\\
38	0.44	\\
39	0.42	\\
40	0.4	\\
41	0.37	\\
42	0.35	\\
43	0.33	\\
44	0.31	\\
45	0.3	\\
46	0.28	\\
47	0.27	\\
48	0.25	\\
49	0.24	\\
};
\addlegendentry{Plant}

\addplot[draw=red,dashed,thick] table[row sep=crcr] {
-1 -1 \\ 
-1 -2 \\
};
\addlegendentry{Controller}

\addplot[draw=gray,thin] table[row sep=crcr] {
	0 0 \\ 
	50 0 \\	
};

\end{axis}
\end{tikzpicture}

%% file: Figures/PlotSimulationDeterministic.tex
\definecolor{mycolorgray}{HTML}{B5B5B5}%
\definecolor{mycolor2}{HTML}{D55E00}%
\definecolor{mycolor1}{HTML}{0072B2}%
\definecolor{mycolor3}{HTML}{009E73}%
\definecolor{mycolor4}{HTML}{F0E442}%
\definecolor{mycolor5}{HTML}{CC79A7}%
\def \w {4}
\def \fsize {\footnotesize}
\def \tsize {\footnotesize}
\begin{tikzpicture}

\begin{axis}[%
width=.8\columnwidth,
height=.27\columnwidth,
at={(0in,0in)},
scale only axis,
bar shift auto,
xmin=0,
xmax=50,
xlabel style={font=\color{white!15!black},yshift=2mm},%,font=\fsize},
xlabel={Time},
xtick={0,25,50},
ymin=-3,
ymax=3,
ylabel style={font=\color{white!15!black},align=center,yshift=-3mm},%,font=\fsize},
ylabel={Output},
axis background/.style={fill=white},
axis x line*=bottom,
axis y line*=left,
legend style={legend cell align=left, align=left, draw=white!15!black}%, font = \fsiz
]
\addplot[draw= black,thick] table[row sep=crcr] {
-5	0	\\
-4	0	\\
-3	0	\\
-2	0	\\
-1	0	\\
0	0	\\
1	0.96	\\
2	1.72	\\
3	2.44	\\
4	2.16	\\
5	2.03	\\
6	1.9	\\
7	1.76	\\
8	1.63	\\
9	1.5	\\
10	1.38	\\
11	1.26	\\
12	1.16	\\
13	1.06	\\
14	0.97	\\
15	0.89	\\
16	0.82	\\
17	0.75	\\
18	0.68	\\
19	0.63	\\
20	0.57	\\
21	0.52	\\
22	0.48	\\
23	0.44	\\
24	0.4	\\
25	0.37	\\
26	0.34	\\
27	0.31	\\
28	0.28	\\
29	0.26	\\
30	0.24	\\
31	0.22	\\
32	0.2	\\
33	0.18	\\
34	0.17	\\
35	0.15	\\
36	0.14	\\
37	0.13	\\
38	0.12	\\
39	0.11	\\
40	0.1	\\
41	0.09	\\
42	0.08	\\
43	0.08	\\
44	0.07	\\
45	0.06	\\
46	0.06	\\
47	0.05	\\
48	0.05	\\
49	0.05	\\
};

\addplot[draw= red,thick,dashed] table[row sep=crcr] {
-5	0	\\
-4	0	\\
-3	0	\\
-2	0	\\
-1	0	\\
0	-2.2	\\
1	-2.41	\\
2	-2.58	\\
3	-0.54	\\
4	-0.47	\\
5	-0.45	\\
6	-0.42	\\
7	-0.39	\\
8	-0.36	\\
9	-0.33	\\
10	-0.3	\\
11	-0.28	\\
12	-0.26	\\
13	-0.23	\\
14	-0.21	\\
15	-0.2	\\
16	-0.18	\\
17	-0.17	\\
18	-0.15	\\
19	-0.14	\\
20	-0.13	\\
21	-0.12	\\
22	-0.11	\\
23	-0.1	\\
24	-0.09	\\
25	-0.08	\\
26	-0.07	\\
27	-0.07	\\
28	-0.06	\\
29	-0.06	\\
30	-0.05	\\
31	-0.05	\\
32	-0.04	\\
33	-0.04	\\
34	-0.04	\\
35	-0.03	\\
36	-0.03	\\
37	-0.03	\\
38	-0.03	\\
39	-0.02	\\
40	-0.02	\\
41	-0.02	\\
42	-0.02	\\
43	-0.02	\\
44	-0.02	\\
45	-0.01	\\
46	-0.01	\\
47	-0.01	\\
48	-0.01	\\
49	-0.01	\\
};

\addplot[draw=gray,thin] table[row sep=crcr] {
	0 0 \\ 
	50 0 \\	
};

\end{axis}
\end{tikzpicture}

%% file: Figures/PlotSimulationStochastic.tex
\definecolor{mycolorgray}{HTML}{B5B5B5}%
\definecolor{mycolor2}{HTML}{D55E00}%
\definecolor{mycolor1}{HTML}{0072B2}%
\definecolor{mycolor3}{HTML}{009E73}%
\definecolor{mycolor4}{HTML}{F0E442}%
\definecolor{mycolor5}{HTML}{CC79A7}%
\def \w {4}
\def \fsize {\footnotesize}
\def \tsize {\footnotesize}
\begin{tikzpicture}

\begin{axis}[%
width=.8\columnwidth,
height=.27\columnwidth,
at={(0in,0in)},
scale only axis,
bar shift auto,
xmin=0,
xmax=50,
xlabel style={font=\color{white!15!black},yshift=2mm},%,font=\fsize},
xlabel={Time},
ymin=-3,
ymax=3,
ylabel style={font=\color{white!15!black},align=center,yshift=-3mm},%,font=\fsize},
ylabel={Output},
axis background/.style={fill=white},
axis x line*=bottom,
axis y line*=left,
legend style={legend cell align=left, align=left, draw=white!15!black}%, font = \fsiz
]
\addplot[draw= black,thick] table[row sep=crcr] {
-5	0	\\
-4	0	\\
-3	0	\\
-2	0	\\
-1	0	\\
0	0	\\
1	0.96	\\
2	0.8204	\\
3	0.5472	\\
4	-0.6835	\\
5	-0.7852	\\
6	-0.7175	\\
7	-0.6167	\\
8	-0.513	\\
9	-0.4169	\\
10	-0.3329	\\
11	-0.2514	\\
12	-0.1931	\\
13	-0.1476	\\
14	-0.1055	\\
15	-0.0779	\\
16	-0.0575	\\
17	-0.0423	\\
18	-0.0209	\\
19	-0.0107	\\
20	-0.0058	\\
21	-0.0029	\\
22	-0.0012	\\
23	-0.0002	\\
24	0.0003	\\
25	0.0005	\\
26	0.0006	\\
27	0.0008	\\
28	0.0007	\\
29	0.0006	\\
30	0.0005	\\
31	0.0004	\\
32	0.0004	\\
33	0.0003	\\
34	0.0002	\\
35	0.0001	\\
36	0.0001	\\
37	0.0001	\\
38	0.0001	\\
39	0	\\
40	0	\\
41	0	\\
42	0	\\
43	0	\\
44	0	\\
45	0	\\
46	0	\\
47	0	\\
48	0	\\
49	0	\\
};

\addplot[draw= red,thick,dashed] table[row sep=crcr] {
-5	0	\\
-4	0	\\
-3	0	\\
-2	0	\\
-1	0	\\
0	-11.5579	\\
1	-12.6675	\\
2	-12.5062	\\
3	-0.6324	\\
4	0.79	\\
5	0.9076	\\
6	0.8293	\\
7	0.7128	\\
8	0.5929	\\
9	0.4819	\\
10	0.3848	\\
11	0.2906	\\
12	0.2232	\\
13	0.1706	\\
14	0.122	\\
15	0.09	\\
16	0.0664	\\
17	0.0488	\\
18	0.0242	\\
19	0.0124	\\
20	0.0067	\\
21	0.0034	\\
22	0.0014	\\
23	0.0003	\\
24	-0.0003	\\
25	-0.0006	\\
26	-0.0007	\\
27	-0.0009	\\
28	-0.0009	\\
29	-0.0008	\\
30	-0.0006	\\
31	-0.0005	\\
32	-0.0004	\\
33	-0.0003	\\
34	-0.0002	\\
35	-0.0002	\\
36	-0.0001	\\
37	-0.0001	\\
38	-0.0001	\\
39	0	\\
40	0	\\
41	0	\\
42	0	\\
43	0	\\
44	0	\\
45	0	\\
46	0	\\
47	0	\\
48	0	\\
49	0	\\
};

\addplot[draw=gray,thin] table[row sep=crcr] {
	0 0 \\ 
	50 0 \\	
};

\end{axis}
\end{tikzpicture}

%% file: Figures/PlotSimulationUnstable.tex
\definecolor{mycolorgray}{HTML}{B5B5B5}%
\definecolor{mycolor2}{HTML}{D55E00}%
\definecolor{mycolor1}{HTML}{0072B2}%
\definecolor{mycolor3}{HTML}{009E73}%
\definecolor{mycolor4}{HTML}{F0E442}%
\definecolor{mycolor5}{HTML}{CC79A7}%
\def \w {4}
\def \fsize {\footnotesize}
\def \tsize {\footnotesize}
\begin{tikzpicture}

\begin{axis}[%
width=.75\columnwidth,
height=.27\columnwidth,
at={(0in,0in)},
scale only axis,
bar shift auto,
xmin=0,
xmax=50,
xlabel style={font=\color{white!15!black},yshift=2mm},%,font=\fsize},
xlabel={Time},
ymin=-3,
ymax=3,
ylabel style={font=\color{white!15!black},align=center,yshift=-3mm},%,font=\fsize},
ylabel={Output},
ytick={-3,-2,-1,0,1,2,3},
yticklabels={$-3e6$,$-2e6$,$-1e6$,$0$,$1e6$,$2e6$,$3e6$},
axis background/.style={fill=white},
axis x line*=bottom,
axis y line*=left,
legend style={legend cell align=left, align=left, draw=white!15!black}%, font = \fsiz
]
\addplot[draw= black,thick] table[row sep=crcr] {
-5	0	\\
-4	0	\\
-3	0	\\
-2	0	\\
-1	0	\\
0	0	\\
1	0.00000096	\\
2	-0.00001535	\\
3	-0.00003351	\\
4	-0.00002422	\\
5	0.00003455	\\
6	0.00007801	\\
7	0.00001886	\\
8	-0.00011797	\\
9	-0.0001535	\\
10	0.000051	\\
11	0.00026071	\\
12	0.00017679	\\
13	-0.00027807	\\
14	-0.00073955	\\
15	-0.00026362	\\
16	0.00103307	\\
17	0.00151798	\\
18	0.00199661	\\
19	0.00018208	\\
20	-0.00334664	\\
21	-0.00374423	\\
22	0.00206361	\\
23	0.0086898	\\
24	0.00523143	\\
25	-0.0099538	\\
26	-0.0193578	\\
27	-0.02299847	\\
28	0.01079917	\\
29	0.05158468	\\
30	0.03367208	\\
31	-0.0563847	\\
32	-0.11674167	\\
33	-0.01986481	\\
34	-0.07477359	\\
35	-0.03882308	\\
36	-0.00160653	\\
37	0.0682729	\\
38	0.07331242	\\
39	-0.04489557	\\
40	-0.17454714	\\
41	-0.09892555	\\
42	-0.01777342	\\
43	0.15810827	\\
44	0.193696	\\
45	-0.08000094	\\
46	-0.4220773	\\
47	-0.28959942	\\
48	0.44731143	\\
49	0.9658428	\\
};

\addplot[draw= red,thick,dashed] table[row sep=crcr] {
-5	0	\\
-4	0	\\
-3	0	\\
-2	0	\\
-1	0	\\
0	-0.00018	\\
1	-0.00019728	\\
2	0.0000963	\\
3	0.00060326	\\
4	0.00043598	\\
5	-0.00062192	\\
6	-0.00140421	\\
7	-0.00033954	\\
8	0.00212352	\\
9	0.00276302	\\
10	-0.00091808	\\
11	-0.00469287	\\
12	-0.00318214	\\
13	0.00500534	\\
14	0.01331186	\\
15	0.0047451	\\
16	-0.01859533	\\
17	-0.02732362	\\
18	-0.03593906	\\
19	-0.00327745	\\
20	0.06023946	\\
21	0.06739621	\\
22	-0.03714497	\\
23	-0.15641635	\\
24	-0.09416567	\\
25	0.17916845	\\
26	0.34844042	\\
27	0.41397249	\\
28	-0.19438509	\\
29	-0.92852427	\\
30	-0.60609742	\\
31	1.01492465	\\
32	2.10135	\\
33	0.35756656	\\
34	1.34592469	\\
35	0.69881544	\\
36	0.02891756	\\
37	-1.22891219	\\
38	-1.31962359	\\
39	0.80812028	\\
40	3.14184861	\\
41	1.78065988	\\
42	0.31992148	\\
43	-2.84594886	\\
44	-3.48652804	\\
45	1.44001698	\\
46	7.59739136	\\
47	5.2127896	\\
48	-8.0516058	\\
49	-17.38517032	\\
};

\addplot[draw=gray,thin] table[row sep=crcr] {
	0 0 \\ 
	50 0 \\	
};

\end{axis}
\end{tikzpicture}